\documentclass[12pt]{article}
\usepackage[utf8]{inputenc}
\usepackage[english]{babel}
\usepackage{comment}
\usepackage[dvipsnames,svgnames,x11names]{xcolor}
\usepackage{tikz}
\usepackage{bm}
\usetikzlibrary{shapes.multipart, positioning}
\usepackage{multicol}
\usepackage{amsfonts,amsthm,amsmath,amssymb}
\usepackage{caption}
\usepackage{graphicx,float,subcaption}
\usepackage{verbatim}
\usepackage{amsthm}
\usepackage{graphicx}
\usepackage{caption}
\usepackage{booktabs}
\usepackage{hyperref}
\usepackage{threeparttable}
\graphicspath{{./immagini/}}
\hypersetup{hypertex=true,
colorlinks=true,
linkcolor=black,
anchorcolor=black,
citecolor=black}

\newtheorem{theorem}{Theorem}[section]
\newtheorem{lemma}{Lemma}[section]
\newtheorem{proposition}[theorem]{Proposition}

\newtheorem{remark}{Remark}[section]

\newtheorem{corollary}[theorem]{Corollary}

\numberwithin{equation}{section}
\allowdisplaybreaks
\date{}

\begin{document}
\title{The discrete inverse conductivity problem solved by the
weights of an interpretable neural network\thanks{The work of M. Deng is supported by the Hong Kong PhD fellowship scheme, and that of B. Jin is supported by Hong Kong RGC General Research Fund (Projects
14306423 and 14306824).}}
\author{Elena Beretta\thanks{Division of Science, New York University Abu Dhabi, Saadiyat Island, Abu Dhabi, United Arab Emirates (\texttt{eb147@nyu.edu,ag189@nyu.edu})} \and Maolin Deng\thanks{Department of Mathematics, The Chinese University of Hong Kong, Shatin, N.T., Hong Kong (\texttt{mldeng@link.cuhk.edu.hk, b.jin@cuhk.edu.hk})} \and Alberto Gandolfi\footnotemark[2] \and Bangti Jin\footnotemark[3]}
%\author[1]{Elena Beretta}
%\author[2]{Maolin Deng}
%\author[1]{Alberto Gandolfi}
%\author[2]{Bangti Jin}

%\affil[1]{Division of Science, New York University Abu Dhabi, Saadiyat Island, Abu Dhabi, United Arab Emirates (\texttt{eb147@nyu.edu,ag189@nyu.edu})}
%\affil[2]{Department of Mathematics, The Chinese University of Hong Kong, Shatin, N.T., Hong Kong (\texttt{mldeng@link.cuhk.edu.hk, b.jin@cuhk.edu.hk})}

\maketitle

\begin{abstract}
In this work, we develop a novel neural network (NN) approach to solve the discrete inverse conductivity problem of recovering the conductivity profile on network edges from the discrete Dirichlet-to-Neumann map on a square lattice.

The novelty of the approach lies in the fact that the sought-after conductivity is not provided directly as the output of the NN but is instead encoded in the weights of the post-trainig NN in the second layer. Hence the weights of the trained NN acquire a clear physical meaning, which contrasts with most existing neural network approaches, where the weights are typically not interpretable. This work represents a step toward designing NNs with interpretable post-training weights.

Numerically, we observe that the method outperforms the conventional Curtis-Morrow algorithm for both noisy full and partial data.
\vskip 3truemm
\noindent\textbf{Key words}: discrete inverse conductivity problem, interpretable neural network, global minimum, sensitivity analysis
\end{abstract}

\section{Introduction}

Neural Networks (NNs) have been extensively applied across various domains of science, engineering and real-life tasks \cite{bishop1994neural}, achieving remarkable success in solving complex problems that were previously deemed intractable \cite{sejnowski2020unreasonable}. However, a persistent challenge in the application of NNs is the intelligibility of the final status of trained models \cite{GaoGuan:2023}. The weights of a post-training NN, which are central to its decision-making process, often remain opaque, obscuring the model's internal logic and reasoning. This opacity stands in contrast to traditional methods of scientific inquiry, for which understanding the mechanism of action is as crucial as the outcomes themselves. The quest for models
whose post-training status can be better understood
is actually not merely academic but has practical implications for trust, reliability, and furtherance of knowledge
\cite{rudin2022interpretable}.
This has given rise to the search of NNs which are human-understandable, either by attempting to explain
some of the values of their trained weights
\cite{alvarez2018towards,blazek2021explainable}, or by designing models
in which some of the trained weights can be
directly interpreted
\cite{akhtar2023survey,rudin2019stop,rudin2022interpretable,Zhang:2021}.
 In fact, although it might appear at first sight that
 interpretability and performance of NN models are
inversely correlated \cite{ akhtar2023survey},
several findings seem to contrast this view \cite{rudin2019stop}, as we also realize in the present research.

Recent
works have devised interpretable approaches mainly along two directions: (a) designing interpretable models from scratch, and (b) augmenting black-box models with interpretable components.
This work follows the first route and tackles the challenge of making NNs more interpretable by designing a network architecture that directly links the trained weights to the physical parameters of interest in the context of the discrete inverse conductivity problem. In the conductivity inverse problem, also known as electrical impedance tomography (EIT), the goal is to determine the electrical conductivity profile inside an object using measurements taken on the boundary of the given object. This typically involves applying an electrical current on the boundary and measuring the resulting voltage on the boundary, or vice versa. The measurements are then used to infer the internal conductivity profile, which is crucial for several applications including medical imaging, nondestructive testing of materials, and geophysical exploration \cite{Borcea}.

Alessandrini  \cite{Aless1988,A1989} proved that with a priori smoothness assumptions on the unknown conductivity, logarithmic stability holds. Mandache \cite{M} proved  that with a priori smoothness assumptions on the unknown conductivity, a conditional stability estimate of logarithmic type is the best possible. This indicates that the conductivity inverse problem is severely ill-posed, meaning that even a very low level of noise may lead to significant errors in the reconstruction, necessitating the use of regularization techniques to stabilize the solution procedure. One natural question that then arises when solving the problem is whether it is possible to find remedies to the ill-posedness by replacing the a priori smoothness assumptions with different ones that are physically relevant and may lead to better, possibly optimal Lipschitz stability. Indeed, it can be proved that if the conductivity belongs to a finite-dimensional subspace or manifold of a Banach space, then Lipschitz stability holds (see, e.g.,  \cite{curtis1990determining,AV,BF}), leading to a robust reconstruction \cite{BerMicPerSan18}. However, it is hard to establish the Lipschitz constant's dependence on the distribution of the unknowns. Therefore, the question is how to parameterize the unknown conductivity in numerical inversion so that we can control its stability without needing excessive regularization. In this direction there is an extensive literature on a discrete version of the inverse conductivity problem, i.e., resistor networks, starting from the pioneering works of Curtis and Morrow \cite{curtis1990determining,curtis1991DNmap}. In the works \cite{borcea2008electrical,BorceaDruskin:2010pyramid,BorceaMamonov:2010,Borcea:2013}, by viewing the resistor networks as the reduced-order
model of EIT, Borcea et al thoroughly explored the connection between circular planar networks and EIT for full and partial data, especially in determining optimal grids (i.e., grids that give spectrally accurate approximations of the data).
Other important theoretical contributions on resistor networks can be found in, e.g.,  \cite{CurtisMorrow:1994,Ingerman:2010,LamPylyavskyy:2012,BoyerGarzella:2016,BlastenIsozaki:2023}.

Here, to regularize the problem, we follow the approach and some of the results of Curtis and Morrow in \cite{curtis1990determining,curtis1991DNmap} on  the discrete inverse conductivity problem (posed on a network of a square lattice) and then propose a specific NN to recover the conductivity profile. With diverse modalities, NNs have been recently used extensively for EIT imaging  \cite{HamiltonHauptmann:2018,seo2019learning,li2019novel,Fan2020,GuoJiang:2021,CenJin:2023}. These techniques can be roughly classified into two categories, postprocessing versus direct inversion. The works \cite{HamiltonHauptmann:2018,GuoJiang:2021,CenJin:2023} employ convolutional neural networks (more precisely U-net \cite{ronneberger2015u} and its variants) to post-process an initial estimate provided by a direct reconstruction algorithm (e.g., D-bar method, direct sampling method and Calder\'{o}n's method).  Fan and Ying \cite{Fan2020} constructed a neural network based approach for EIT using a low-rank representation of the forward map. See the survey \cite{TanyuMaass:2023} for relevant methods and a comparative study of several deep learning based methods.

The present work contributes to the literature in the following two aspects. First,
we demonstrate that NNs can effectively solve non-linear inverse problems, e.g., the discrete inverse conductivity problem, with precision higher than traditional methods; see the numerical illustrations in Section \ref{NUME}. We also provide a sensitivity analysis of the minimizer in the small noise regime. Second, and more importantly, we show that by carefully designing the network architecture, the weights of the trained NN can be made interpretable, offering valuable insights into the problem structure and solution process. Specifically, the solution is encoded in the weights of the NN's second layer (see Theorem \ref{Thm1} for the precise statement), providing a clear interpretation of these weights in terms of  the conductivities in the inverse problem. This approach addresses the challenge of supervised disentanglement in NNs \cite{rudin2019stop} by ensuring that the internal representations of the network are transparent and meaningful. Additionally, by imposing constraints on the weights and ensuring sparsity, the method enhances the interpretability and stability of the solution, contributing to developing interpretable machine learning approaches for solving nonlinear inverse problems.

Now we situate the work in the context of interpretable machine learning. Rudin et al \cite{rudin2022interpretable} identify ten critical challenges, including supervised disentanglement of NNs, the application of constraints to encourage sparsity, and the characterization of the ``Rashomon set" of models, with which this work aligns well; see Remarks \ref{rmk:interp1} and \ref{rmk:interp2} for detailed connections. The work \cite{alvarez2018towards} proposes the concept of self-explaining NNs to improve the interpretability of NNs by ensuring that the model's predictions can be directly attributed to specific, understandable components of the input. We take this further by designing an NN where the second layer's weights are not just interpretable but directly correspond to physical parameters of the problem, thereby contributing to developing NNs that are both robust and interpretable. In the literature, there are several works on designing NNs for matrix operations, e.g., structured trainable networks for matrix algebra \cite{wang1990structured}, feed-forward NN for matrix computation \cite{al2001feed}
and second-order NN for computing Moore-Penrose inverse of matrices \cite{li2022efficient}. {%We advance the discourse on making NN weights understandable post-training, and leverage NNs to solve a discrete inverse conductivity problem with the weights of the final layer directly relating to physical parameters.}
Our study indicates that the NN is not just a computational tool but also a means of embedding and extracting interpretable mathematical structures, thereby advancing its use in solving complex problems \cite{chakraborty2024explainable}.

The rest of the paper is organized as follows. In Section \ref{CondProb}, we describe the problem of recovering the conductivity profile from the discrete Dirichlet to Neumann data on a square lattice. Then in Section \ref{ARCH}, we propose a novel feed-forward neural network architecture to solve the discrete inverse conductivity problem directly by the trained weights in the second layer. In Section \ref{MainResults} we show that if the sample size is large enough, all minima of the loss function have the same value of the weights in the second layer. In Section \ref{NUME}, we present numerical results to illustrate the NN approach, including a comparison with the Curtis-Morrow algorithm. Last, in Section \ref{SENS}, we provide a sensitivity analysis of the minimizer to the problem under a small noise assumption.

\section{The discrete conductivity problem}\label{CondProb}
Now we describe the problem of recovering
the conductivity profile from boundary measurements, i.e.,
from the Dirichlet to Neumann data on
a square grid of the form
\begin{eqnarray}  \label{Domain}
D=\{p \, |\, p=(i,j), \, i,j=1,\dots,n\},
\end{eqnarray}
with the boundary points
\begin{eqnarray} \label{boundary_points}
\partial D =\{q \, |\,  q=(0,j),\, j=1, \dots, n  \text{ or } q=(n+1,j),\, j=1, \dots, n \\
\quad \text{ or } q=(i,0),\, i=1, \dots, n \text{ or } q=(i,n+1),\, i=1, \dots, n
 \}.\nonumber
\end{eqnarray}
Let \(\overline D= D \cup \partial D, \)
and for each $p \in \overline D$, we indicate by $\mathcal N(p)$ the set of adjacent nodes to $p$.
The points in $D$ have four adjacent points, and the points in $\partial D$ have only one adjacent point in the interior (in $D$). The bonds between nearest neighbor nodes are given by
\begin{eqnarray}  \label{BoundaryIndices}
B =\{\{p,q\} \, |\, p=(i,j), q=(i',j'), \,
(i,j), (i',j') \in \overline D,
|i-i'|+|j-j'|=1\}.\nonumber
\end{eqnarray}
Below we will use $ij$ in the subscript to indicate $(i,j)$ and $ijhk$ instead of $\{(i,j),(h,k)\}$. See Fig. \ref{fig:network}(a) for a schematic illustration.

Given a conductivity profile \( \bm\gamma=[{\gamma}_{\{p,q\}}]_{\{p,q\} \in B} \), the discrete conductivity problem with the Dirichlet boundary condition $ {\bf \overline u}=[\overline u_q]_{q \in \partial D}$ consists of solving the following discrete conductivity equation for the variable ${\bf u}=[u_p]_{p \in\overline{D}}$:
  \begin{eqnarray} \label{EQ1}
     \sum_{q \in \mathcal N(p)} {\gamma}_{\{q,p\}}
(u_{q} - u_{p}) =0,\quad p\in D,
 \end{eqnarray}
subject to the following Dirichlet boundary condition
 \begin{equation}\label{EQ2}
 u_q=\overline u_{q},\quad q\in \partial D.
 \end{equation}
Given the conductivity profile $ {\bm \gamma}$, by \cite[Proposition 2.4]{curtis1990determining}, there is a unique solution $\mathbf{u}$ of problem \eqref{EQ1}--\eqref{EQ2}. For any edge $\{p,q\}\in B$, the induced current $I_{pq}$ is defined by
\begin{equation*}
  I_{pq}={\gamma}_{\{q,p\}} (u_{q} - u_{p}), \quad  \forall p,q\in\overline{D},
\end{equation*}
Similarly, the discrete Neumann problem satisfies equation \eqref{EQ1} and takes the boundary condition ${\bf \overline v}= [\overline v_q]_{q  \in  \partial D}$
such that
\begin{eqnarray} \label{Eq3}
\overline v_{q}={\gamma}_{\{q,p\}} (u_{q} - u_{p}), \quad
 q  \in  \partial D,
 p \in \mathcal N(q).
 \end{eqnarray}
The Dirichlet to Neumann (DtN) data consists of a pair
 \( {\bf d} =( {\bf \overline u},{\bf \overline v}    )  \),
where ${\bf \overline u}$ is a Dirichlet boundary data,
and ${\bf \overline v}$ consists of the corresponding flux measurement
on $\partial D$. Since the map is linear, we denote by $\Lambda_{{\bm\gamma}}$
the matrix that maps the vector ${\bf \overline u}$
to the vector ${\bf \overline v}$.

\begin{figure}[hbt!]
\begin{center}
\setlength{\tabcolsep}{0pt}
\begin{tabular}{cc}
 \begin{tikzpicture}
    % Set up the grid dimensions
    \newcommand{\n}{6} % Number of grid points on one side

    % Draw the grid
    \foreach \x in {1,...,\n} {
        \foreach \y in {1,...,\n} {
            \filldraw (\x,\y) circle (2pt); % Draw nodes
        }
    }

    % Draw the boundary
    \foreach \x in {1,...,\n }
        {\filldraw (\x,0) circle (2pt); % Draw nodes
        }
        % Draw the boundary
    \foreach \x in {1,...,\n }
        {\filldraw (\x,\n+1) circle (2pt); % Draw nodes
        }
    % Draw the boundary
    \foreach \y in {1,...,\n }
        {\filldraw (0,\y) circle (2pt); % Draw nodes
        }
         % Draw the boundary
    \foreach \y in {1,...,\n }
        {\filldraw (\n+1,\y) circle (2pt); % Draw nodes
        }

    % Draw the horizontal lines
    \foreach \y in {1,...,\n} {
        \draw[thick] (0,\y) -- (\n+1,\y);
    }

    % Draw the vertical lines
    \foreach \x in {1,...,\n} {
        \draw[thick] (\x,0) -- (\x,\n+1);
    }

    % Add labels on the right side
    \node[below] at (\n,0) {$u_{n+1,n}$};
    \node[right] at (\n,1.2) {$u_{n,n}$};
    \node[right] at (\n,2) {$\vdots$};

    % Add labels on the left side

    \node[left] at (0,2) {$\vdots$};

    % Add labels on the top side
    \node[left] at (0,5.9) {$u_{1,0}$};
    \node[above] at (1.3,6.9) {$u_{0,1}$};
    \node[above] at (1.3,5.9) {$u_{1,1}$};
    \node[above] at (2,\n) {$\cdots$};

    % Add labels on the bottom side

    \node[below] at (1,0) {$u_{n+1,1}$};
    \node[below] at (2,0) {$\cdots$};

    % Add label w next to one of the horizontal lines
    \node[above] at (2.5,3) {$w$};  % Adjust the position and the value of \y as needed
\end{tikzpicture} &
\includegraphics[width=0.5\linewidth]{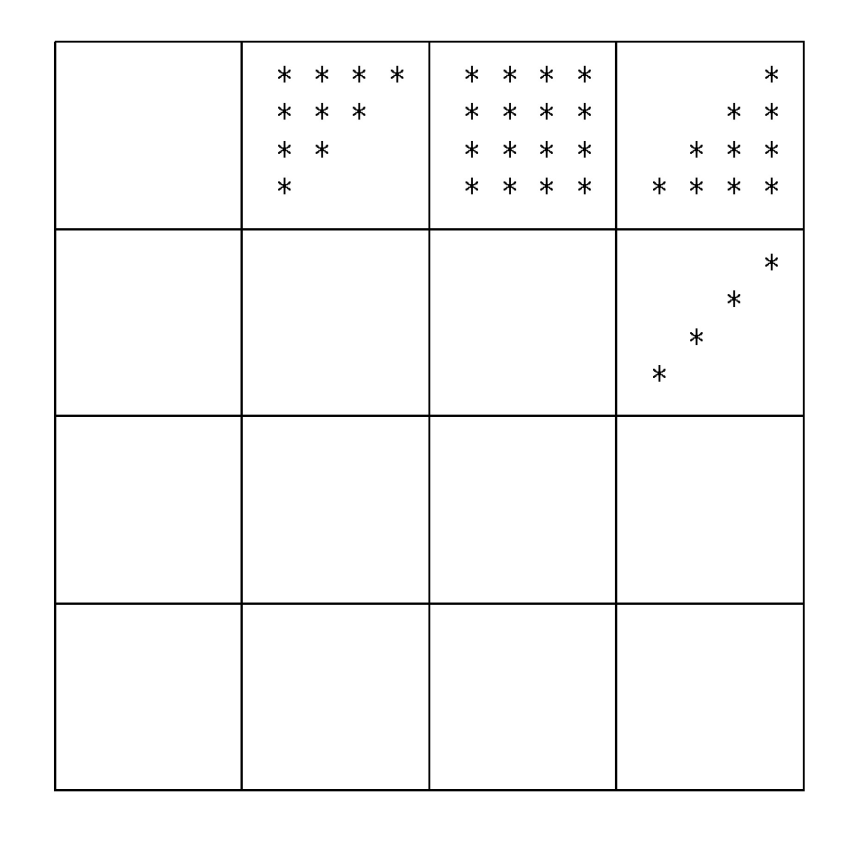}\\
(a) network on square lattice & (b) partial DtN data
\end{tabular}
\end{center}
\caption{Schematic illustrations of (a) the network and (b) partial DtN data.}
    \label{fig:network}
\end{figure}
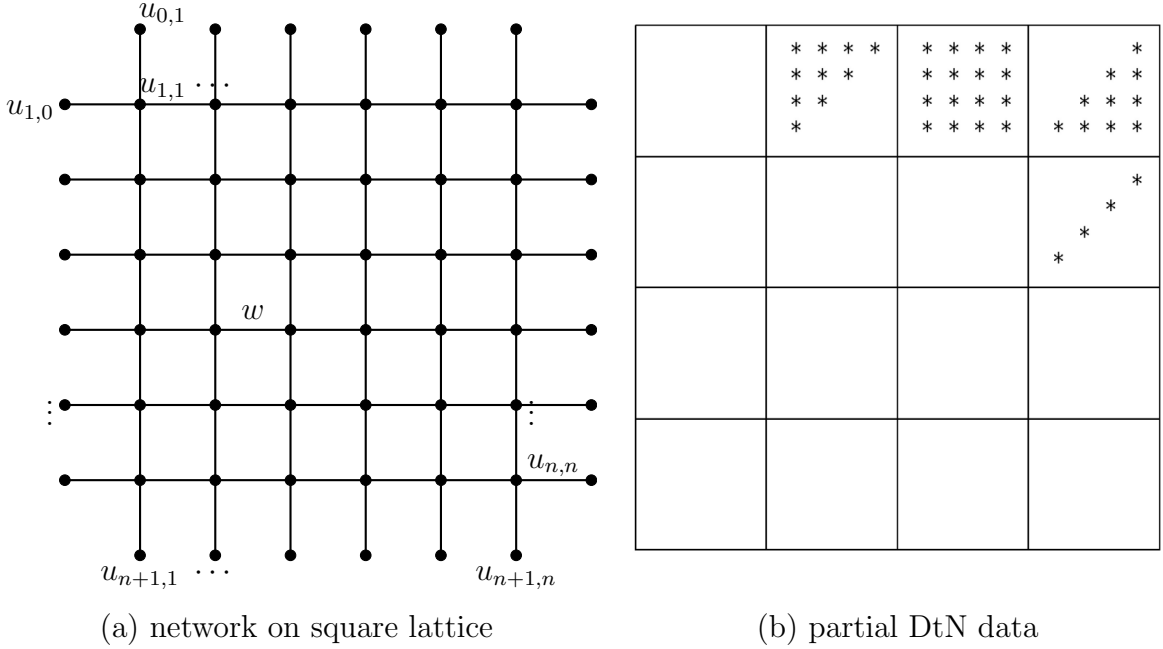

The inverse problem consists of retrieving
 the conductivity $\bm\gamma$ from
 $m$ pairs of Cauchy data $\{ {\bf d}(k)\}_{k=1}^m$.
By \cite[Theorem 3.2]{curtis1990determining}, a knowledge of the full DtN map $\Lambda({\bm\gamma}):{\bf\overline u} \rightarrow
 {\bf\overline v}$
 uniquely determines the conductivity $\bm{\gamma}$.
This a priori would require $4n$ pairs \( {\bf d} (k)=( {\bf \overline u}(k),{\bf \overline v}(k)    )  \) where the vectors
$\{{\bf \overline u}(k)\}_{k=1}^m$ form a basis of the vector space of the Dirichlet boundary data. Remarkably, Curtis and Morrow
\cite[Theorem 5.1]{curtis1991DNmap} showed that the DtN matrix $\Lambda_{\bm{\gamma}}$
actually depends only on $2n(n+1) $ parameters,
which are specific entries of $\Lambda_{\bm{\gamma}}$ (see Fig. \ref{fig:network}(b) for the pattern),
and from their positions one can deduce that actually $3n$ measurements which include $3n$ specific vectors of the basis of the vector space of the Dirichlet boundary data are sufficient for the uniqueness.

\section{The neural network architecture}\label{ARCH}

In this section, we propose a novel feed-forward neural network
with the property that with an enough number of samples of the DtN data, it solves the
discrete inverse conductivity problem directly by the learned weights in the second layer.

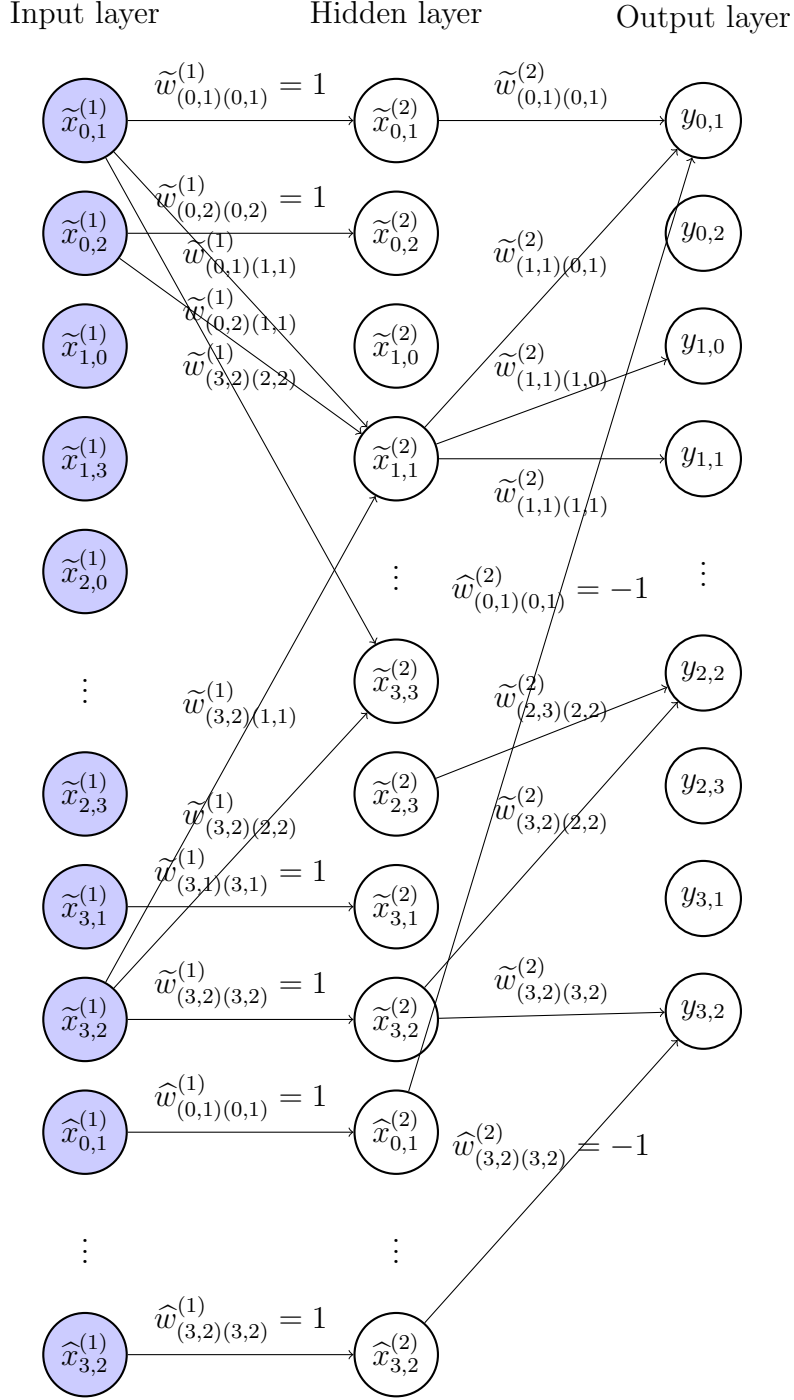
\begin{figure}[h!]
\centering
\begin{tikzpicture}[
    node distance=1.5cm,
    layer/.style={rectangle, draw, thick},
    neuron/.style={circle, draw, thick, minimum size=1cm,inner sep=2pt },
    weight/.style={auto}
]

%\node[circle, draw, fill=blue!20, minimum size=1.2cm] (x1) at (0,4) {$x_1$};

% Input Layer Nodes
\node[neuron,fill=blue!20] (i1) {\makebox[0.7cm][c]{$\widetilde x^{(1)}_{ 0,1 }$}};
\node[neuron, below of=i1,fill=blue!20] (i2) {\makebox[0.7cm][c]{$\widetilde x^{(1)}_{0,2}$}};
\node[neuron, below of=i2,fill=blue!20] (i3){\makebox[0.7cm][c] {$\widetilde x^{(1)}_{1,0}$}};
\node[neuron, below of=i3,fill=blue!20] (i4) {\makebox[0.7cm][c]{$\widetilde x^{(1)}_{1,3}$}};
\node[neuron, below of=i4,fill=blue!20] (i5) {\makebox[0.7cm][c]{$\widetilde x^{(1)}_{2,0}$}};
\node[below=0.5cm of i5] (dots1) {$\vdots$};
\node[neuron, below =0.5cm of dots1,fill=blue!20] (i6) {\makebox[0.7cm][c]{$\widetilde x^{(1)}_{2,3}$}};
\node[neuron, below of=i6,fill=blue!20] (i7) {\makebox[0.7cm][c]{$\widetilde x^{(1)}_{3,1}$}};
\node[neuron, below of=i7,fill=blue!20] (i8) {\makebox[0.7cm][c]{$\widetilde x^{(1)}_{3,2}$}};
\node[neuron, below of=i8,fill=blue!20] (i9) {\makebox[0.7cm][c]{$\widehat{x}^{(1)}_{0,1}$}};
\node[below=0.5cm of i9] (dots1) {$\vdots$};
\node[neuron, below =0.5cm of dots1,fill=blue!20] (i10) {\makebox[0.7cm][c]{$\widehat{ x}^{(1)}_{3,2}$}};

% Hidden Layer Nodes
\node[neuron, right=3cm of i1] (h1) {\makebox[0.7cm][c]{$\widetilde x^{(2)}_{0,1}$}};
\node[neuron, below of=h1] (h2) {\makebox[0.7cm][c]{$\widetilde x^{(2)}_{0,2}$}};
\node[neuron, below of=h2] (h3) {\makebox[0.7cm][c]{$\widetilde x^{(2)}_{1,0}$}};
\node[neuron, below of=h3] (h4) {\makebox[0.7cm][c]{$\widetilde x^{(2)}_{1,1}$}};
\node[below=0.5cm of h4] (dots1) {$\vdots$}; % Dots to indicate continuation
\node[neuron, below=0.5cm of dots1] (h5) {\makebox[0.7cm][c]{$\widetilde x^{(2)}_{3,3}$}};

\node[neuron, below of=h5] (h6) {\makebox[0.7cm][c]{$\widetilde x^{(2)}_{2,3}$}};
\node[neuron, below of=h6] (h7) {\makebox[0.7cm][c]{$\widetilde x^{(2)}_{3,1}$}};
\node[neuron, below of=h7] (h8) {\makebox[0.7cm][c]{$\widetilde x^{(2)}_{3,2}$}};
\node[neuron, below of=h8] (h9) {\makebox[0.7cm][c]{$\widehat x^{(2)}_{0,1}$}};
\node[below=0.5cm of h9] (dots1) {$\vdots$}; % Dots to indicate continuation
\node[neuron, below=0.5cm of dots1] (h10) {\makebox[0.7cm][c]{$\widehat x^{(2)}_{3,2}$}};

% Output Layer Nodes
\node[neuron, right=3cm of h1] (o1) {\makebox[0.7cm][c]{$y_{0,1}$}};
\node[neuron, below of=o1] (o2) {\makebox[0.7cm][c]{$y_{0,2}$}};
\node[neuron, below of=o2] (o3) {\makebox[0.7cm][c]{$y_{1,0}$}};
\node[neuron, below of=o3] (o4) {\makebox[0.7cm][c]{$y_{1,1}$}};
\node[below=0.5cm of o4] (dots1) {$\vdots$}; % Dots to indicate continuation
\node[neuron, below=0.5cm of dots1] (o5) {\makebox[0.7cm][c]{$y_{2,2}$}};

\node[neuron, below of=o5] (o6) {\makebox[0.7cm][c]{$y_{2,3}$}};
\node[neuron, below of=o6] (o7) {\makebox[0.7cm][c]{$y_{3,1}$}};
\node[neuron, below of=o7] (o8) {\makebox[0.7cm][c]{$y_{3,2}$}};

% Connections Input to Hidden Layer
\draw[->] (i1) -- (h1) node[midway, above] {$\widetilde w^{(1)}_{(0,1)(0,1)}=1$};
\draw[->] (i2) -- (h2) node[midway, above] {$\widetilde w^{(1)}_{(0,2)(0,2)}=1$};
\draw[->] (i1) -- (h4) node[midway, above] {$\widetilde w^{(1)}_{(0,1)(1,1)}$};
\draw[->] (i2) -- (h4) node[midway, above] {$\widetilde w^{(1)}_{(0,2)(1,1)}$};
\draw[->] (i8) -- (h4) node[midway, above] {$\widetilde w^{(1)}_{(3,2)(1,1)}$};
\draw[->] (i1) -- (h5) node[midway, above] {$\widetilde w^{(1)}_{(3,2)(2,2)}$};
\draw[->] (i8) -- (h5) node[midway, above] {$\widetilde w^{(1)}_{(3,2)(2,2)}$};
\draw[->] (i7) -- (h7) node[midway, above] {$\widetilde w^{(1)}_{(3,1)(3,1)}=1$};
\draw[->] (i8) -- (h8) node[midway, above] {$\widetilde w^{(1)}_{(3,2)(3,2)}=1$};
\draw[->] (i9) -- (h9) node[midway, above] {$\widehat {w}^{(1)}_{(0,1)(0,1)}=1$};
\draw[->] (i10) -- (h10) node[midway, above] {$\widehat{w}^{(1)}_{(3,2)(3,2)}=1$};
% Add more connections as needed

% Connections Hidden to Output Layer
\draw[->] (h1) -- (o1) node[midway, above] {$\widetilde w^{(2)}_{(0,1)(0,1)}$};
\draw[->] (h4) -- (o4) node[midway, below] {$\widetilde w^{(2)}_{(1,1)(1,1)}$};
\draw[->] (h4) -- (o1) node[midway ,above ] {$\widetilde w^{(2)}_{(1,1)(0,1)}$};
\draw[->] (h4) -- (o3) node[midway, above] {$\widetilde w^{(2)}_{(1,1)(1,0)}$};
% \draw[->] (h1) -- (o4) node[midway,above] {$w^{(2)}_{(0,1)(1,1)}$};
% \draw[->] (h2) -- (o2) ;
% \draw[->] (h2) -- (o1) ;
% \draw[->] (h3) -- (o3) ;
% \draw[->] (h3) -- (o4) ;
\draw[->] (h8) -- (o8) node[midway, above] {$\widetilde w^{(2)}_{(3,2)(3,2)}$};
\draw[->] (h8) -- (o5) node[midway, above] {$\widetilde w^{(2)}_{(3,2)(2,2)}$};
\draw[->] (h6) -- (o5) node[midway, above] {$\widetilde w^{(2)}_{(2,3)(2,2)}$};
\draw[->] (h9) -- (o1) node[midway, above] {$\widehat { w}^{(2)}_{(0,1)(0,1)}=-1$};
\draw[->] (h10) -- (o8) node[midway, above] {$\widehat {w}^{(2)}_{(3,2)(3,2)}=-1$};

% % Connections Hidden to Output Layer
% \draw[->] (h1) -- (o1) node[midway, above] {$w^{(2)}_{(0,1)(0,1)}=-0.51$};
% \draw[->] (h4) -- (o4) node[midway, below] {$w^{(2)}_{(1,1)(1,1)}=-2.48$};
% \draw[->] (h4) -- (o1) node[midway ,above ] {$w^{(2)}_{(1,1)(0,1)}=0.51$};
% \draw[->] (h4) -- (o3) node[midway, above] {$w^{(2)}_{(1,1)(1,0)}=0.55$};
% % \draw[->] (h1) -- (o4) node[midway,above] {$w^{(2)}_{(0,1)(1,1)}=0.51$};
% % \draw[->] (h2) -- (o2) ;
% % \draw[->] (h2) -- (o1) ;
% % \draw[->] (h3) -- (o3) ;
% % \draw[->] (h3) -- (o4) ;
% \draw[->] (h8) -- (o8) node[midway, above] {$w^{(2)}_{(3,2)(3,2)}=-0.99$};
% \draw[->] (h8) -- (o5) node[midway, above] {$w^{(2)}_{(3,2)(2,2)}=0.99$};
% \draw[->] (h6) -- (o5) node[midway, above] {$w^{(2)}_{(2,3)(2,2)}=0.91$};

% Add more connections as needed

% Labels for some weights (example)
\node[above=0.5cm of i1] (label) {Input layer};
\node[above=0.5cm of h1] {Hidden layer};
\node[above=0.5cm of o1] {Output layer};

% Propagation arrows
%\draw[->, thick, color={rgb:red,1;green,73;blue,70}] (-1, 2) -- (8, 2) node[above, pos=0.5] {Information forward propagation};
%\draw[->, thick, color={rgb:red,143;green,52;blue,48}] (8, -13) -- (-1, -13) node[below, pos=0.5] {Error back propagation};
\end{tikzpicture}
\caption{The architecture of the proposed FNN.}
\label{fig:FNN}
\end{figure}

The proposed architecture has three layers, and there are \(8 n\), \(n^2+8n\), and \(n^2+4n\) neurons in the first, second and third layers, respectively; see Fig. \ref{fig:FNN} for an illustration of the architecture.
Each of the first two layers and the associated weight matrix are divided
into two parts, indicated by the symbols $~\widetilde{}~$  and $~\widehat{}~$; the first part involves the solution of the direct problem with Dirichlet boundary
conditions, and the second part only transmits the Neumann data to the loss.

The neurons are labeled $\widetilde x^{(1)}_{p}$, $p \in \partial D$ and
${\widehat  x}^{(1)}_{p}$, $p \in \partial D$
in the first layer (input);
  $\widetilde x^{(2)}_{p}$, $p \in \overline D$, and
${\widehat x}^{(2)}_{p}$, $p \in \partial D$   in the second layer;
$y_{p}$, $p \in \overline D$  in the third layer. Physically, $\widetilde x_p^{(1)}$ and $ {\widehat x}_p^{(1)}$ denote the input voltage (Dirichlet boundary data) and the measured current data
(Neumann boundary data) on the boundary $\partial D$,
respectively;
$\widetilde x_p^{(2)}$  denotes the potential over the entire lattice,
while ${\widehat x}_p^{(2)}$ is still the measured current data; and $y_p$ reflects Kirchhoff law and the
equality between the boundary flux and the Neumann data.

We denote by $W^{(1)}$ and $W^{(2)}$ the weight matrices in the first and second layers, respectively, and construct them as follows.
Between the first and second layers,
we fix weights \(\widetilde w^{(1)}_{qq}\equiv 1, q \in \partial D\), weights  \(\widetilde w^{(1)}_{qp}\)
between each neuron $\widetilde x^{(1)}_{q}, q \in \partial D$, in the first layer and each neuron $\widetilde x^{(2)}_{p}, p \in D$, in the hidden
layer. Physically, these two conditions respectively reflect the Dirichlet boundary conditions and discrete Green's kernel representing the mapping from the Dirichlet boundary data to the interior nodes. Next we fix weights \({\widehat w}^{(1)}_{qq}\equiv 1, q \in \partial D\) between each neuron ${\widehat x}^{(1)}_{q}, q \in \partial D$,
in the first layer and each neuron ${\widehat x}^{(2)}_{q}, q \in \partial D$, in the hidden
layer, and fix the weights to zero otherwise. This determines a matrix $\widetilde W^{(1)}\in \mathbb{R}^{(n^2+4n)\times 4n}$, and a
matrix $\widehat { W}^{(1)}=I_{4n}\in \mathbb{R}^{4n\times 4n}$, where $I_{4n}$ is the identity matrix in $\mathbb R^{4n}$. These two matrices can be combined in a block matrix
\begin{equation}
    W^{(1)} = \begin{bmatrix}
        \widetilde W^{(1)} & 0 \\
        0 & I_{4n}
    \end{bmatrix}.
\end{equation}
The output of the \(p\)-th neuron of the hidden layer  is then given by
\begin{align}
\left\{\begin{aligned}
\widetilde  x^{(2)}_{p} &= \sum_{q \in \partial D} \widetilde w^{(1)}_{qp}
\widetilde x^{(1)}_{q},\quad p \in \overline D, \\
   {\widehat x}^{(2)}_{p} &=
   {\widehat x}^{(1)}_{p},\quad p \in \partial D.
\end{aligned}\right.
   \label{Eq2}
  \end{align}

Between the hidden layer and the output layer,  there is a weight \(\widetilde w^{(2)}_{pr}\)
between the neurons  $\widetilde x^{(2)}_{p}, p \in \overline D$ and  $y_{r}, r \in \overline D$ if  \(\{p,r\} \in B\) or if $p=r$  with the following constraints: (i) the weights are symmetric, i.e., $\widetilde w^{(2)}_{pr}=
\widetilde w^{(2)}_{rp}$; (ii) $\widetilde w^{(2)}_{rr}=   - \sum_{p \in \mathcal N(r)}  \widetilde w^{(2)}_{pr}$ for $r\in D$. Physically, these conditions mimic the symmetry and conservation property of the conductivities on network edges.
 Furthermore, there is a weight
$\widehat w^{(2)}_{pp}=-1$ between the neurons  $ {\widehat  x}^{(2)}_{p}$ and  $y_{p}$ for $ p \in \partial D$. The contributions from $\widetilde{x}_p^{(2)}$ and $\widehat{x}_p^{(2)}$ form the difference between the given Neumann data and the one predicted from the NN.
All other weights are fixed at zero. The weights between the
second and third layers can be summarized as
\begin{equation}
    W^{(2)} = \begin{bmatrix}
        \widetilde W^{(2)} & - I_{4n}
    \end{bmatrix}.
\end{equation}
The output of  the \(r\)-th  neuron of the output layer is given by
\begin{equation}
y_{r} = \begin{cases}
\sum_{p \in \mathcal N(r)}  \widetilde w^{(2)}_{pr}
(\widetilde x^{(2)}_{p} - \widetilde x^{(2)}_{r}), \quad r \in D,\\
\widetilde w^{(2)}_{pr}
(\widetilde x^{(2)}_{r} -\widetilde  x^{(2)}
_{p})-\widehat {x}^{(2)}_{r}, \quad r \in \partial D,
p \in \mathcal N(r).
\end{cases}
\end{equation}

Below we denote $\widetilde W=\{\widetilde W^{(1)}, \widetilde W^{(2)}\}$, which is learned from the training data. The activation functions are simply the identity. Note that in the architecture, some of
the weights are fixed at prescribed values
like  $1$  or $0$ \cite{han2015learning}, and
some other pairs of weights are forced
to take the same values \cite{yuan2006model}.

Given the input vector pairs \({\bf {\widetilde x}}(k) = [{\widetilde x}_{q}(k)]^\top_{q\in \partial D}\in\mathbb{R}^{4n}\), \({\bf{\widehat {x}}}(k) = [{\widehat x}_{q}(k)]^\top_{q\in \partial D}\in\mathbb{R}^{4n}\),
with $k=1,\dots, m$ ($m$ is the sample size), let
\begin{equation*}
{\bf x}(k)
= \left[\begin{aligned} {{\widetilde{\textbf{x}}}^{(1)}}(k)\\ {{\widehat {\textbf{x}}}}^{(1)}(k)\end{aligned}\right] \in \mathbb{R}^{8n}.
\end{equation*}
Then with the weight matrices  \(W^{(1)}\) and \(W^{(2)}\)
between the  layers,
the output vector \({\bf y}(k) = [y_{p}(k)]^\top_{p\in
\overline D}\), with $k=1,\dots, m$,
obtained from the forward propagation, can be expressed in a matrix-vector form:
\begin{equation}
{\bf y}(k) = W^{(2)} (W^{(1)} {\bf  x}(k)).
\end{equation}

A collection of the DtN data  $\{({{\widetilde {\bf x}}^{(1)}}(k), {{\widehat{\bf x}}}^{(1)}(k))\}_{k=1}^m$ is used as the  training data, where \({ {\widetilde {\bf x}}^{(1)}}(k)={\bf \overline u}(k)\)
is the Dirichlet boundary data and \({ \widehat {\bf x}^{(1)}}(k)={\bf \overline v}(k),\,q \in
\partial D\) is the corresponding Neumann
data. Then we employ the mean square error as the loss
\begin{equation} \label{loss}
\begin{aligned}
 C_{\alpha}(\widetilde W)
&= \frac{1}{2m} \sum_{k=1}^{m}
\left(\sum_{p
\in \partial D} ( y_{p}(k))^2
+ \alpha \sum_{p
\in  D} ( y_{p}(k))^2
\right),
\end{aligned}
\end{equation}
where the scalar $\alpha>0$ can be suitably chosen to improve the training of the weights $\widetilde W$.

In Section \ref{MainResults} below, we show that for exact data,
if the sample size $m$ is large enough, then every
global minimizer
$ \widetilde W_*=(\widetilde W^{(1)}_*, \widetilde W^{(2)}_*)$
of the loss $C_\alpha(\widetilde W^{(1)},\widetilde W^{(2)})$ is such that
$C_\alpha(\widetilde W^{(1)}_*, \widetilde W^{(2)}_*)=0$ and
$\widetilde W^{(2)}_*$
is equivalent to the conductivity $\bm \gamma^\dag$ {corresponding to the DtN data  $\{({ \widetilde{\bf x}^{(1)}}(k), {{\widehat{\bf x}}}^{(1)}(k))\}_{k=1}^m$}  .
If the sample size $m$ increases further,
then there is a unique minimum,
and the weights in the matrix $\widetilde  W_*^{(1)}$
correspond to the coefficients
(which are nonlinear in the conductivity
) of the
linear transformation which
determines the solution of
the Dirichlet problem given the
boundary condition.

\begin{remark}\label{rmk:interp1}
The proposed NN architecture incorporates specific constraints, e.g., fixing some weights to predefined values or forcing some weights to be identical. These constraints not only help stabilize the solution process but also make the model more interpretable by reducing the number of free parameters, which aligns with the goal of encouraging sparsity and improving interpretability. Moreover, it directly ties the trained weights of the NN to physical parameters of the problem, i.e., the conductivity $\bm \gamma$. By its design, the NN satisfies the property that the weights of the second layer correspond to the conductivity values, which effectively ``disentangle" internal representations of the NN. This makes the NN more interpretable since the weights have a clear, understandable meaning tied to the physical domain. This is a crucial step in making NNs more transparent and easier to analyze in high-stake applications where understanding the model behavior is critical.
\end{remark}
\begin{remark}\label{rmk:interp12}
The weights in the first layer mimic the
 matrix that maps the Dirichlet boundary values to the solution
in the interior, i.e., discrete Green kernel,
which is a nonlinear function of the conductivity $\bm\gamma$.
Corollary \ref{corollary} below shows that with enough data, the weights in $\widetilde W_*^{(1)}$ coincide with the discrete Green kernel, and it is thus a nonlinear
function of $\widetilde W_*^{(2)}$. However,
it looks like that the NN takes advantage of not
assuming such dependence during training; especially
for noisy or missing data, the NN seems to exploit
the additional flexibility to achieve better performance
{\rm(}see Remark \ref{rmk:interp21}{\rm)}.
\end{remark}

\section{Optimal weights for noiseless data} \label{MainResults}
We now identify some properties of
the minimizers of the loss $C_\alpha(\widetilde{W})$.
First we show that if the sample size
$m \geq 3n$
and the Dirichlet data are
suitably chosen, then all minima
 $\widetilde W_*=(\widetilde W^{(1)}_*, \widetilde W^{(2)}_*)$ of the loss $C_\alpha(\widetilde W)$ have the same value of $ \widetilde W_*^{(2)}$ and  all of
the  conductivities $\bm\gamma$ which
solve the inverse problem with the given DtN data can be read from the entries of
$ \widetilde W_*^{(2)}$.
More precisely, we have the following result.

\begin{theorem} \label{Thm1} Fix $\alpha>0$. Given a conductivity profile
$\bm{\gamma}=[\gamma_{\{p,q\}}]_{\{p,q\}\in B}$, consider
$m \geq 3n$ Cauchy data pairs
\begin{equation*}
{\bf d}(k)=( \overline{\bf u}(k), \overline{\bf v}(k)), \quad k=1,\dots, m,
\end{equation*}
solving problem \eqref{EQ1} for the given
${\bm\gamma}$ such that
\(\{\overline{\bf u}(k)\}_{k=1, \dots, m}\)
contain a basis of the vector space of the
subspace of $\mathbb R^{4n}$
generated by the last \(3n\) coordinates.
Next, consider the FNN in \eqref{Eq2} and
\eqref{Eq3}. Then for the input data $
{\bf \widetilde x}^{(1)}(k)= \overline{\bf u}(k)$ and $ {\bf {\widehat x}}^{(1)}(k)=\overline {\bf v}(k)$, $k=1,\ldots,m$, there are minima  ${\widetilde W}_*$
of the loss $C_\alpha(\widetilde W)$ in \eqref{loss}
such that $C_\alpha({\widetilde W}_*)=0$. Moreover, for all such minima, there holds
\begin{eqnarray} \label{EQ10}
    \widetilde w^{(2)}_{pq }=
    \begin{cases}
        {\gamma}_{\{ p,q \}},  \qquad &\text{ if }
        ( p,q ) \in B, p,q \in D, \\
        -{\gamma}_{\{ p,q \}},  \qquad &\text{ if }
        ( p,q ) \in B, \{p,q\} \cap \partial D \neq \emptyset, \\
        0, \qquad &\text{ if } ( p,q ) \notin B, p \neq q, \\
        -\sum_{r\in\mathcal{N}(p)}
        {\gamma}_{\{r,p\}},
        \qquad &\text{ if } p = q \in D,\\
        {\gamma}_{\{ r,p\}},
        \qquad &\text{ if } p = q \in \partial D,
        (r,p) \in B.
    \end{cases}
\end{eqnarray}

\end{theorem}
\begin{proof}
In the proof, we fix $\alpha>0$, since the argument
below holds regardless of the value of $\alpha$. Given $\bm{\gamma}$ and for each fixed $k\in \{1,\dots m\}$,  there exists a unique solution
%to Neumann
%pair \({\bf d}(k)=({\bf \overline u}(k),{\bf \overline v}(k)) \) there exists a unique solution
\({\bf  u}(k)\)
to problem \eqref{EQ1}-\eqref{EQ2} with \({\bf \overline u}(k)\)
as the boundary datum \cite[Proposition 2.4]{curtis1990determining}. Consider the
matrix $A$ mapping the Dirichlet boundary datum $\overline{\bf u}$
into the solution ${\bf  u}$. The existence of $A$ is then guaranteed by the linearity of  problem \eqref{EQ1}-\eqref{EQ2}.
If we take $\widetilde W_*^{(1)}$
equal to $A$, so that ${\bf \widetilde x}_*^{(2)} = {\bf u}(k)$ for each $k$,
and $\widetilde W_*^{(2)}$ as in \eqref{EQ10},
then clearly $C_\alpha(\widetilde  W_*)=0$, and
the weights $W_*$ minimize the loss $C_\alpha(\widetilde W)$ with a loss value 0.

Next, suppose that  $\widetilde W= ({\widetilde  W}^{(1)},{\widetilde W}^{(2)})$ satisfies $C_\alpha({\widetilde  W})=0$, and
consider the values \(\widetilde {\bf x}^{(2)}(k)\) of the neurons in the hidden layer.
Necessarily, for each $k$, \( \widetilde{\bf x}^{(2)}(k)\)
is a solution of problem \eqref{EQ1}-\eqref{EQ2}
with the Dirichlet datum \({\bf \overline u}(k)\)
and conductivity profile
$\gamma_{\{q, p\}}=\widetilde w^{(2)}_{qp}$ for any $p,q \in D$ if $\{p,q\}\in B$,  and
$\gamma_{\{q,p\}}=-\widetilde w^{(2)}_{qp}$ for $q  \in  \partial D,
 p \in \mathcal N(q)$ for each $k$, and in addition
the flux corresponding to ${\bf {\widetilde x}}^{(2)}(k)$ on $\partial D$ is given by
the Neumann datum \({\bf \overline v}(k)\). In terms of the weight matrix $\widetilde W^{(2)}$, equations \eqref{EQ1}, \eqref{EQ2} and \eqref{Eq3} take the form
\begin{equation}\label{eq11}
\left\{\begin{aligned}\Big( \sum_{q \in \mathcal N(p)} \widetilde w^{(2)}_{qp}
\widetilde x^{(2)}_{q}\Big) - \widetilde{w}^{(2)}_{{pp}}\widetilde{x}^{(2)}_{p} &=0,\quad p\in D,\\
\widetilde x^{(2)}_q&=\overline u_{q}(k),\quad q\in \partial D,\\
 \widetilde  w^{(2)}_{qp} (\widetilde x^{(2)}_{q} - \widetilde x^{(2)}_{p})
  &=\overline v_{q}(k), \quad
 q  \in  \partial D,
 p \in \mathcal N(q),
\end{aligned}\right.
 \end{equation}
for $k=1,\dots,m$. Clearly, if $\widetilde W^{(2)}$ is as in \eqref{EQ10}, then it solves the system  \eqref{eq11} for every $k$. We only need to show that if $\widetilde W^{(2)}$ is a solution to problem \eqref{eq11} for all $k=1,\dots, m$, then it must be of the form \eqref{EQ10}. Indeed, using the symmetry of the weight matrix $\widetilde W^{(2)}$ and the property
 \begin{equation}\label{cond2}
    \widetilde  w^{(2)}_{rr}=
       - \sum_{p \in \mathcal N(r)}
    \widetilde    w^{(2)}_{pr},
\end{equation}
we obtain
\begin{equation} \label{W2}
\left\{\begin{aligned}
  \sum_{q \in \mathcal N(p)}\widetilde  w^{(2)}_{qp}\left(
\widetilde x^{(2)}_{q} - \widetilde x^{(2)}_{p}\right) &=0,\quad p\in D,\\
\widetilde x^{(2)}_q&=\overline u_{q},\quad q\in \partial D,\\
 \widetilde  w^{(2)}_{qp} (\widetilde x^{(2)}_{q} - \widetilde x^{(2)}_{p})
  &=\overline v_{q}, \quad
 q  \in  \partial D,
 p \in \mathcal N(q),
 \end{aligned}\right.
 \end{equation}
which coincide with system \eqref{EQ1}, \eqref{EQ2} and \eqref{Eq3} with $\gamma_{\{p,q\}}$ replaced by $\widetilde w^{(2)}_{p,q}$.
 Since the given Dirichlet data $\{\overline{\bf u}(k)\}_{k=1}^m$ contains a basis of the vector space
generated by the last $3n$ coordinates in
$\mathbb R^{4n}$,  this determines the
last $3n$ columns of the matrix
$\Lambda_{\bm\gamma}$ that represents the DtN map. By
\cite[Theorem 5.1]{curtis1991DNmap}, these columns contain all the parameters needed
to parameterize the manifold of the DtN matrices, and hence
to reconstruct the entire DtN matrix
$\Lambda_{\bm\gamma}$ and determine uniquely $\bm\gamma$. This
 implies necessarily that
$\widetilde w^{(2)}_{qp}=\gamma_{\{q, p\}}$ for any $p,q \in D$ if $\{p,q\}\in B$,  and from \eqref{cond2} it also follows that  $\widetilde w^{(2)}_{qq}= -\sum_{r\in\mathcal{N}(p)} {\gamma}_{\{q,p\}}$.
Finally, from the Neumann condition
$\widetilde w^{(2)}_{qp}=-\gamma_{\{q,p\}}$ for $q  \in  \partial D,
 p \in \mathcal N(q)$ and $\widetilde x_q^{(2)}=u_q$ for all $q\in\overline{D}$. This completes the proof of the theorem.
\end{proof}

When the Dirichlet data contain a basis for $\mathbb{R}^{4n}$, then for each DtN datum,
the weights of the optimal solution
must return the solution $\bf u$
for the neurons of the hidden layer that
correspond to $q \in D$.
Therefore, they realize the linear map from the Dirichlet data to the solution $\bf u$
in $D$ for all vectors of a basis. Hence
we have the following corollary.
\begin{corollary} \label{corollary}
Let the conditions in Theorem \ref{Thm1} hold with
$m\geq 4n$ such that \({\bf \overline{u}}(k)_{k=1, \dots, m}\)
contains a basis of $\mathbb R^{4n}$.
Then there is a unique minimum ${\widetilde W}_*=(\widetilde W^{(1)}_*,\widetilde W^{(2)}_*)$
of the loss  \eqref{loss}
such that $C_\alpha({\widetilde W}_*)=0$. Moreover,
$\widetilde W_*^{(2)}$ satisfies \eqref{EQ10},
and the matrix $\widetilde W_*^{(1)}$ is such that
$\widetilde W^{(1)}_{*,pq}$ for $p \in \partial D$
and $q \in D$ are the entries of the matrix that
maps the Dirichlet datum to the solution $\bf u$
in $D$.
\end{corollary}

\begin{remark}\label{rmk:interp2}
Theorem \ref{Thm1} and Corollary \ref{corollary} provide the interpretability of the proposed NN approach in the following sense. Since the solution to the inverse problem can be directly read from the NN's weights, and that there exists a unique minimum corresponding to the true conductivities, Theorem \ref{Thm1} contributes to understanding the ``Rashomon set"—the set of nearly equivalent models that perform well \cite{rudin2022interpretable}. In this context, the proposed approach identifies a unique interpretable model within this set, enhancing the trustworthiness and robustness of the solution.
\end{remark}

\section{Numerical experiments and discussions}\label{NUME}
In this section, we present numerical results to illustrate the NN approach, including a comparison with the Curtis-Morrow algorithm \cite{curtis1991DNmap}.

\subsection{Data generation and training}
Recall the notation \eqref{Domain} for
the domain and \eqref{boundary_points}  for  the boundary points.
We generate simulated data following the description in
Section \ref{CondProb}. For any network size $n$, the exact conductivity distribution $\bm\gamma$
is taken to be
\begin{eqnarray}
    \begin{cases}
\bm\gamma_{\{p,q\}} = (2\sin^2(0.2i_q+0.4j_q)+1)/3\quad
\text{ if }\{p,q\} \text{ is horizontal}, \\
\bm\gamma_{\{p,q\}} = (2\sin^2(\tfrac{2}{11}i_q+\tfrac{4}{11}j_q)+1)/3 \quad \text{ if }\{p,q\} \text{ is
vertical},
\end{cases}
\end{eqnarray}
where the indices \( (i_q,j_q) \) refer to the node \( q \) adjacent to the edge \( \{p,q\} \), the right one for horizontal edges and the bottom one for vertical edges. The conductivity values fluctuate between 0.3 and 1, showing wave-like features. Unless otherwise stated, we take $n=10$ in the experiments. %See Fig. \ref{fig:cond} for a schematic illustration of a network of size $n=10$. %where the node with index $(i,j)$ is plotted with $(i,-j)$.
%\begin{figure}[hbt!]
%    \centering
%    \includegraphics[width=.32\linewidth]{immagini/conduct.png}
%    \caption{The reference conductivity with $n=10$.}    \label{fig:cond}
%\end{figure}
For the reconstruction using full DtN data $\Lambda$, we take \( m = 4n \) pairs of Dirichlet and Neumann data. For each $k\in \{1,\ldots,m\}$,
we set
\begin{equation}\label{order_boundary}
    q^{(k)}=\left\{\begin{aligned}
 (k,n+1), &\quad k=1,\ldots,n,\\
  (n+1,n+1-(k-n)),&\quad k=n+1,\ldots,2n,\\
  (n+1-(k-2n),0),&\quad k=2n+1,\ldots,3n,\\
  (0,k-3n),&\quad k=3n+1,\ldots,4n.
\end{aligned}\right.
\end{equation}
and take the $k$th Cartesian coordinate vector
${\bf \overline u}(k)=\mathbf e_{q^{(k)}}$
so that the vectors $\{{\bf \overline u}(k)\}_{k=1}^m$ form a basis
for the whole space $\mathbb{R}^{4n}$.
The Neumann data ${\bf \overline v}(k)$
is obtained by solving
problem \eqref{EQ1}--\eqref{EQ2}, and then computing the flux at the boundary nodes according to \eqref{Eq3}. This effectively computes the DtN matrix $\Lambda_{\bf \gamma}.$ Numerically, we also experiment with incomplete / partial DtN data, which comprises a subset of the DtN matrix $\Lambda$. Then there are fewer samples available for training the NN.

To minimize the loss $C_\alpha(\widetilde W)$, we adopt the Adam algorithm \cite{Kingma2014AdamAM}, which is the most popular training algorithm in deep learning. It dynamically adjusts the learning rates for each component, utilizing the first and second moments of the gradients.
At each iteration, it requires computing the gradients of the loss $C_\alpha$, updating the first and second moments, and then adjusting the weights. It has three adjustable hyperparameters: stepsize and two decay rates for moment estimates. We have adopted a decaying learning rate schedule. During the optimization process, the loss $C_\alpha(\widetilde W)$ tends to stagnate after a number of iterations,
and does not decrease further; then we reduce the step size by a factor of ten with an initial step size 0.002, and increase the inertia of the first and second moments. We have chosen the empirical decay rates for the inertia parameters, i.e., the first and second moments, within the ranges $\{0.9, 0.99, 0.999\}$ and $\{0.999, 0.9999\}$, respectively. In the experiment, we initialize the conductivity $\bm\gamma$ randomly from \(\mathbb{R}_+^{2n(n+1)}\). The python source codes for reproducing the experiments will be made available at the github link \url{https://github.com/cornind/Discrete-EIT.git}.

%In our experiment, to prevent our recovered conductivity from being valued unrealistically, we use parameter $\theta_l$ to encode $\gamma_l$, for $l\in B$, in the following form
%\begin{equation}\label{eqn:encoder}
%    \gamma_l = 0.5\arctan\left(\theta_l\right)+0.25\pi+0.1,
%\end{equation}
%which takes value in $(0.1,0.1+0.5\pi)$.
%\begin{figure}[hbt!]
%    \centering
%    \includegraphics[width=0.4\linewidth]{immagini/desmos-graph (1).png}
%    \caption{Encode Function}    \label{fig:enter-label}
%\end{figure}
\subsection{Numerical results and discussions}

\subsubsection{Neural network approach}
\begin{figure}[hbt!]
    \centering
\begin{tabular}{cc}
\includegraphics[width=0.45\linewidth]{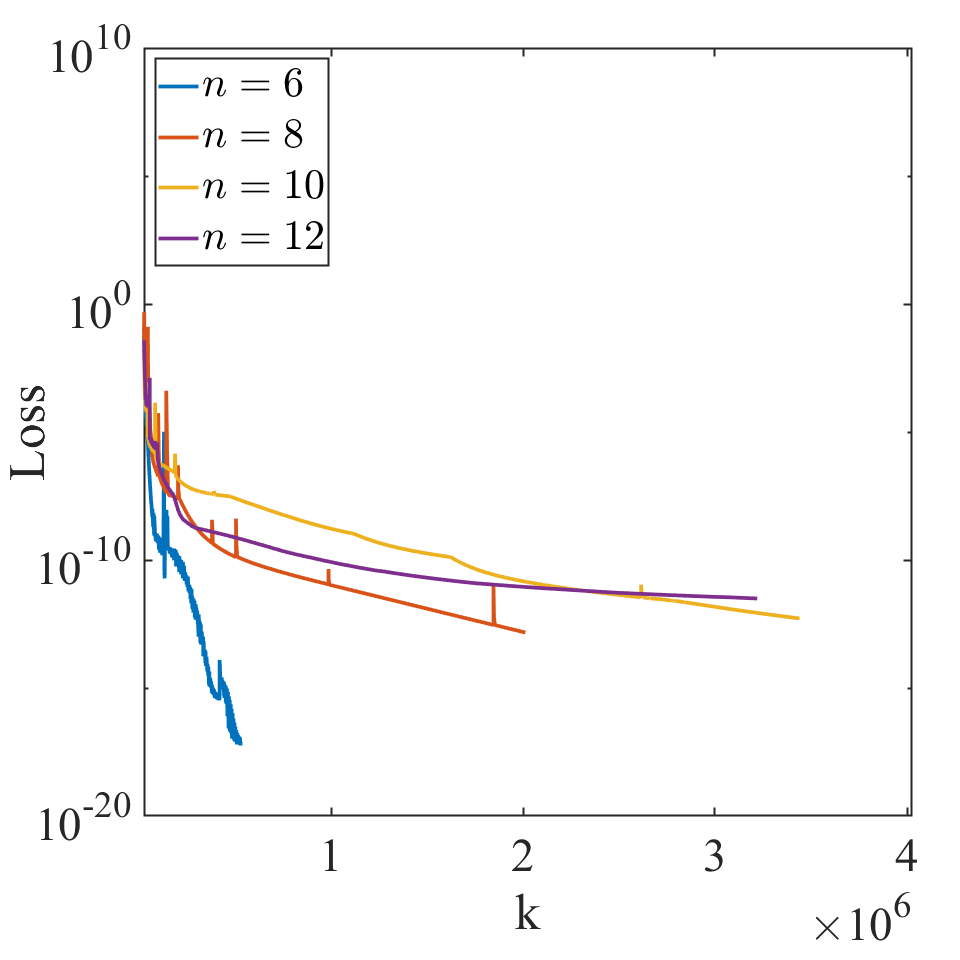}&    \includegraphics[width=0.45\linewidth]{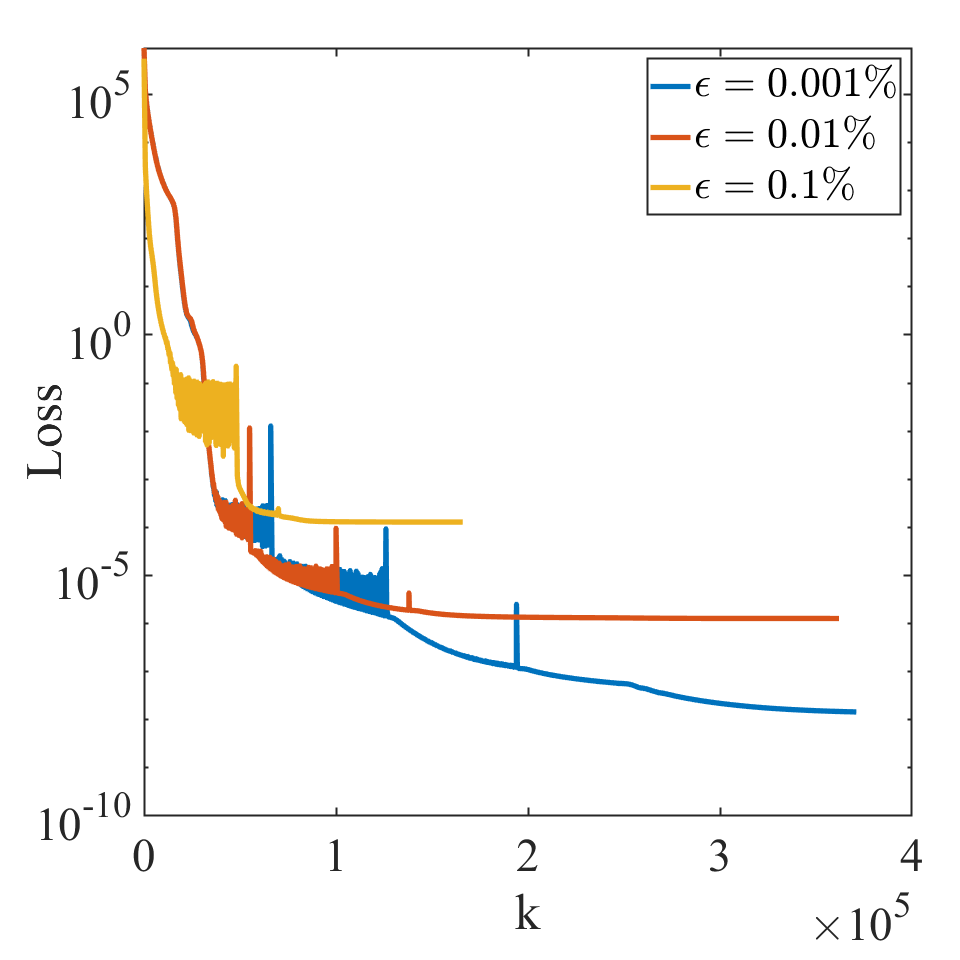}\\
(a) exact data & (b) noisy data
\end{tabular}
\caption{The training dynamics of the algorithm for (a) exact data with $n=6,8,10,12$, $\alpha=1$, random initial conductivity; and (b) noisy data at different noise level with $n=10$, $\alpha=1$. It takes 7.29s, 13.5s, 18.8s and 27.6s for per $10^4$ iterations when $n=6,8,10$ and $12$ for full data, respectively. }    \label{fig:vals}
\end{figure}

First, we present numerical results for exact data with different network sizes $n$. The training dynamics of the NN, illustrated in Fig. \ref{fig:vals}(a), shows a steady decrease in the loss value \( C_\alpha \), ranging between \( 10^{-16} \) and \( 10^{-11} \) by the end of the training process, with only minor fluctuations. When the network size $n$ increases from 6 to 12, the convergence behavior of the algorithm changes slightly, due to the associated exponential ill-conditioning with the number of unknowns, and the algorithm becomes less efficient in decreasing the loss and takes longer to reach steady-state, which results in larger reconstruction errors. The numerical reconstructions for exact data at different sizes $n$ are shown in Fig. \ref{fig:rec-exact}.
Overall, for exact data, the accuracy of the recovered conductivity $\widehat{{\bm\gamma}}$ is very satisfactory, showing a high precision near the boundaries and a slightly lower accuracy in the central region. Thus, the inverse problem exhibits a depth-dependent sensitivity, like EIT \cite{GardeKnudsen:2017}. Additionally, reconstruction errors increase with the network size $n$, particularly in the center. Note that the discrete inverse conductivity problem with a much larger network size $n$, e.g., $n=20$, is very challenging to resolve, and the NN approach (and also the Curtis-Morrow algorithm) could not yield accurate approximations. In fact, when $n$ is large, the computational noise introduced by floating-point arithmetic affects the data and can lead to disproportionately large errors in the recovered conductivities,  due to the severe instability of the problem.

\begin{figure}[hbt!]
\setlength{\tabcolsep}{0pt}
\begin{tabular}{cccc}
\includegraphics[width=.25\linewidth]{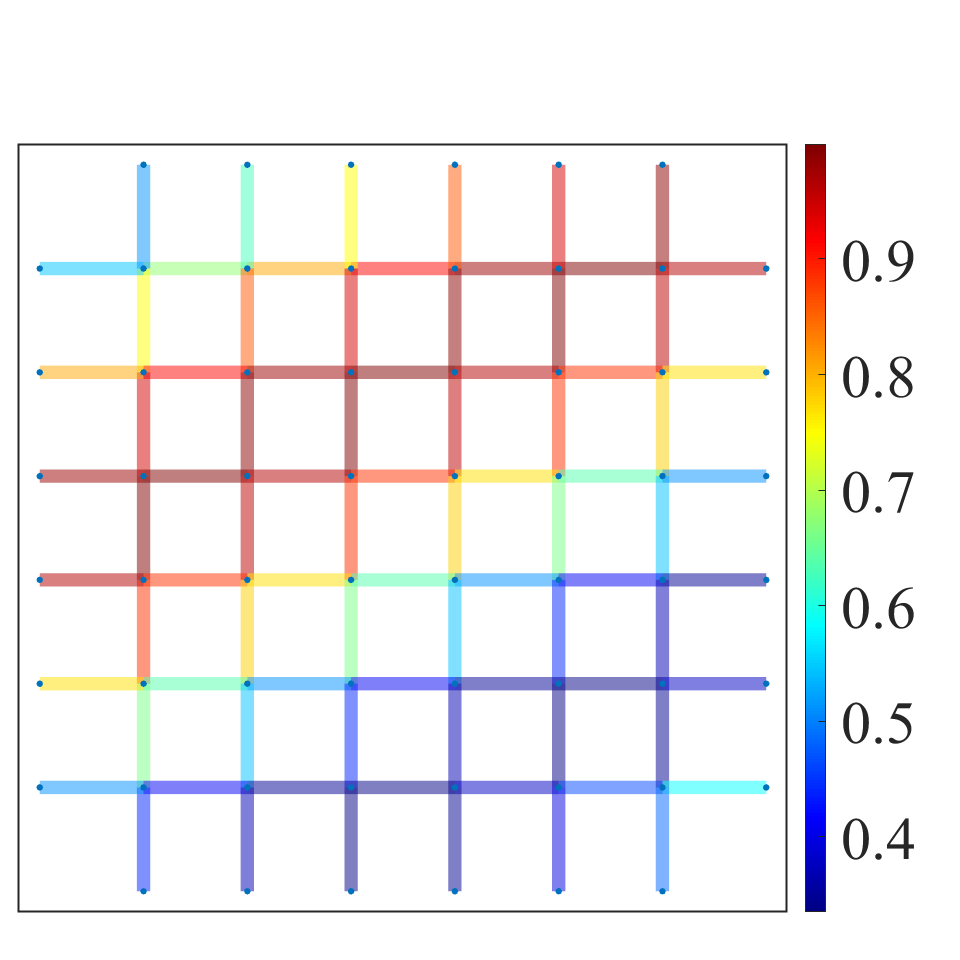} &\includegraphics[width=.25\linewidth]{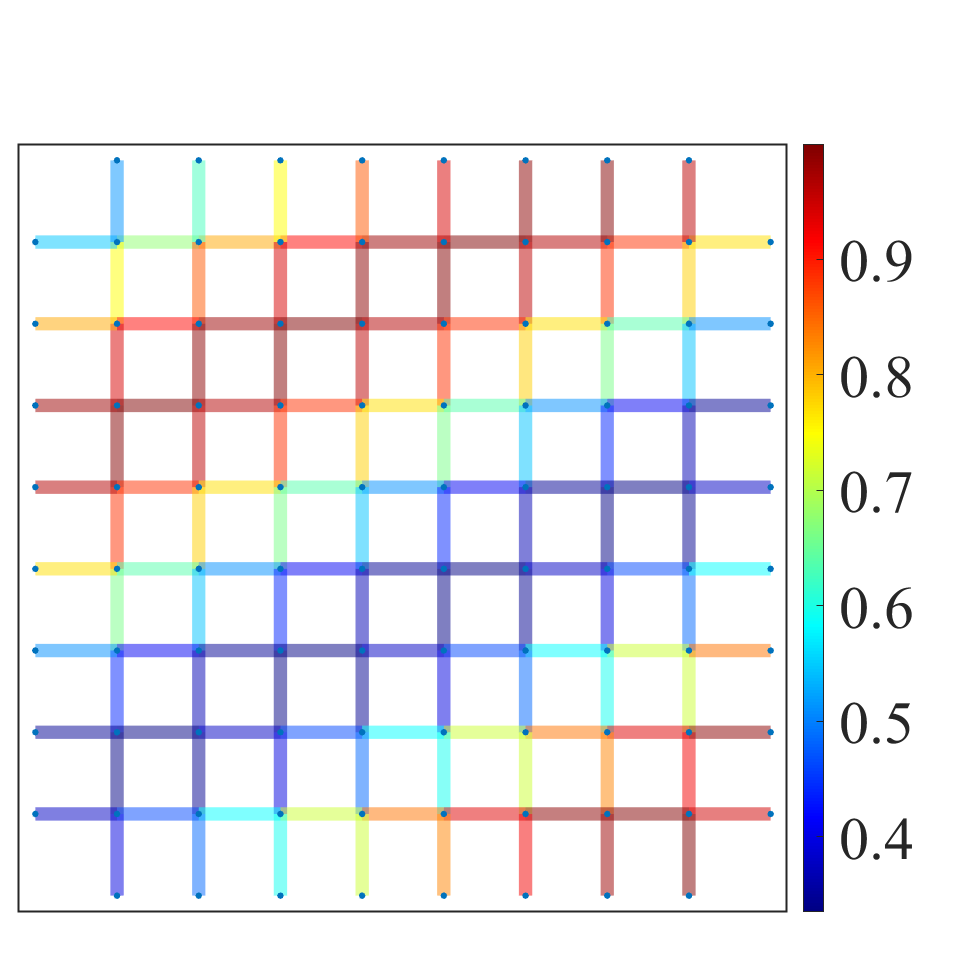} &\includegraphics[width=.25\linewidth]{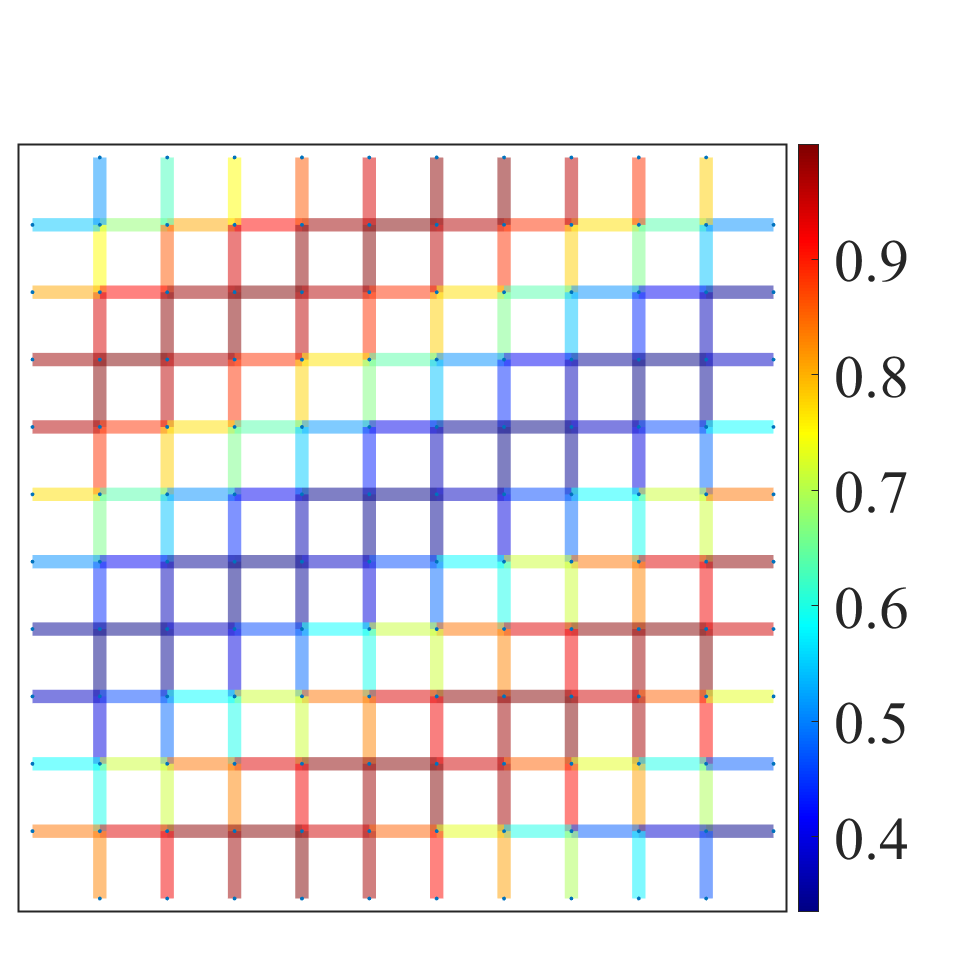} &\includegraphics[width=.25\linewidth]{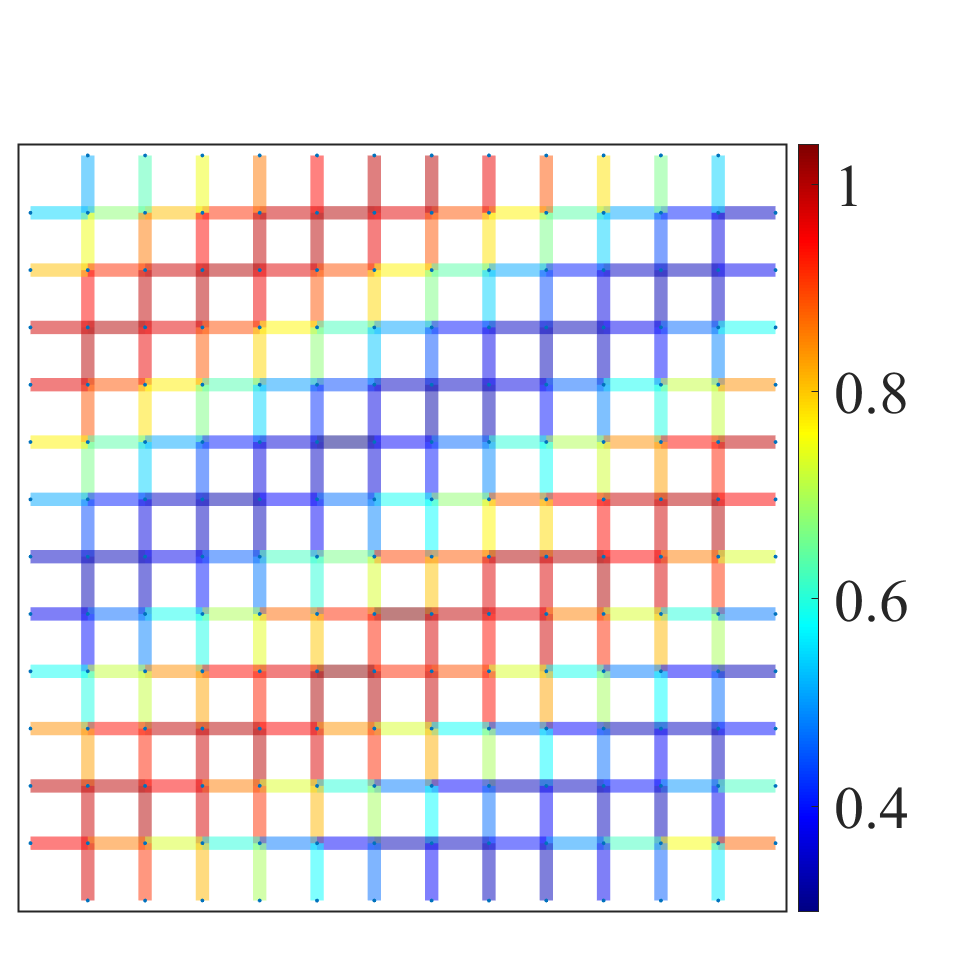} \\   \includegraphics[width=.25\linewidth]{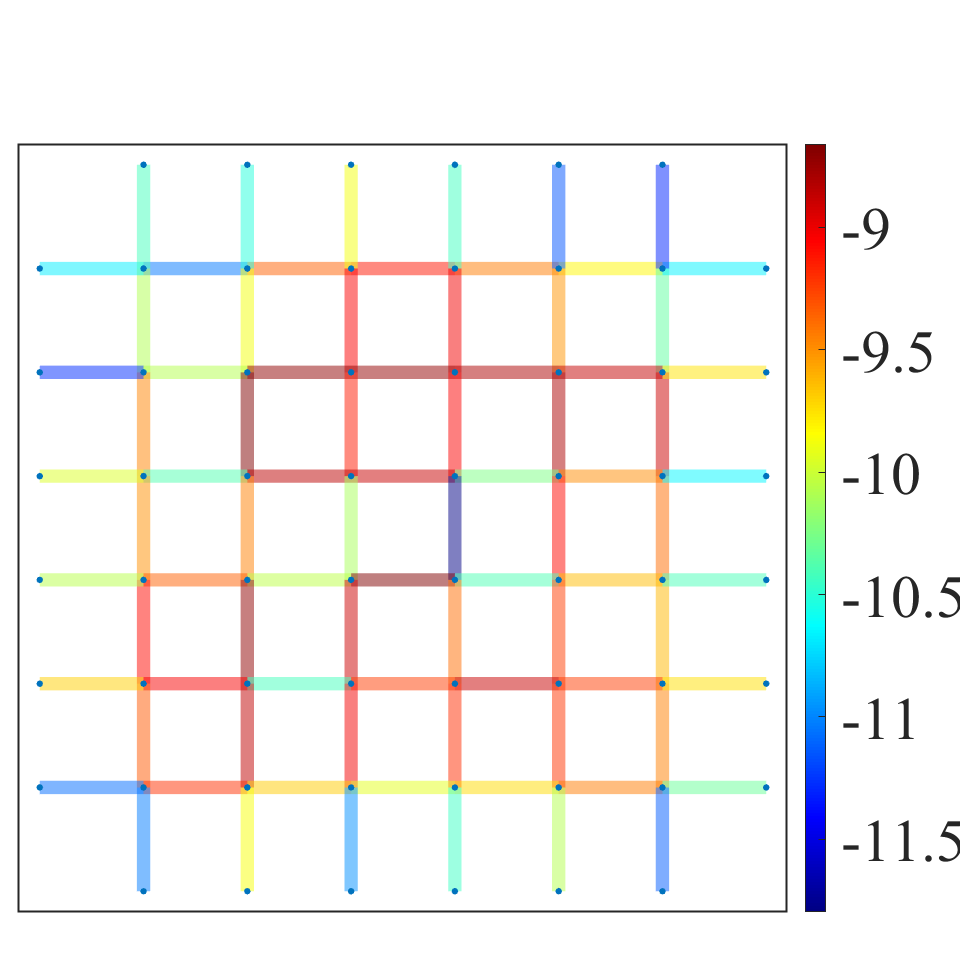} &\includegraphics[width=.25\linewidth]{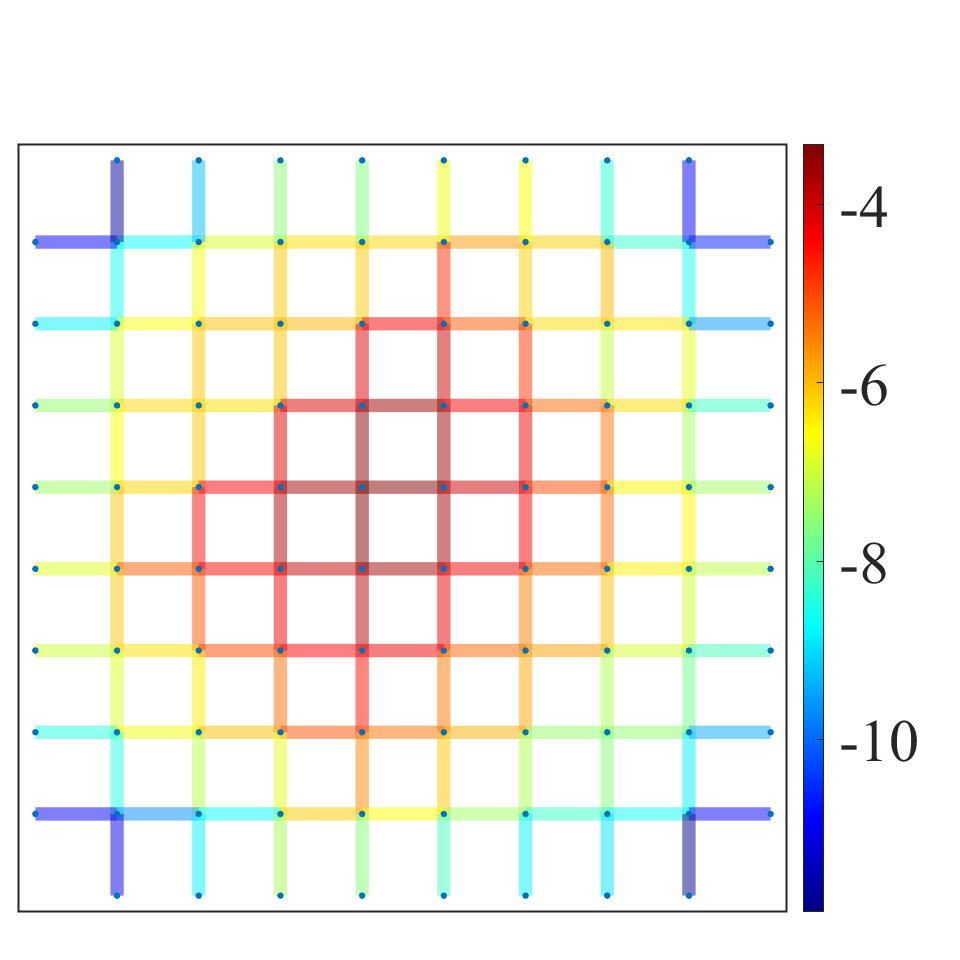} &\includegraphics[width=.25\linewidth]{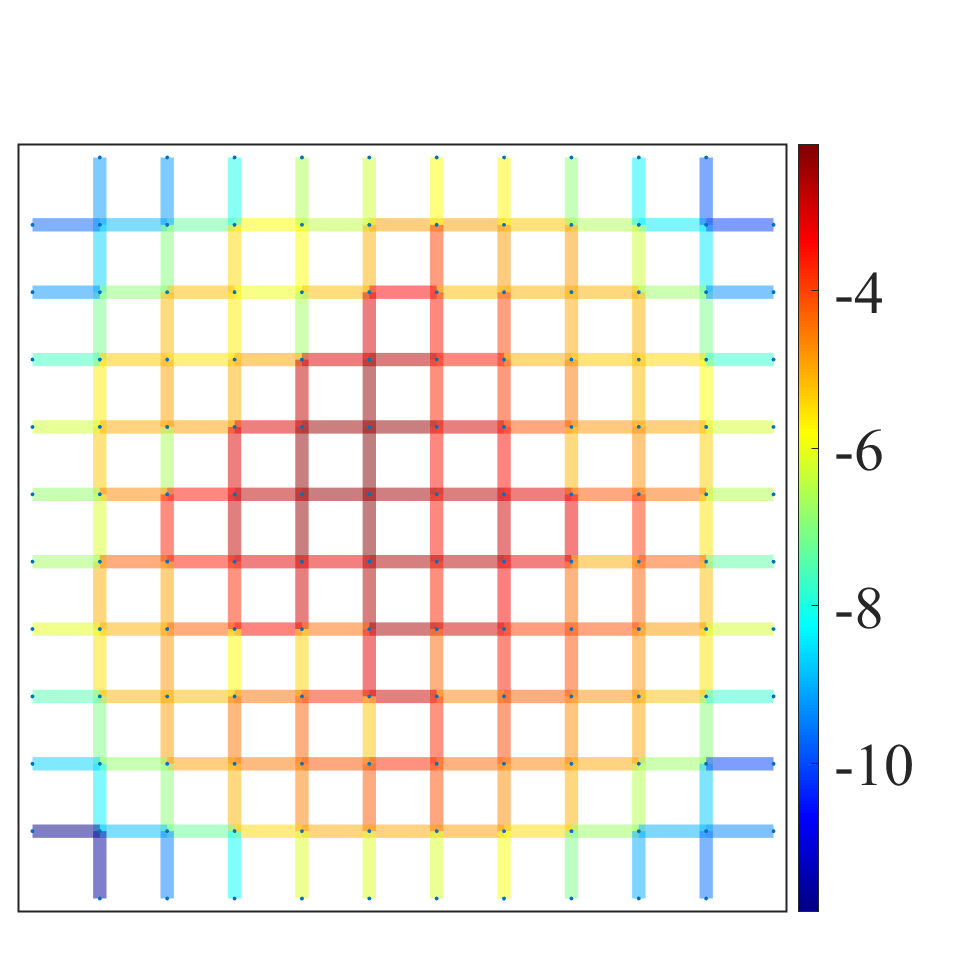}& \includegraphics[width=.25\linewidth]{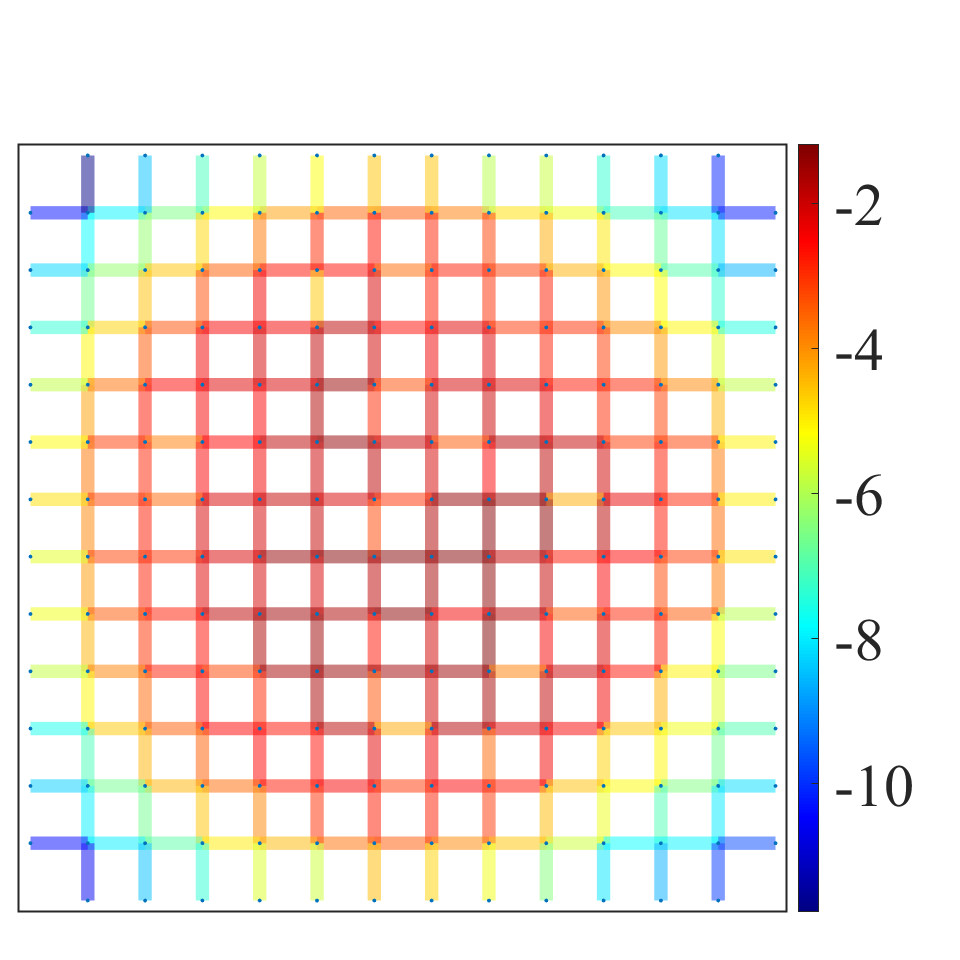}\\
(a) $n=6$&(b) $n=8$&(c) $n=10$&(c) $n=12$
\end{tabular}
\caption{The recovered conductivity for exact data in the case $n=6,8,10,12$ (top), and pointwise error $\log_{10}|\bm e_{\bm\gamma}|$ (bottom) }
\label{fig:rec-exact}
\end{figure}

Next, we illustrate the impact of noise in the data. We generate the noisy data $\Lambda_\epsilon$ (corresponding to the exact DtN matrix $\Lambda^\dag$) by
\begin{equation*}
(\Lambda_\epsilon)_{ij} = \Lambda_{ij} ^\dag+ \epsilon X_{ij}\Lambda_{ij} ^\dag,
\end{equation*}
where $\epsilon$ represents the relative noise level, and $X_{ij}$ follows the i.i.d. standard normal distribution. The training dynamics of the Adam algorithm for noisy data $\Lambda_\epsilon$ at various noise levels $\epsilon$ is shown in Fig. \ref{fig:vals}(b) for the setting with $n=10$ and $\alpha=1$. One observes that with the increase of the noise level $\epsilon$, the algorithm reaches a steady-state in much fewer iterations, i.e., quicker discovery of the minimum. However, the presence of noise also causes the optimizer to stagnate at much larger loss values, e.g., around $10^{-6}$ for the noise level $\epsilon=0.01\%$ (versus $10^{-11}$ for exact data), which also results in much larger errors of the recovered conductivity $\widehat{{\bm\gamma}}$ when compared with that for the exact data. The error associated with the recovered conductivity $\widehat{{\bm\gamma}}$ increases dramatically with the noise level $\epsilon$; see Fig. \ref{fig:rec-noisy} for the reconstructions and pointwise error. Nonetheless, for noisy data, the recovery near the boundary is still quite acceptable, and the dominant error occurs in the central region of the network, similar to the case of exact data.

\begin{figure}[hbt!]
    \centering
    \begin{tabular}{ccc}
    \includegraphics[width=0.3\linewidth]{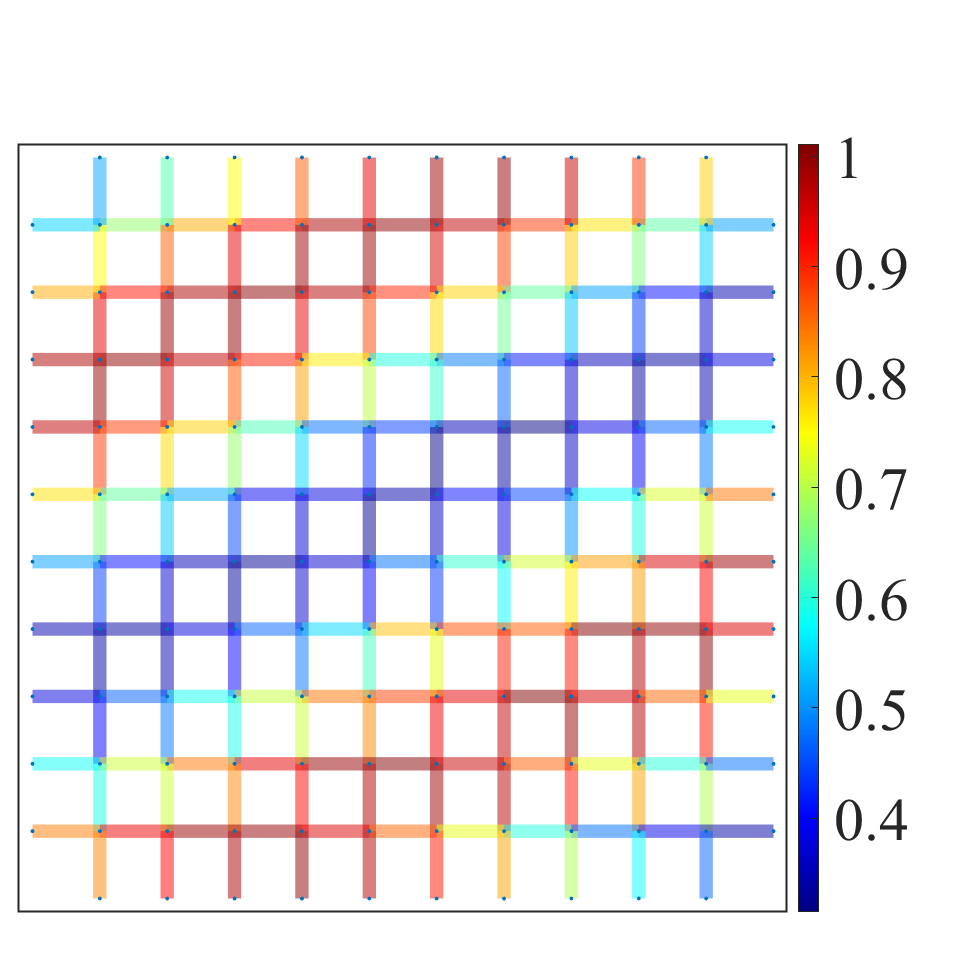}   &\includegraphics[width=0.3\linewidth]{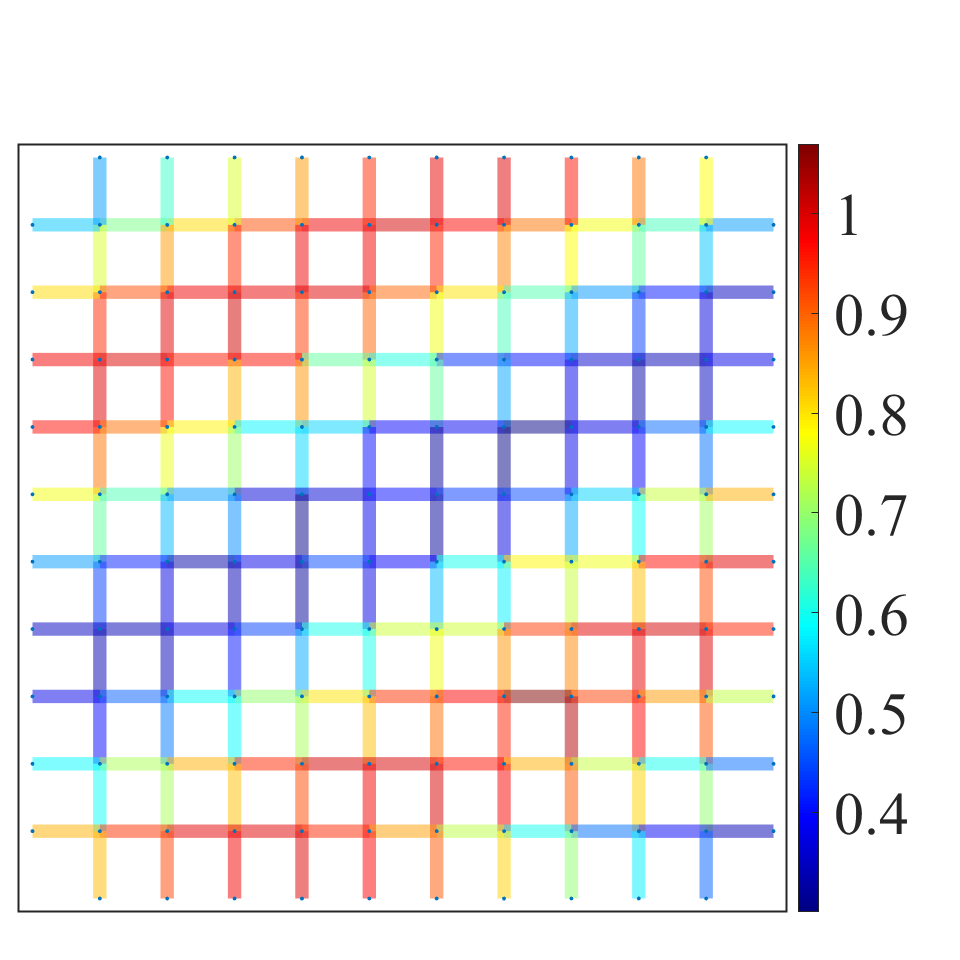} &  \includegraphics[width=0.3\linewidth]{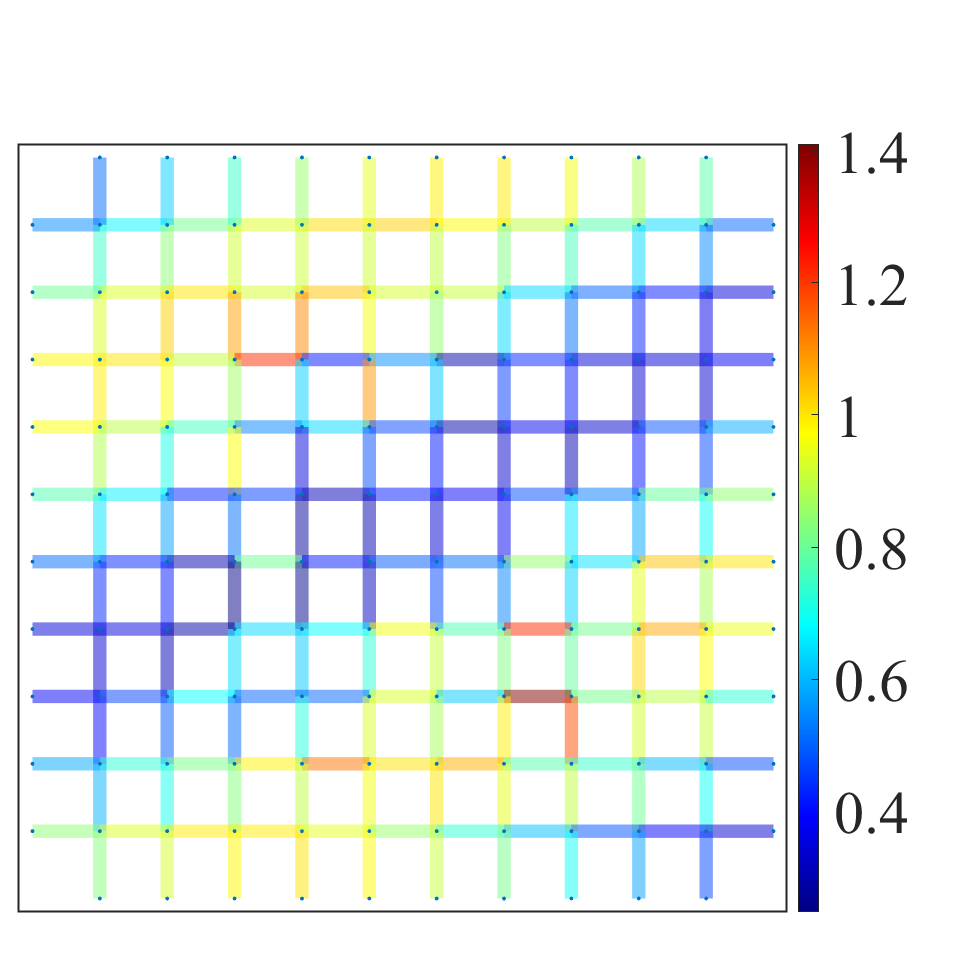} \\
  \includegraphics[width=0.3\linewidth]{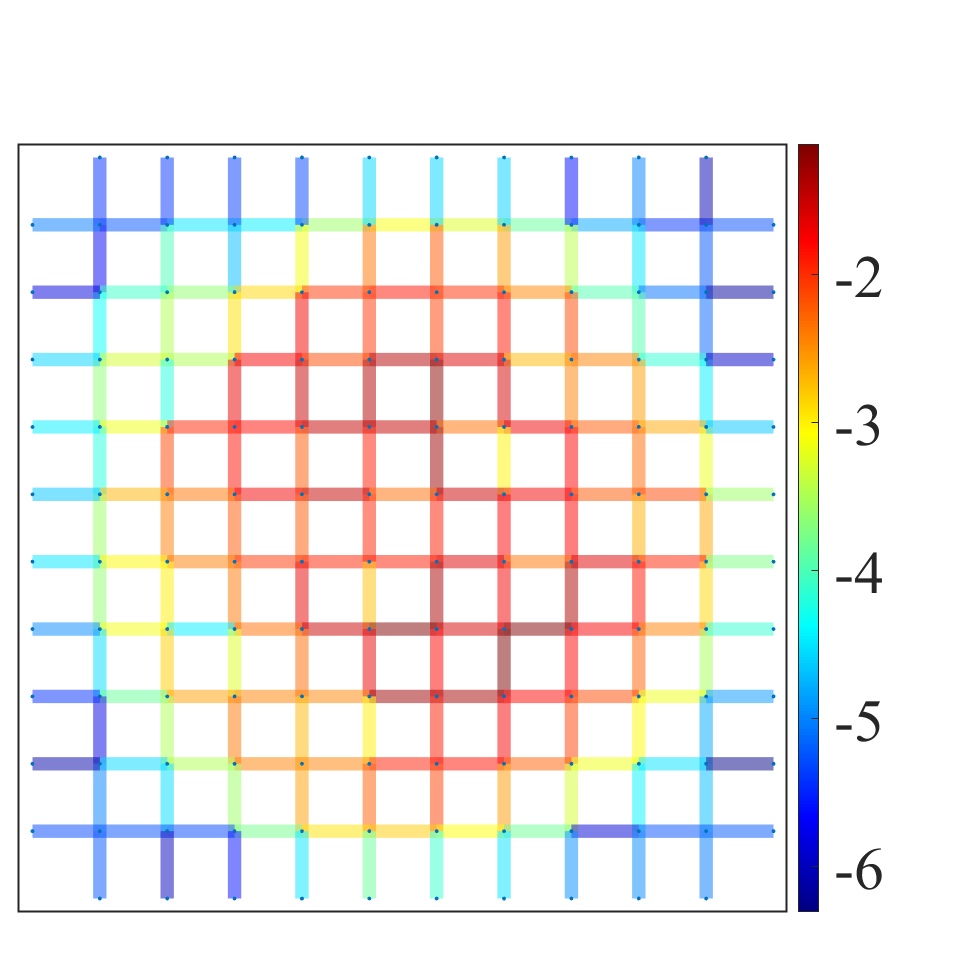}   &\includegraphics[width=0.3\linewidth]{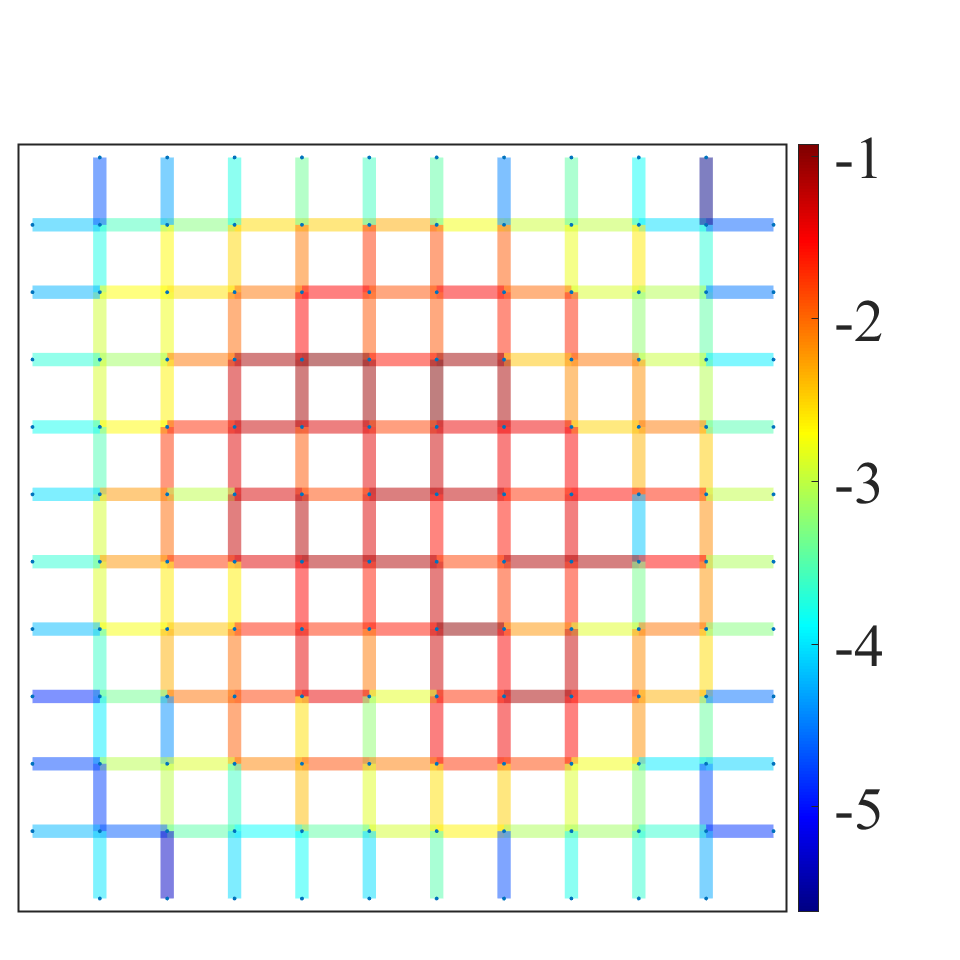} &  \includegraphics[width=0.3\linewidth]{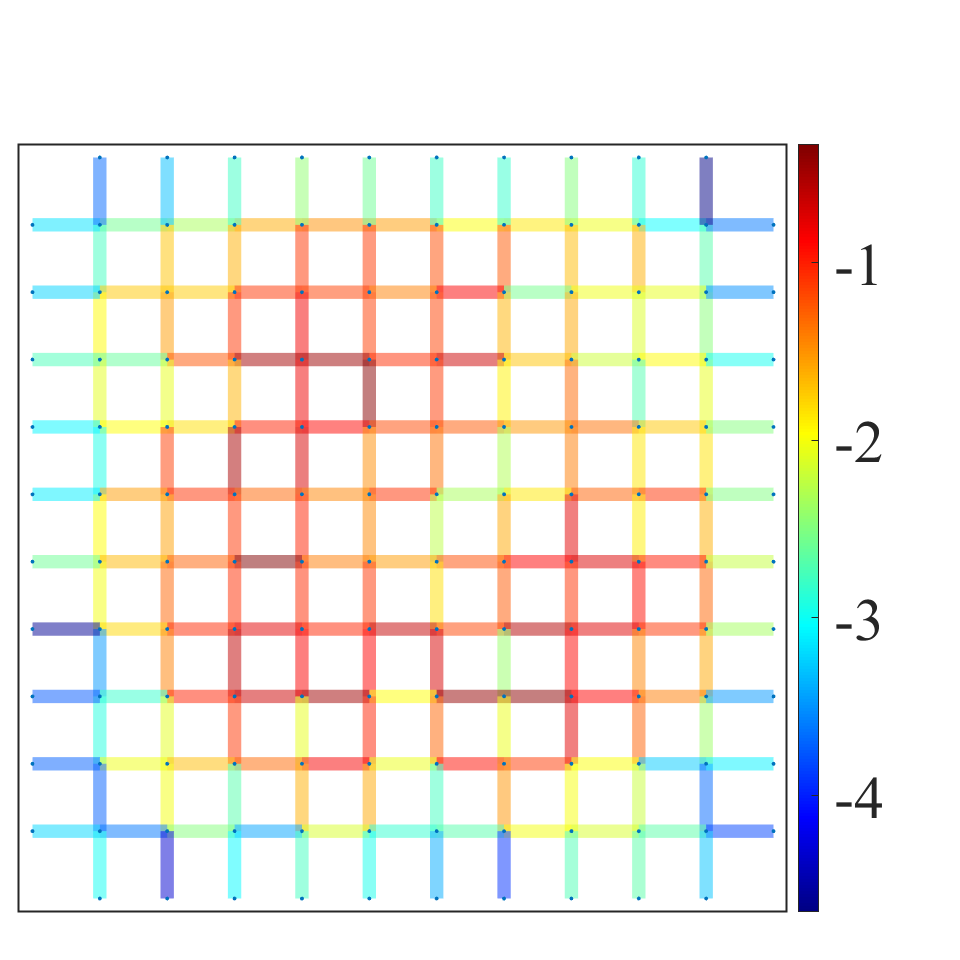} \\
    (a)  $\epsilon=0.001\%$ & (b) $\epsilon = 0.01\%$&(c)  $\epsilon = 0.1\%$

\end{tabular}
    \caption{The recovered conductivity $\widehat{\bm\gamma}$ (top), and the log error $\log_{10} |e_{\bm\gamma}|$ (bottom) at three noise levels. The maximum error $\|\bm e_{\bm\gamma}\|_\infty$ is 0.077, 0.12 and 0.47 for $\epsilon=0.001\%$, $0.01\%$ and $0.1\%$, respectively.}
    \label{fig:rec-noisy}
\end{figure}

One distinct feature of the NN approach is its flexibility in incorporating \textit{a priori} knowledge, e.g., known edge conductivity or box constraint. For example, if the conductivity of one edge is known, e.g., the edge between $(5,5)$ and $(5,4)$, we can fix the edge conductivity during the training process. Fig. \ref{fig:err-fixed-edge} (for $\epsilon=0.01\%$ noise) indicates that incorporating such prior knowledge can reduce the error on the specified edge and its immediate neighbors. However, it does not improve much the overall accuracy of the reconstruction in the central region.

\begin{figure}[hbt!]
\centering
\setlength{\tabcolsep}{0pt}
\begin{tabular}{cc}
    \includegraphics[width=0.35\linewidth]{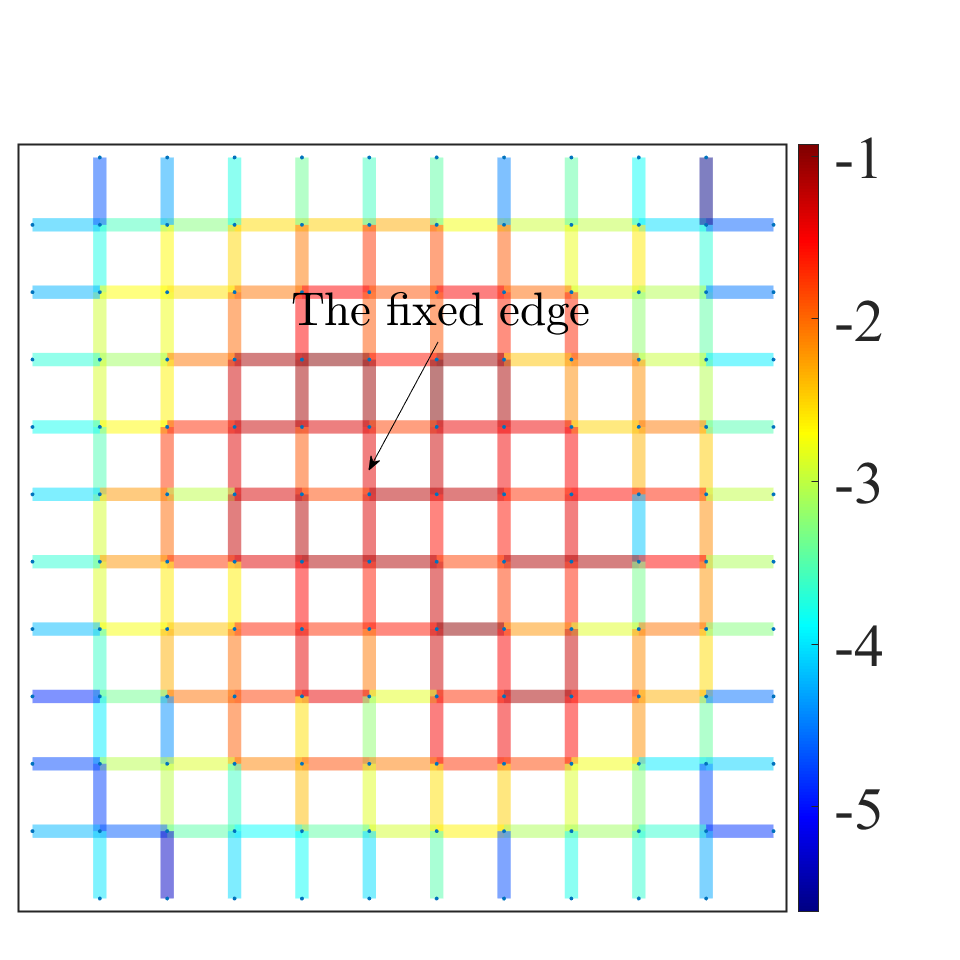} &
	    \includegraphics[width=0.35\linewidth]{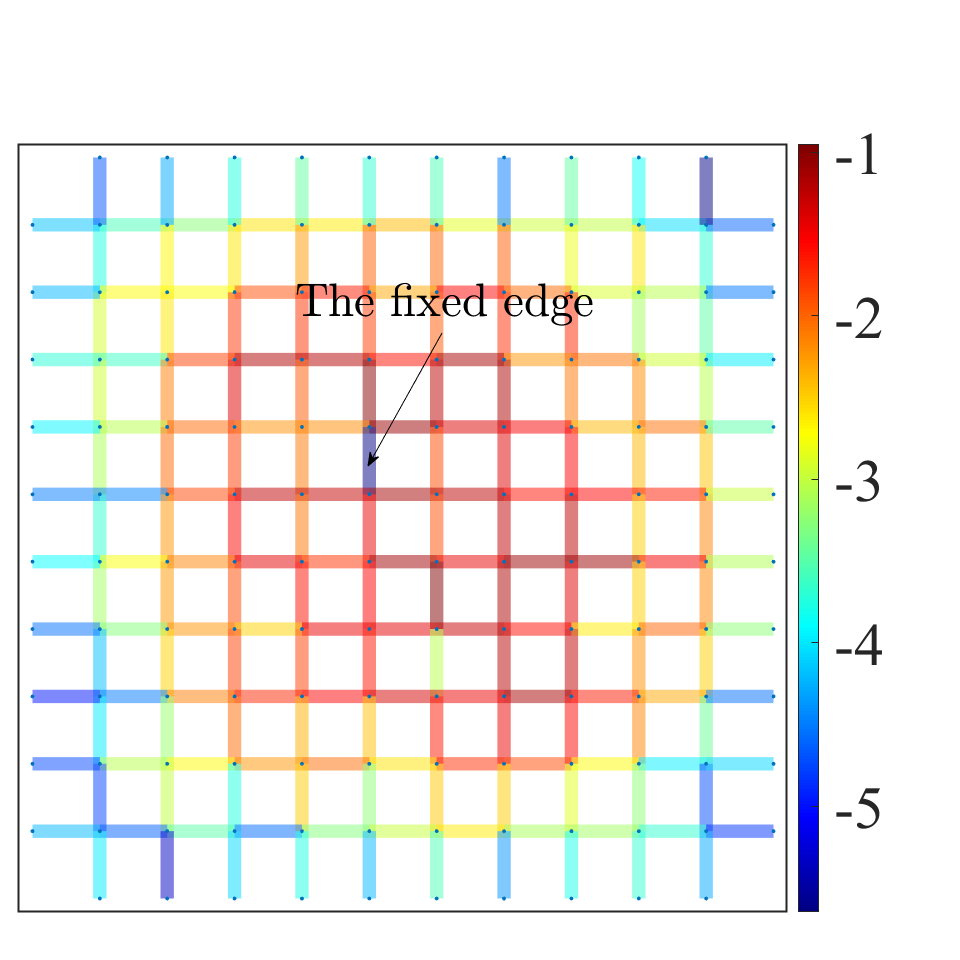}\\
     (a) free edge & (b) fixed edge
\end{tabular}
\caption{Numerical results for the conductivity with $\epsilon=0.01\%$. The maximum error is $0.12$ and $0.11$, respectively, for cases (a) and (b).} \label{fig:err-fixed-edge}
\end{figure}

Last, we illustrate the algorithm in the more challenging case with incomplete DtN data.
By Theorem \ref{Thm1}, the conductivity ${\bm\gamma}$ can be uniquely recovered using the last $3n$ columns of the DtN matrix $\Lambda_{\bm\gamma}$. Thus we numerically investigate the case, and present the results of the NN approach in Fig. \ref{fig:partialNN}. The first column of Fig. \ref{fig:partialNN} shows the convergence of the Adam algorithm in both full data and incomplete data (the last $3n$ columns of $\Lambda$). The plots show that the convergence of the Adam algorithm is largely comparable in both cases, except that the final loss value $C_\alpha(\widetilde{W}_*)$ is slightly smaller for incomplete data.  The recovered conductivity $\widehat{{\bm\gamma}}$ and pointwise error are also largely comparable with that for the full data for both exact and noisy data, indicating the robustness of the NN approach. Although not presented, we also tested the NN approach with partial data over the set $(n+1:4n)\times (n+1:4n)$. Even with exact data, the quarter region near one side of the network suffers from very large errors, due to a lack of data. Further experimental results with partial DtN data will be presented below.

\begin{figure}[hbt!]
\centering \setlength{\tabcolsep}{0pt}
\begin{tabular}{ccc}
\includegraphics[width=0.30\linewidth]{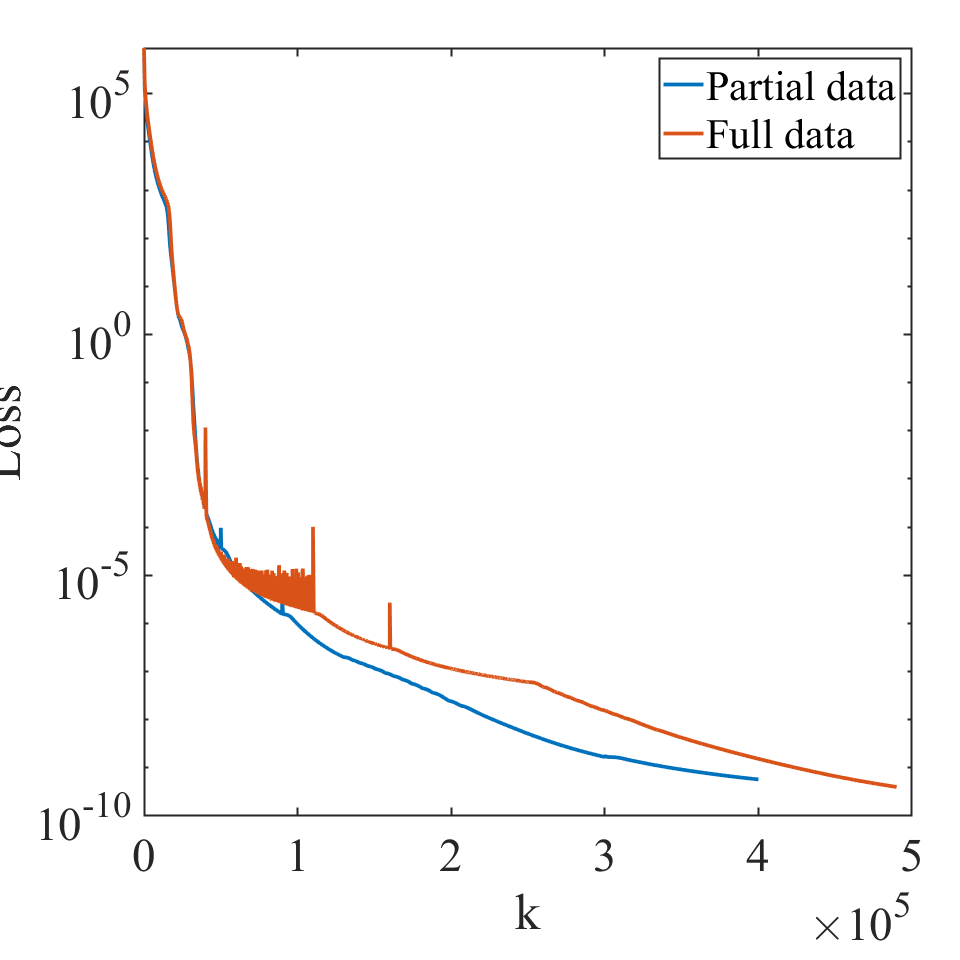}& \includegraphics[width=0.34\linewidth]{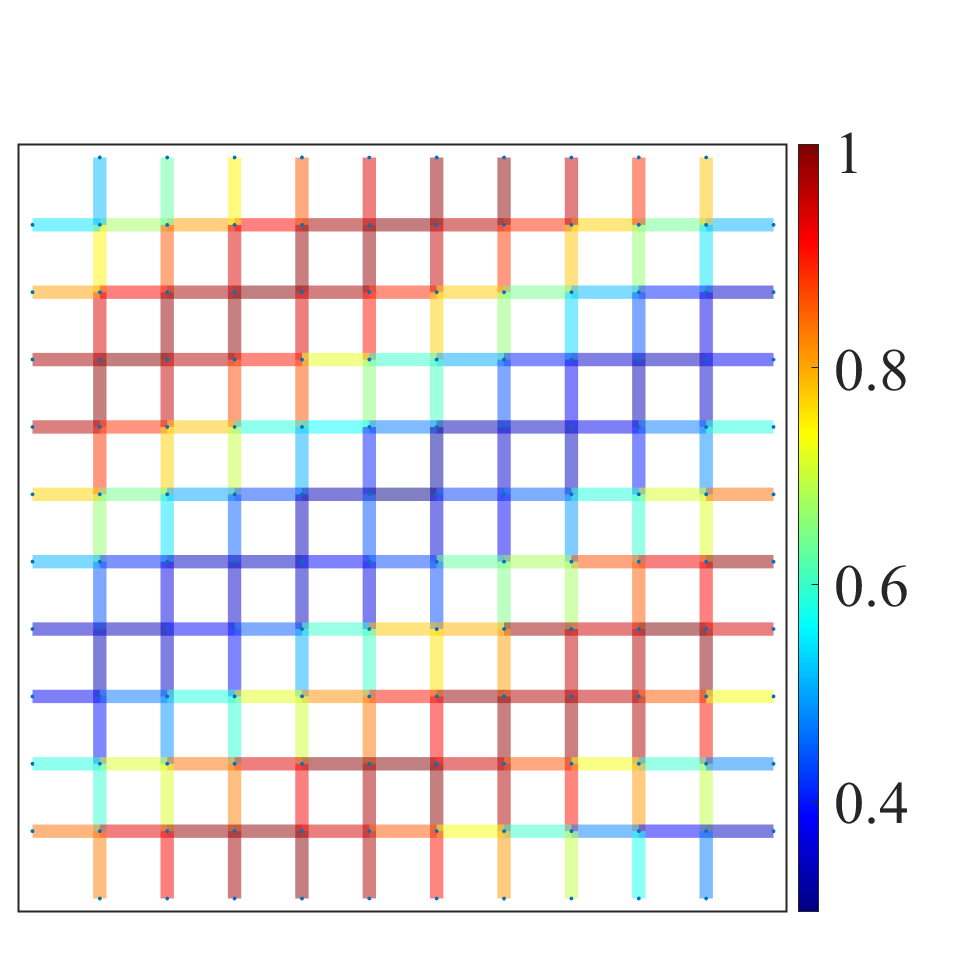}& \includegraphics[width=0.34\linewidth]{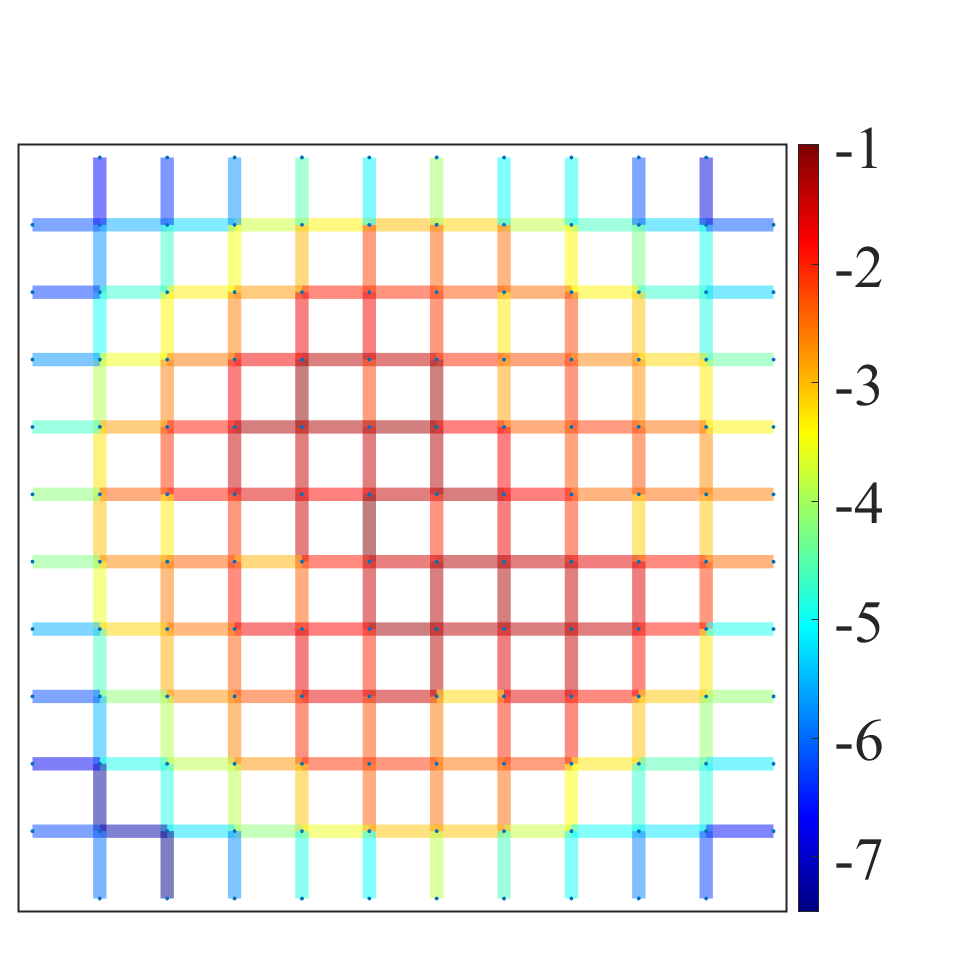}\\
\includegraphics[width=0.30\linewidth]{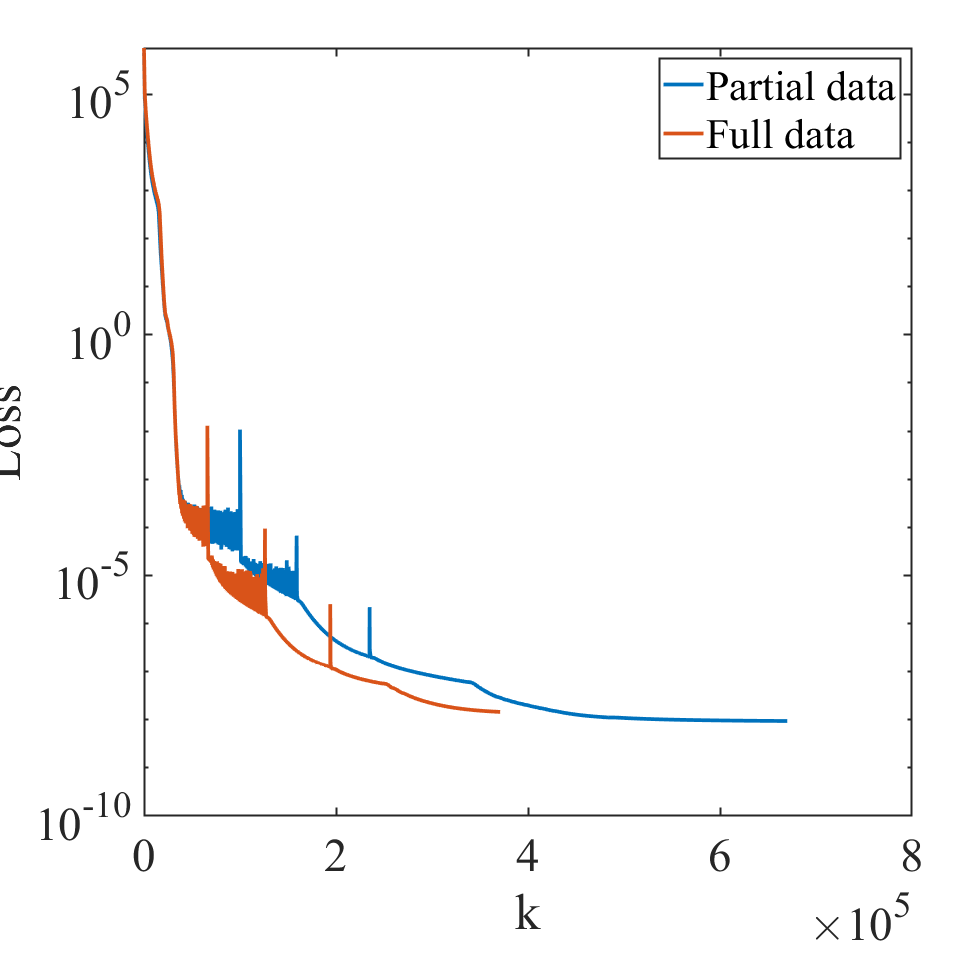}& \includegraphics[width=0.34\linewidth]{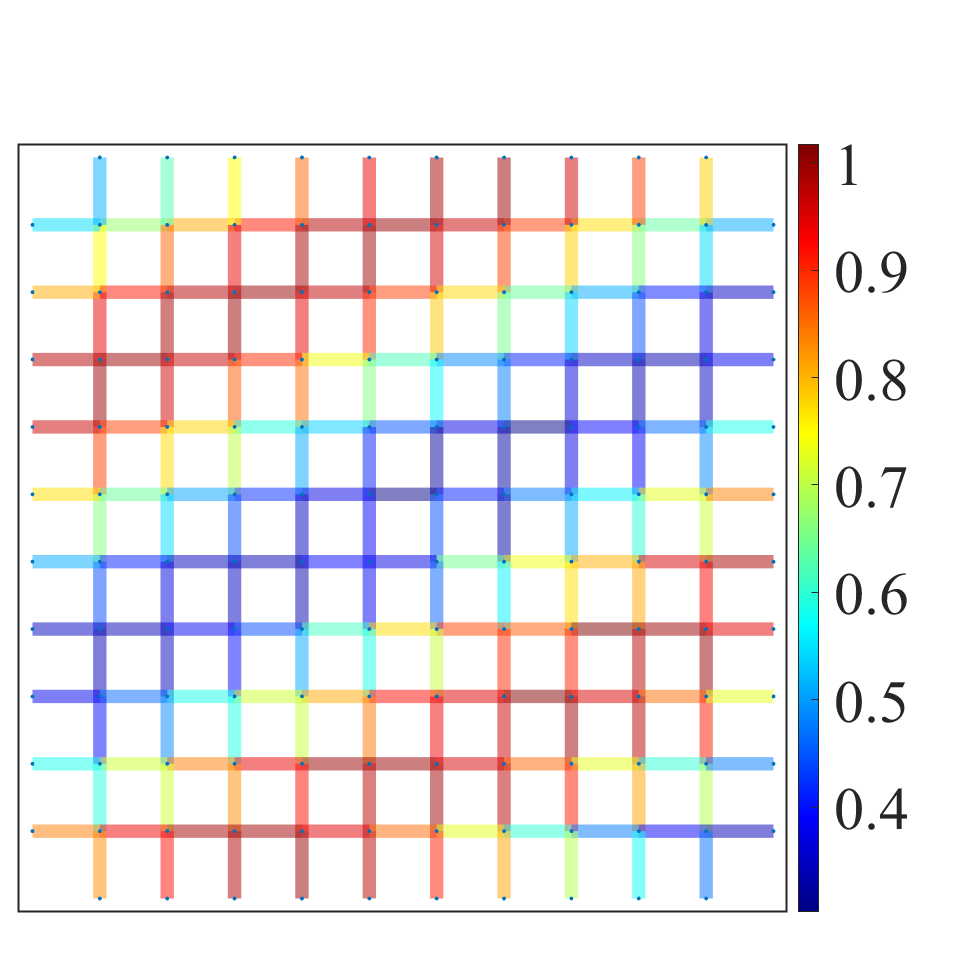}&    \includegraphics[width=0.34\linewidth]{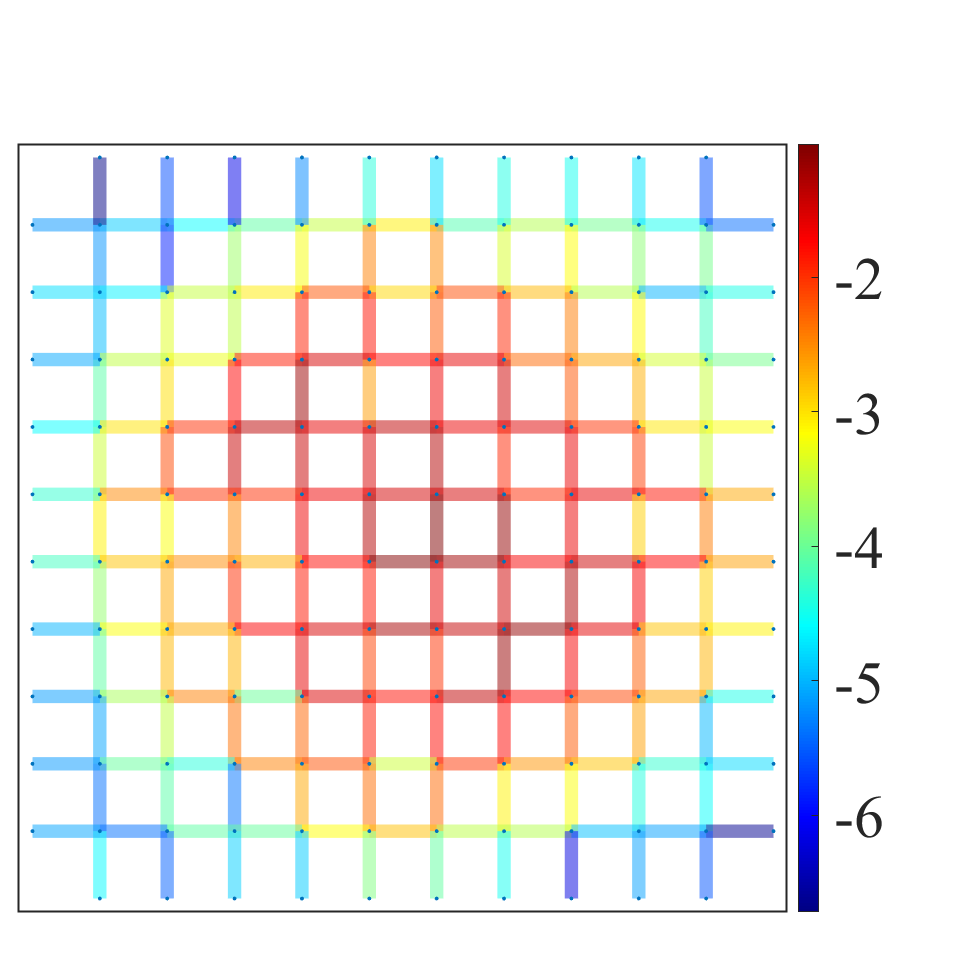}\\
 \includegraphics[width=0.30\linewidth]{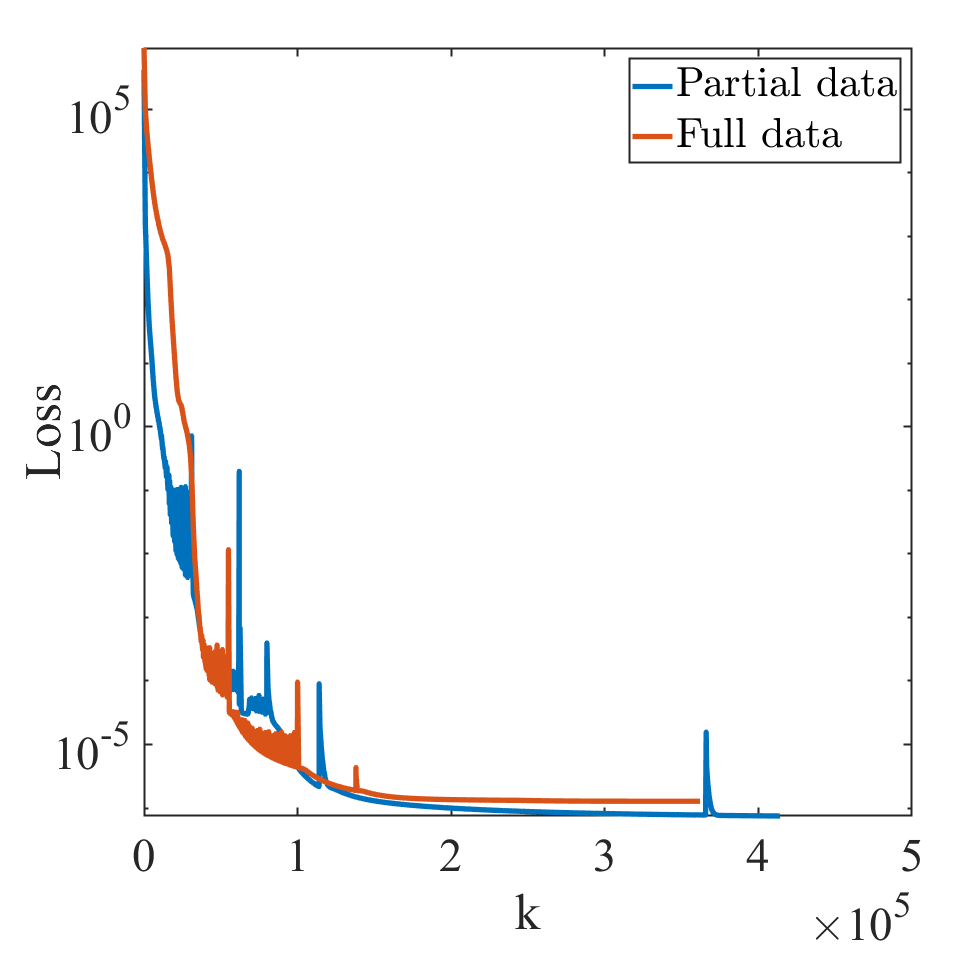}& \includegraphics[width=0.34\linewidth]{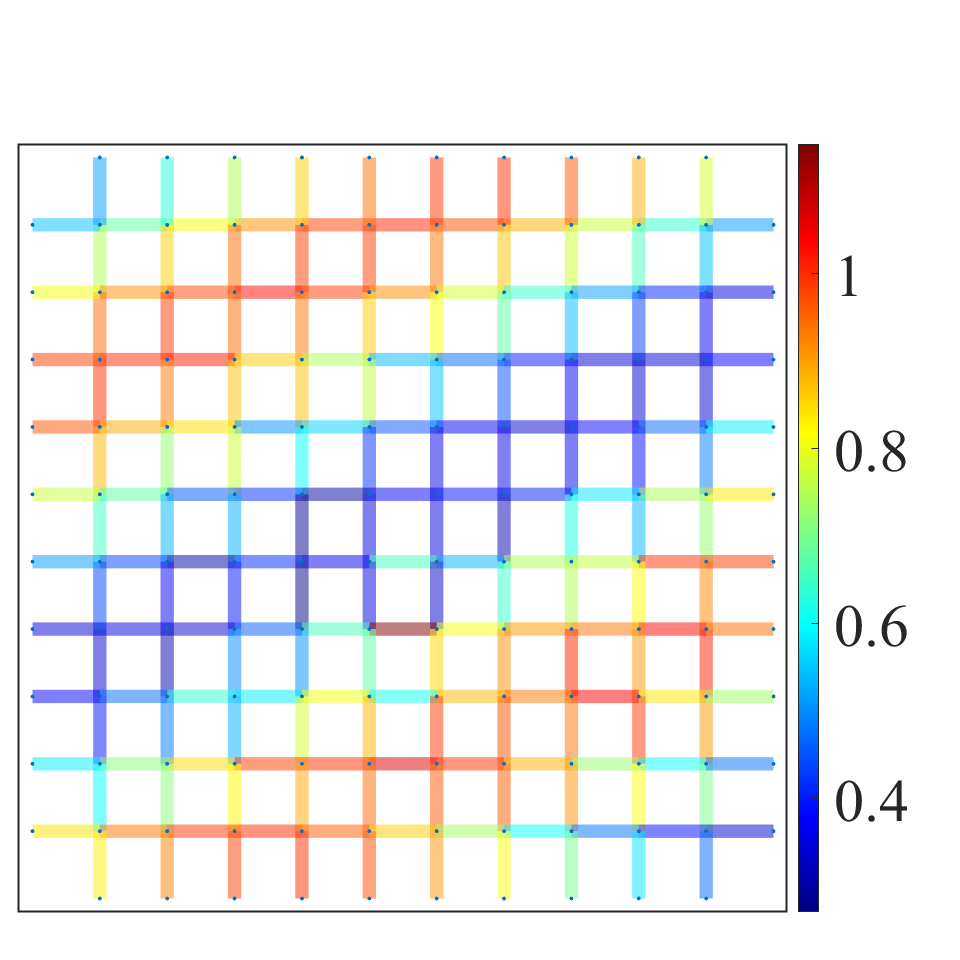}&    \includegraphics[width=0.34\linewidth]{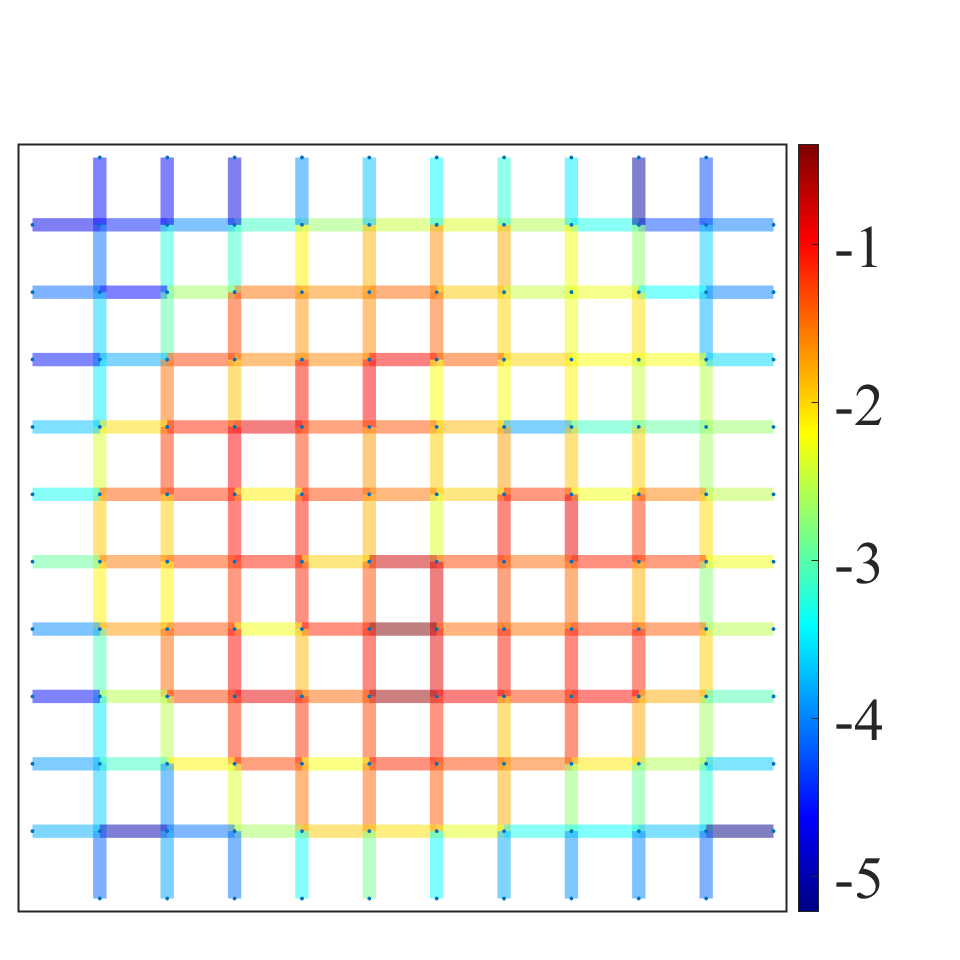}\\
 convergence      &  $\widehat{\bm \gamma}$ & $\log_{10} |e_{{\bm\gamma}}| $
\end{tabular}
\caption{The results of the NN approach using the last $3n$ columns of DtN map, at three noise levels $\epsilon=0.0001\%$ (top), $\epsilon=0.001\%$ (middle) and $\epsilon=0.01\%$ (bottom).} \label{fig:partialNN}
\end{figure}

\subsubsection{Comparison with Curtis-Morrow algorithm}

Now, we compare the NN approach with the Curtis-Morrow algorithm \cite{curtis1991DNmap}, which is an algebraic inversion technique and is implemented as follows. At each step, it involves solving for a specific boundary potential through a formally overdetermined system of linear equations $Ax = b$ (i.e., the matrix $A$ has more rows than columns). For exact data, the unique existence of a solution is guaranteed, and for noisy data, we obtain a solution by \(x = A^\dag b\), where $A^\dag$ denoting the Moore-Penrose pseudoinverse of $A$.

The numerical results in Fig. \ref{fig:CM} show that the Curtis-Morrow algorithm works fairly well for exact data, and can achieve high pointwise accuracy. However, it works poorly for noisy data. Indeed, the error in the central region blows up. It remains unknown how to properly regularize the algorithm \cite[p. 428]{Borcea:2013}. Intuitively, this can be viewed as a manifestation of the severe ill-posedness of EIT  \cite{M}. The comparison shows the high robustness of the NN approach in handling data noise, even though no explicit regularization is incorporated into the loss $C_\alpha(\widetilde{W})$. Such implicit regularizing effect is commonly observed for NN based inversion schemes \cite{JinLiQuanZhou:2024}.

\begin{figure}[hbt!]
\centering\setlength{\tabcolsep}{0pt}
\begin{tabular}{cc|cc}
\toprule
\multicolumn{2}{c}{Curtis-Morrow algorithm}&\multicolumn{2}{c}{NN approach} \\
 \midrule
 \includegraphics[width=0.25\linewidth]{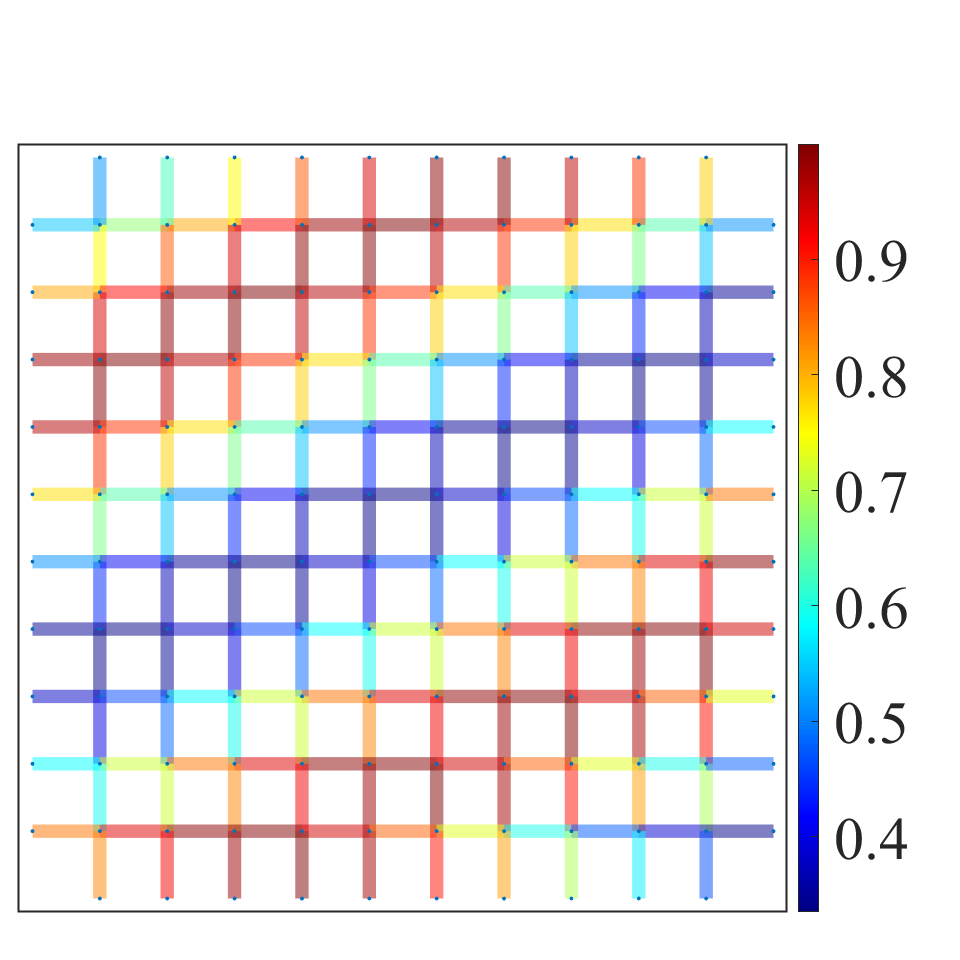} &\includegraphics[width=0.25\linewidth]{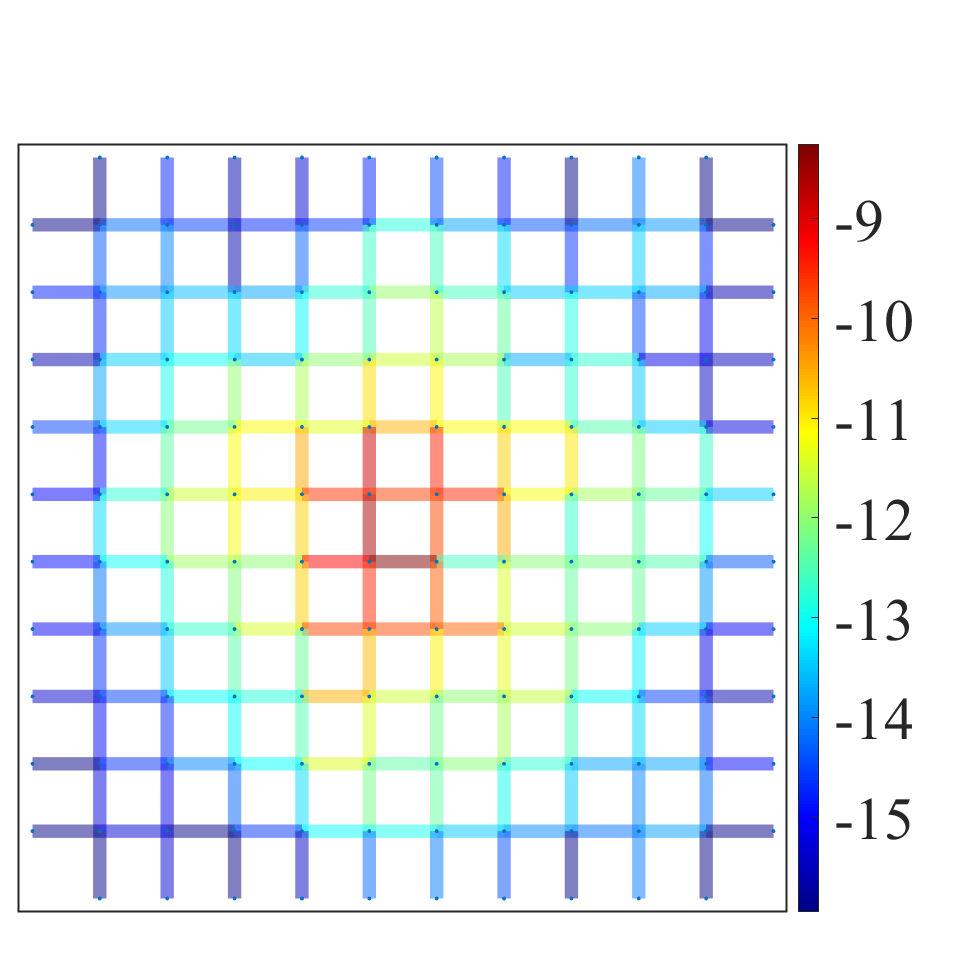}&\includegraphics[width=0.25\linewidth]{conduct_n10_noNoise.png} &\includegraphics[width=0.25\linewidth]{err_n10_noNoise.png}\\ 
 \includegraphics[width=0.25\linewidth]{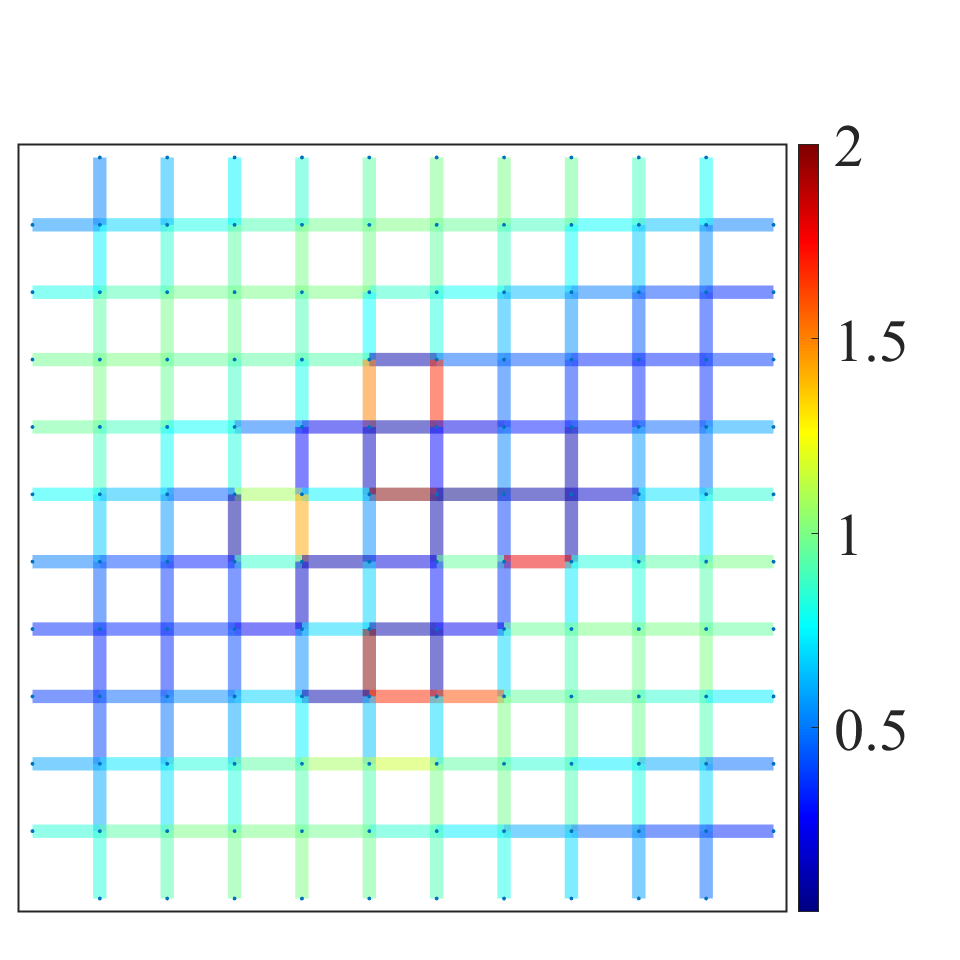}&\includegraphics[width=0.25\linewidth]{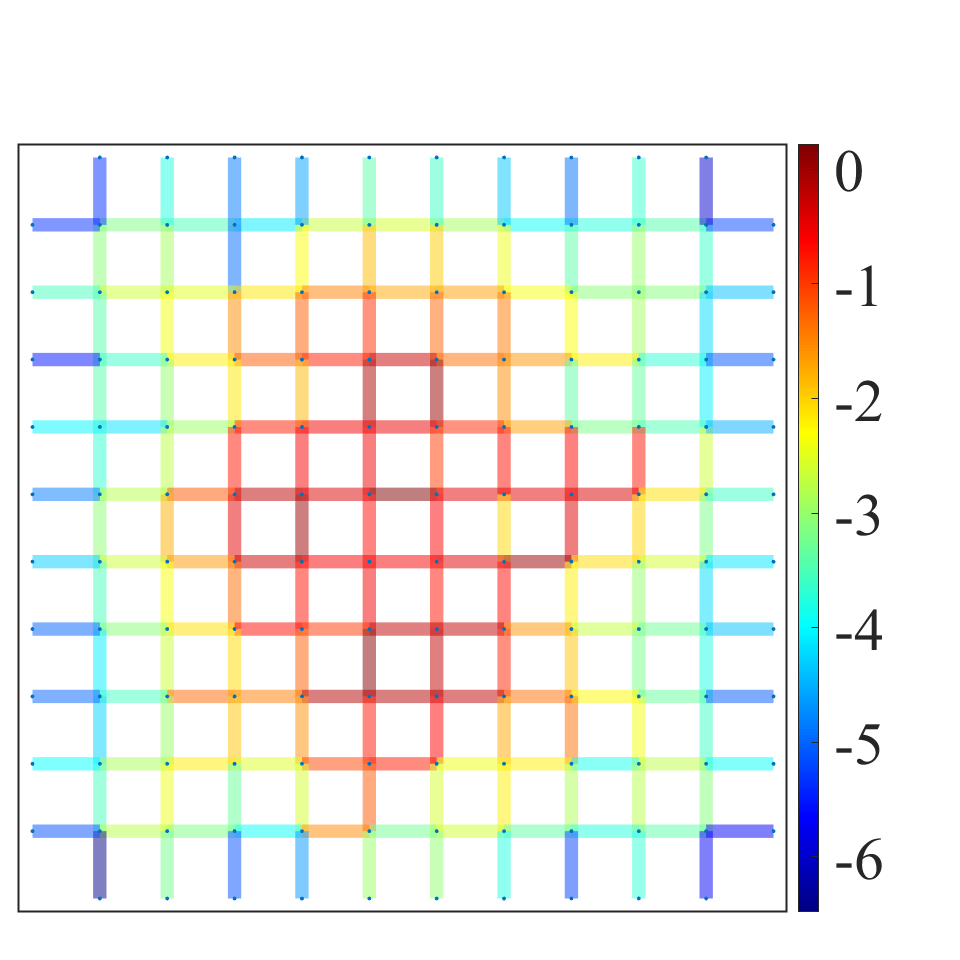}&\includegraphics[width=0.25\linewidth]{conduct_n10_Noise1e-5.png} &\includegraphics[width=0.25\linewidth]{err_n10_Noise1e-5.png}\\ 
 \includegraphics[width=0.25\linewidth]{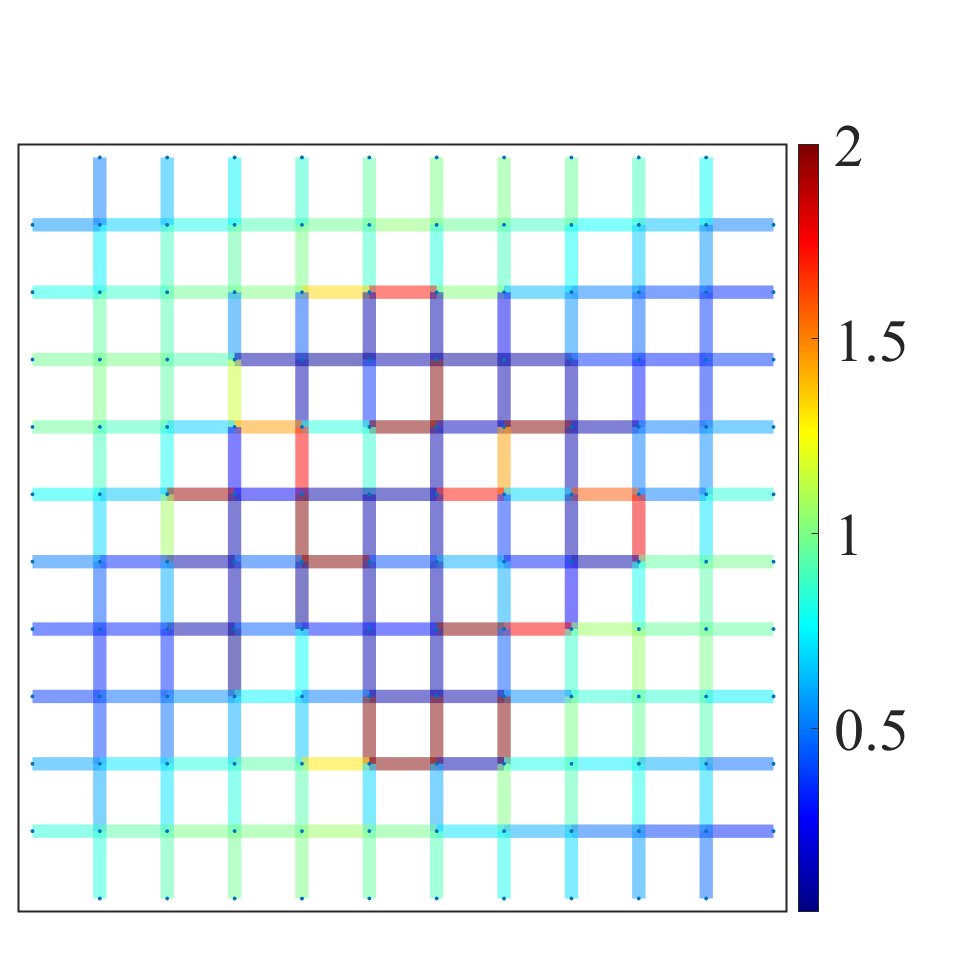} &\includegraphics[width=0.25\linewidth]{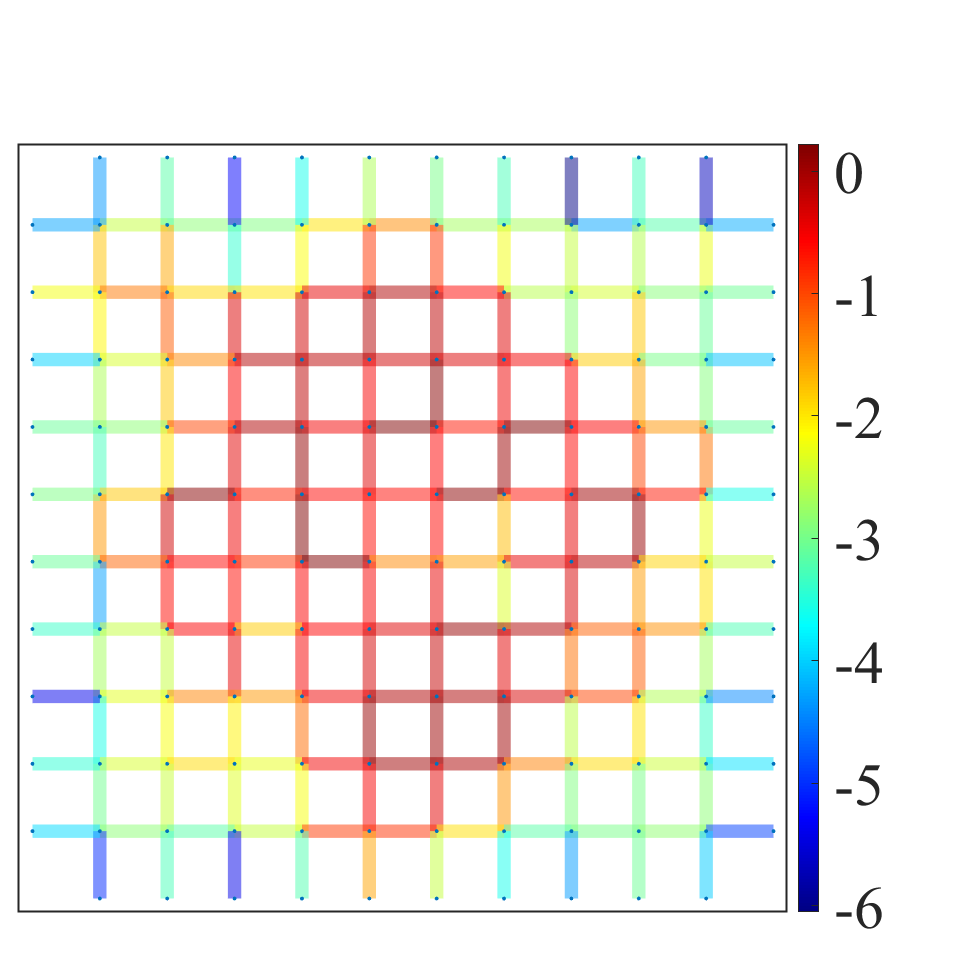}&\includegraphics[width=0.25\linewidth]{conduct_n10_Noise1e-4.png} &\includegraphics[width=0.25\linewidth]{err_n10_Noise1e-4.png} \\
    $\widehat {\bm\gamma}$ & $\log_{10}|e_{\bm\gamma}|$ &  $\widehat {\bm\gamma}$ & $\log_{10}|e_{\bm\gamma}|$

    \end{tabular}
    \caption{Numerical results by Curtis-Morrow algorithm (left) and the NN approach (right) at three noise levels with full data, $\epsilon=0\%$ (top), $\epsilon=0.001\%$ (middle) and $\epsilon=0.01\%$ (bottom).}
    \label{fig:CM}
\end{figure}

Next, we compare the NN approach with the Curtis-Morrow algorithm on the more challenging case with  fewer measurements.
We present numerical results for several incomplete / partial data sets: the last three columns ($C_{3n}$), the set $(1:3 n)\times(n+1:4 n)$ ($T_{\rm sq}$), and the set $(1:2n)\times(n+1:4 n)$ ($T_{\rm re}$). To employ the Curtis-Morrow algorithm, we first reconstruct the full DtN matrix $\Lambda$ from the incomplete one, for which we adopt the methodology outlined in \cite{curtis1991DNmap}.
Since the DtN matrix $\Lambda$ is symmetric, for the case $C_{3n}$, it is sufficient to reconstruct only the leading $n \times n$ block. The first column of Fig. \ref{fig:partialCM} shows the error between the recovered DtN matrix (using the last $3n$ columns) and the exact one at different noise levels, in the leading $n\times n$ block. The recovered block of the DtN matrix suffers from pronounced errors when $\epsilon$ is not very small. This error persists in the recovery $\widehat{{\bm\gamma}}$ by the Curtis-Morrow algorithm. Nonetheless, the recovery $\widehat{\bm\gamma}$ is still reasonably accurate for edges very close to $\partial D$, but fails to capture the wave shape of the exact conductivity profile, and the error in the central part is quite significant.

\begin{figure}[hbt!]
\centering\setlength{\tabcolsep}{0pt}
\begin{tabular}{ccc}
\includegraphics[width=0.33\linewidth]{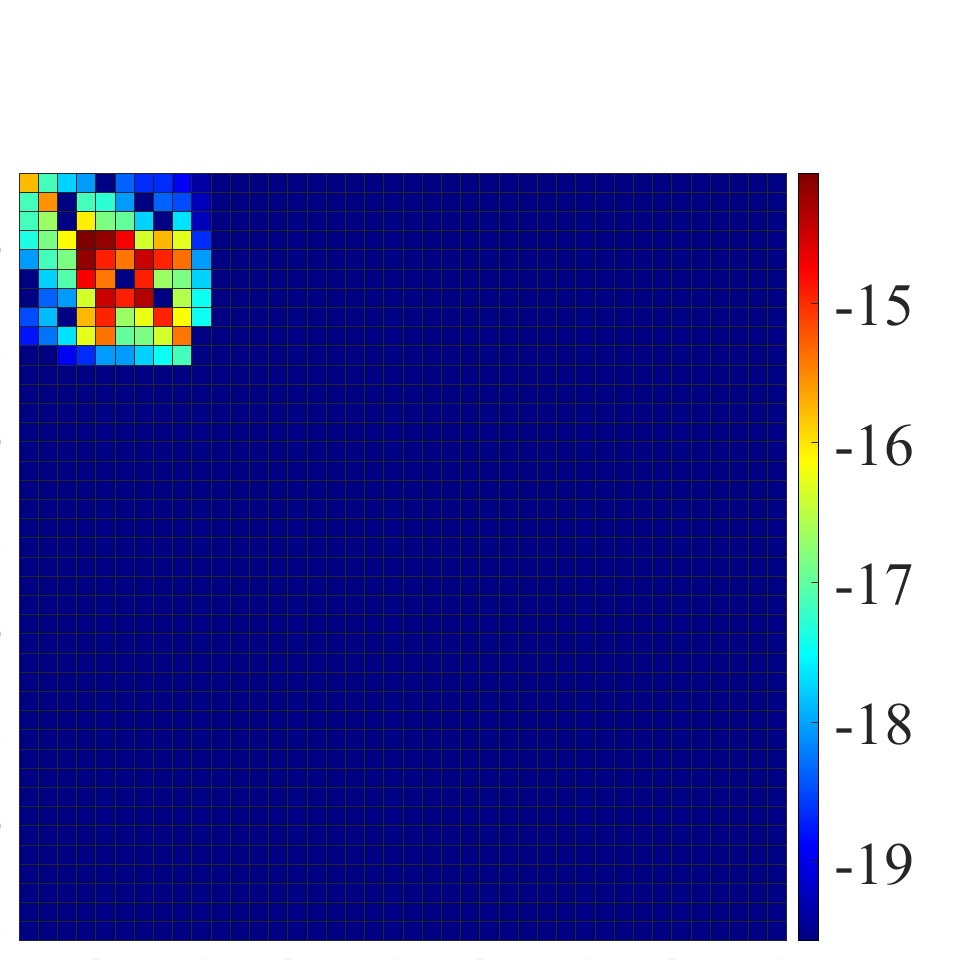}&\includegraphics[width=0.32\linewidth]{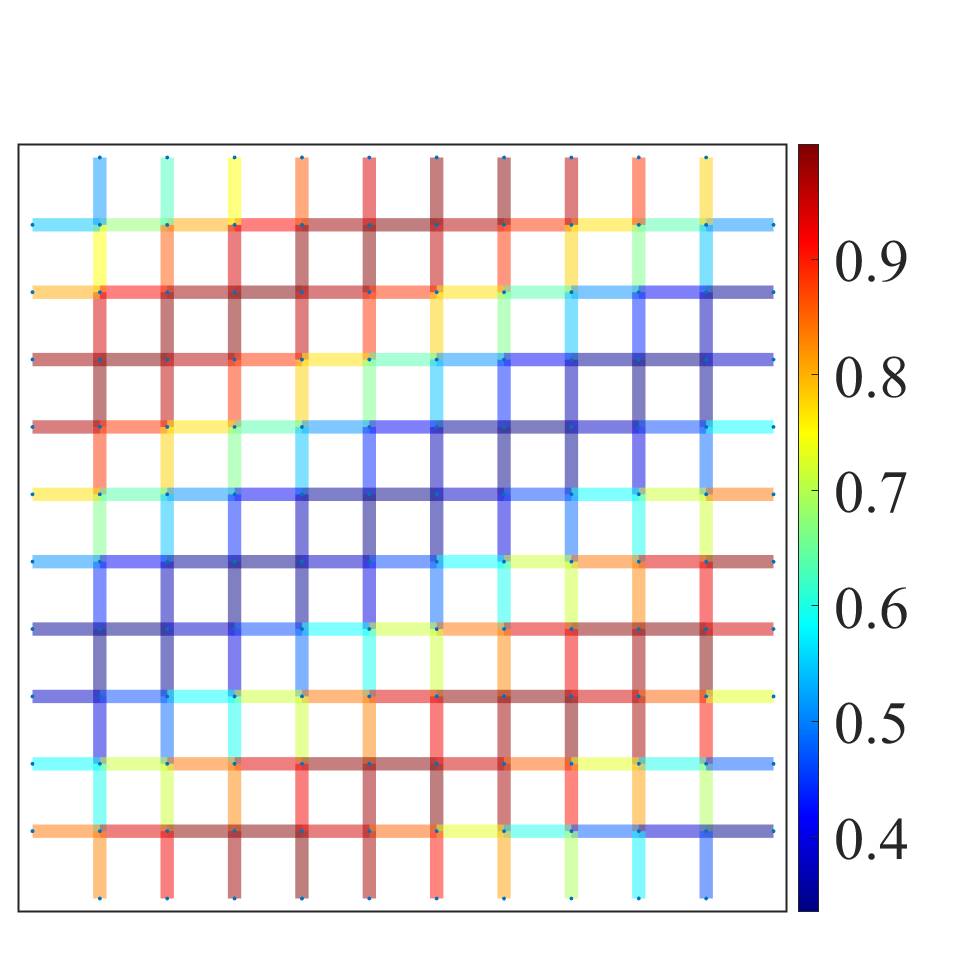}  &\includegraphics[width=0.32\linewidth]{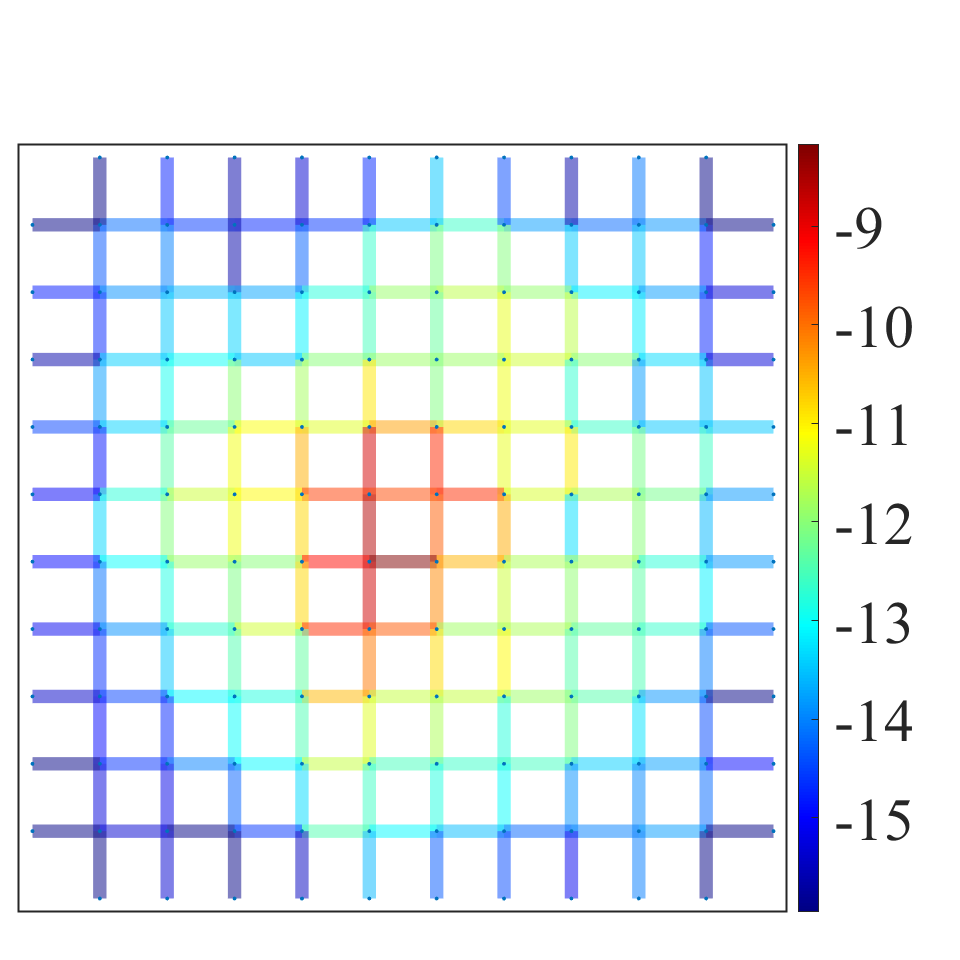}\\
\includegraphics[width=0.33\linewidth]{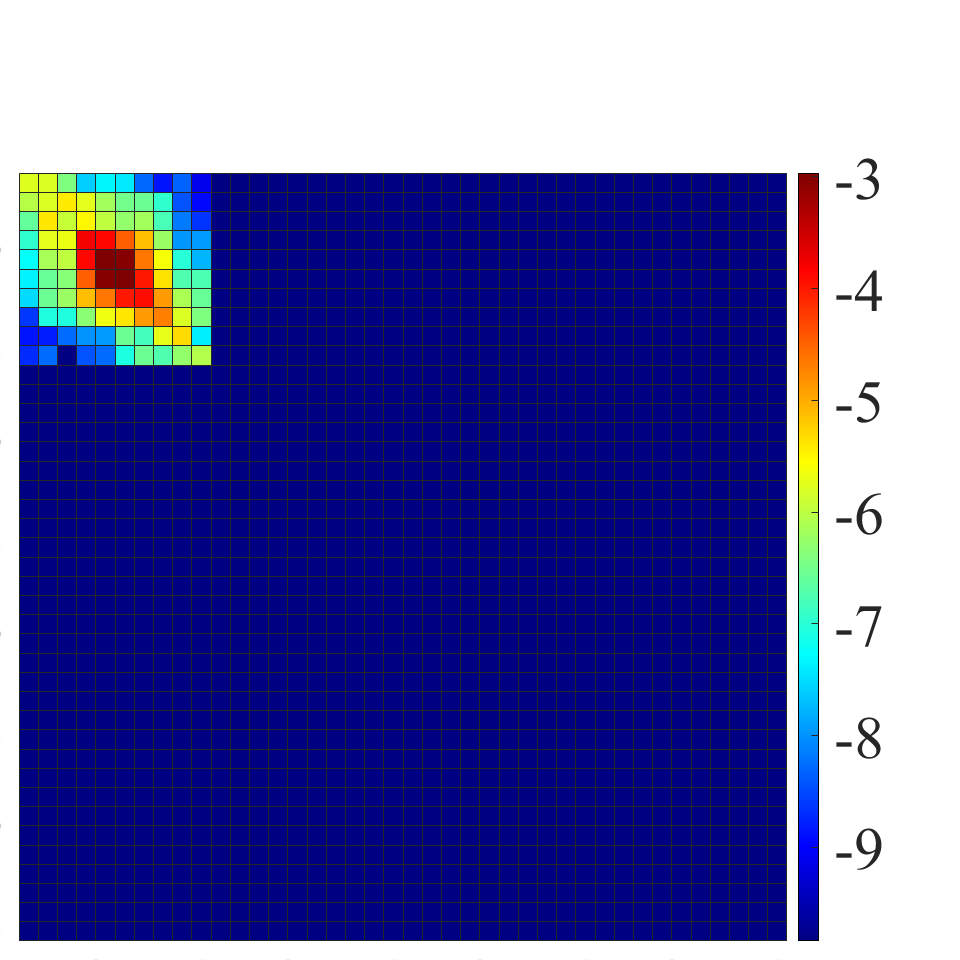}
  &\includegraphics[width=0.32\linewidth]{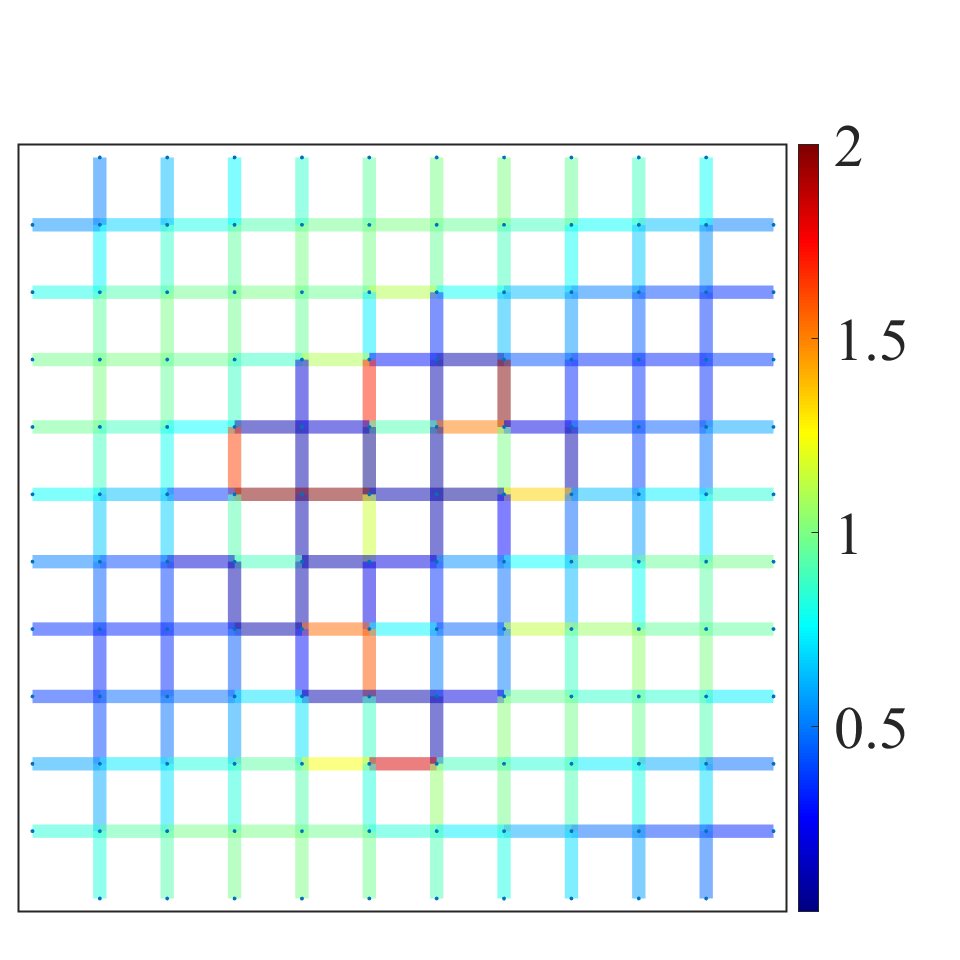}  &\includegraphics[width=0.32\linewidth]{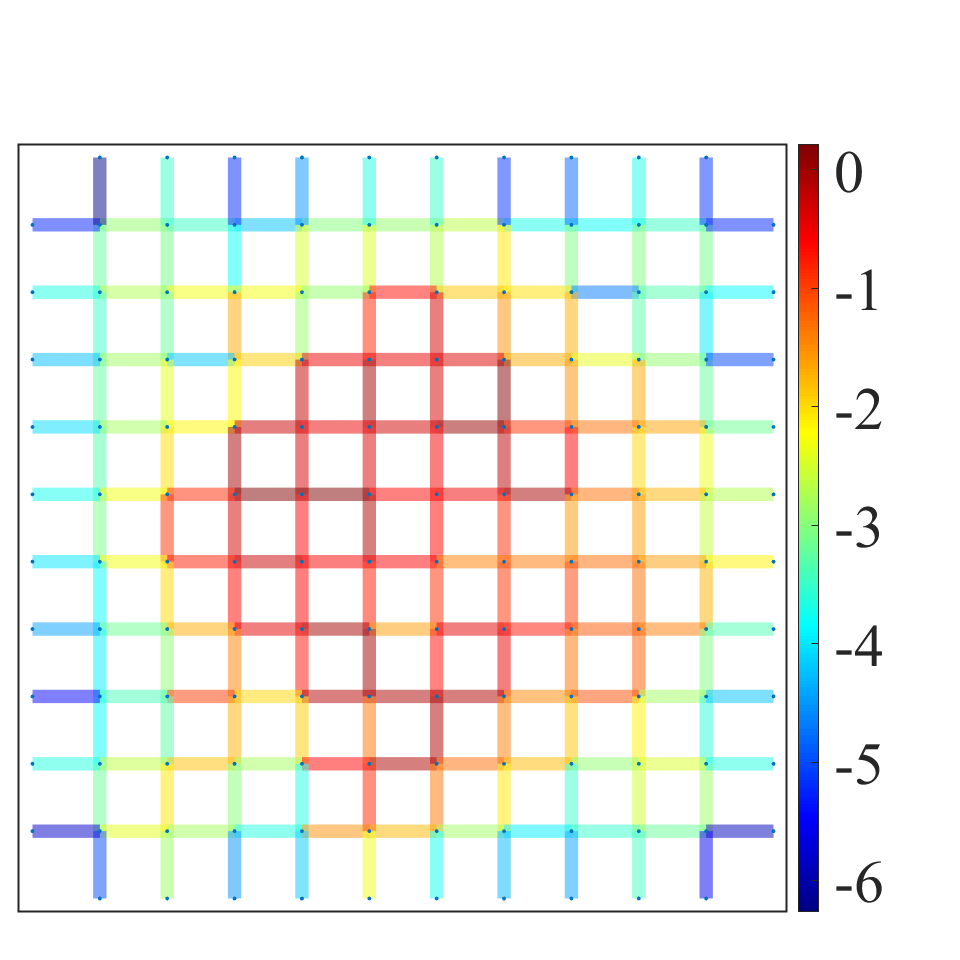}\\
    \includegraphics[width=0.33\linewidth]{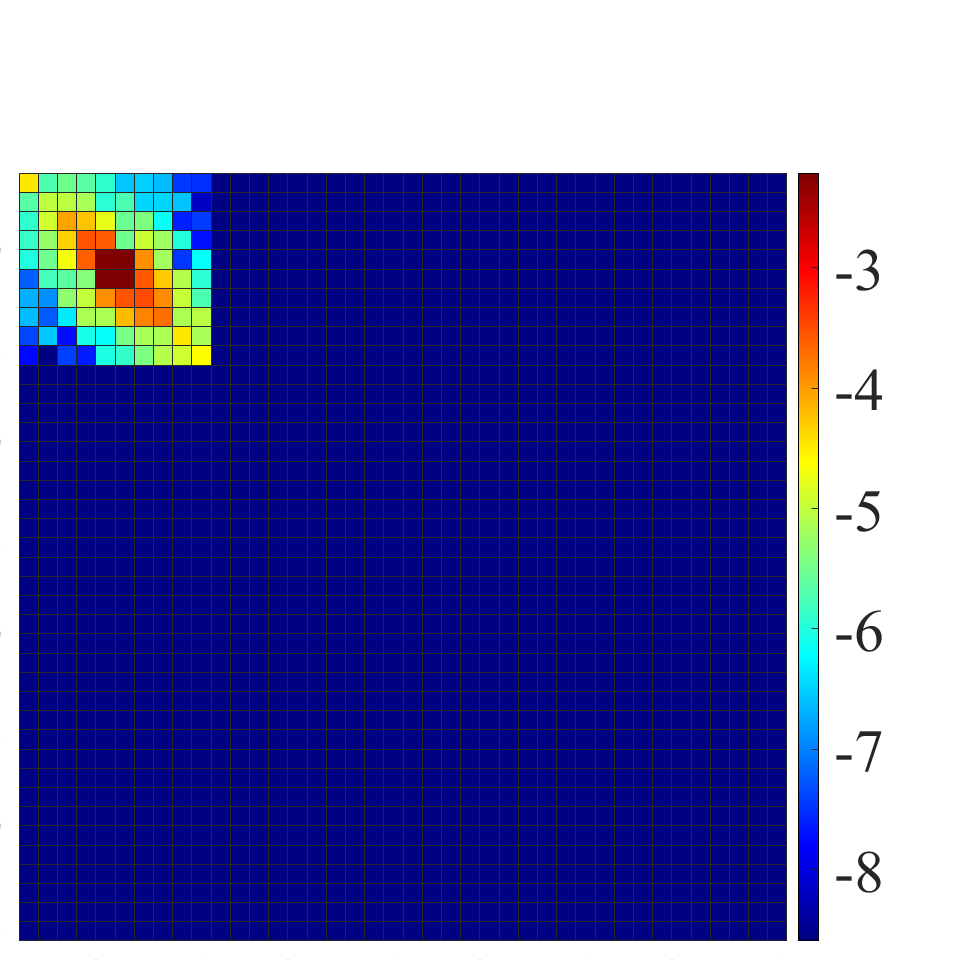}
  &\includegraphics[width=0.32\linewidth]{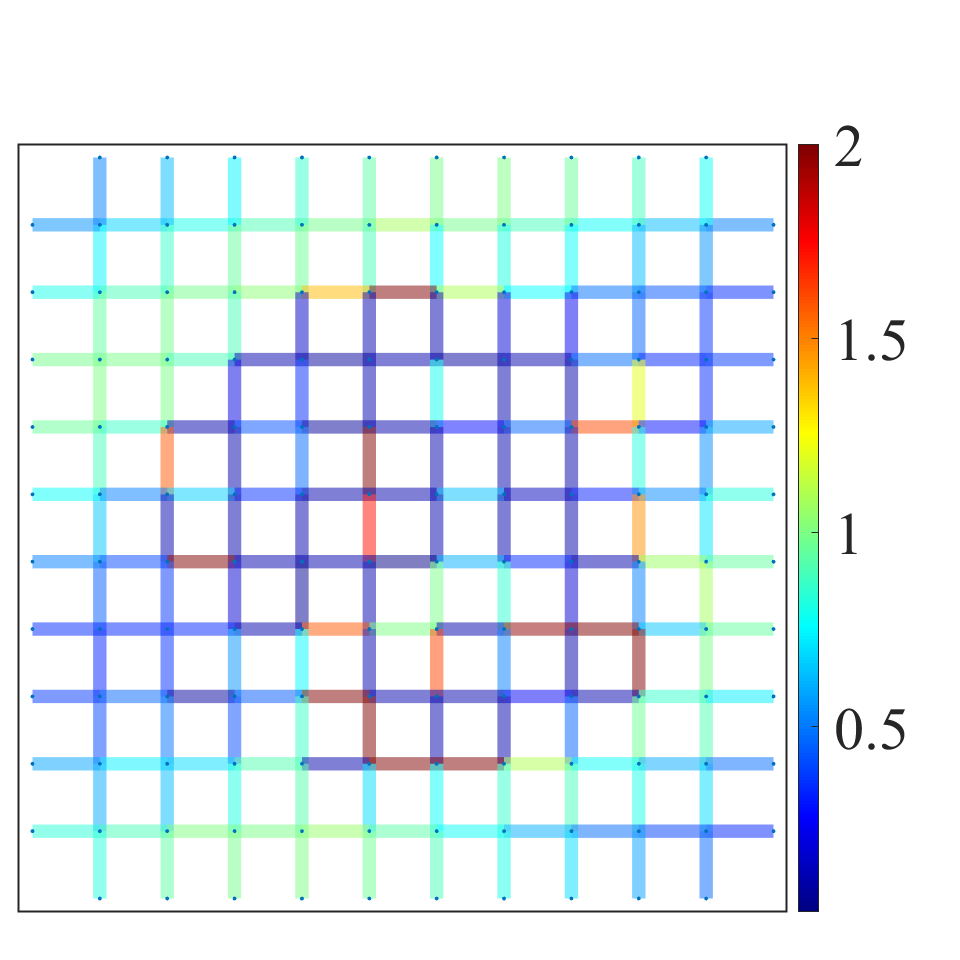}  &\includegraphics[width=0.32\linewidth]{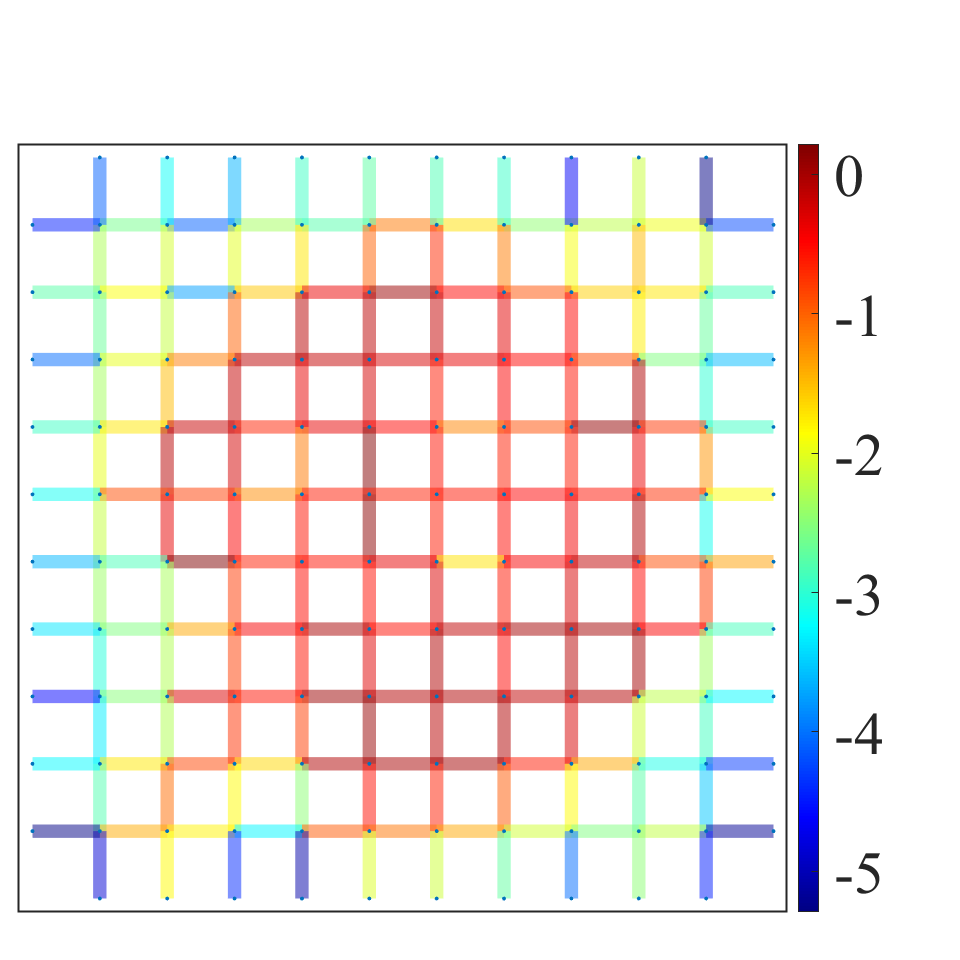}\\
    $\log_{10} |e_{\Lambda}|$  & $\widehat{{\bm\gamma}}$ &$\log_{10} |e_{\bm\gamma}|$
\end{tabular}
\caption{The results with Curtis-Morrow algorithm \cite{curtis1991DNmap} on partial data using the last $3n$ columns of DtN map, at three noise levels $\epsilon=0\%$ (top), $\epsilon=0.001\%$ (middle) and $\epsilon=0.01\%$ (bottom). $ |e_{\Lambda}|$ denotes the component  error between $\Lambda_\varepsilon$ and the recovery $\Lambda_\varepsilon^r$ by the Curtis-Morrow algorithm.} \label{fig:partialCM}
\end{figure}

\begin{figure}[hbt!]
    \centering
    \setlength{\tabcolsep}{0pt}
    \begin{tabular}{ccc|cc}
 \toprule
\multicolumn{3}{c}{ Curtis-Morrow algorithm} & \multicolumn{2}{c}{NN approach}\\
 \midrule
    \includegraphics[width=0.21\linewidth]{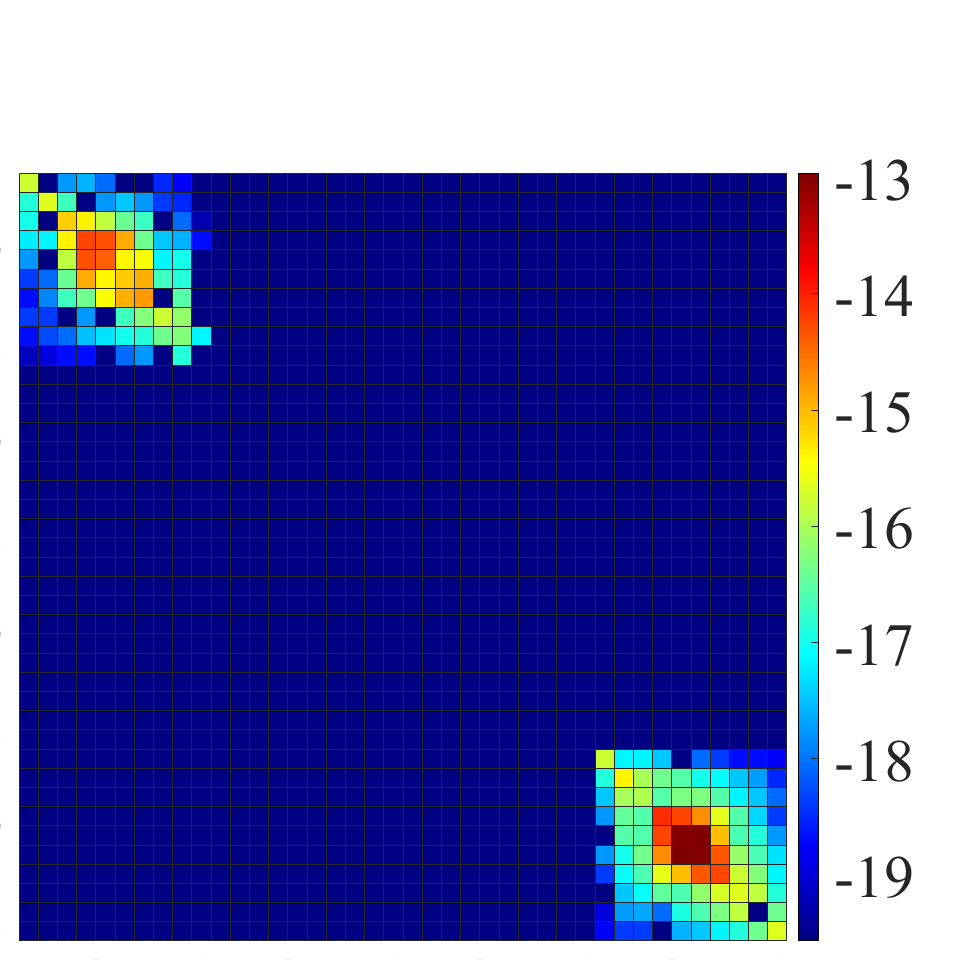} &\includegraphics[width=0.2\linewidth]{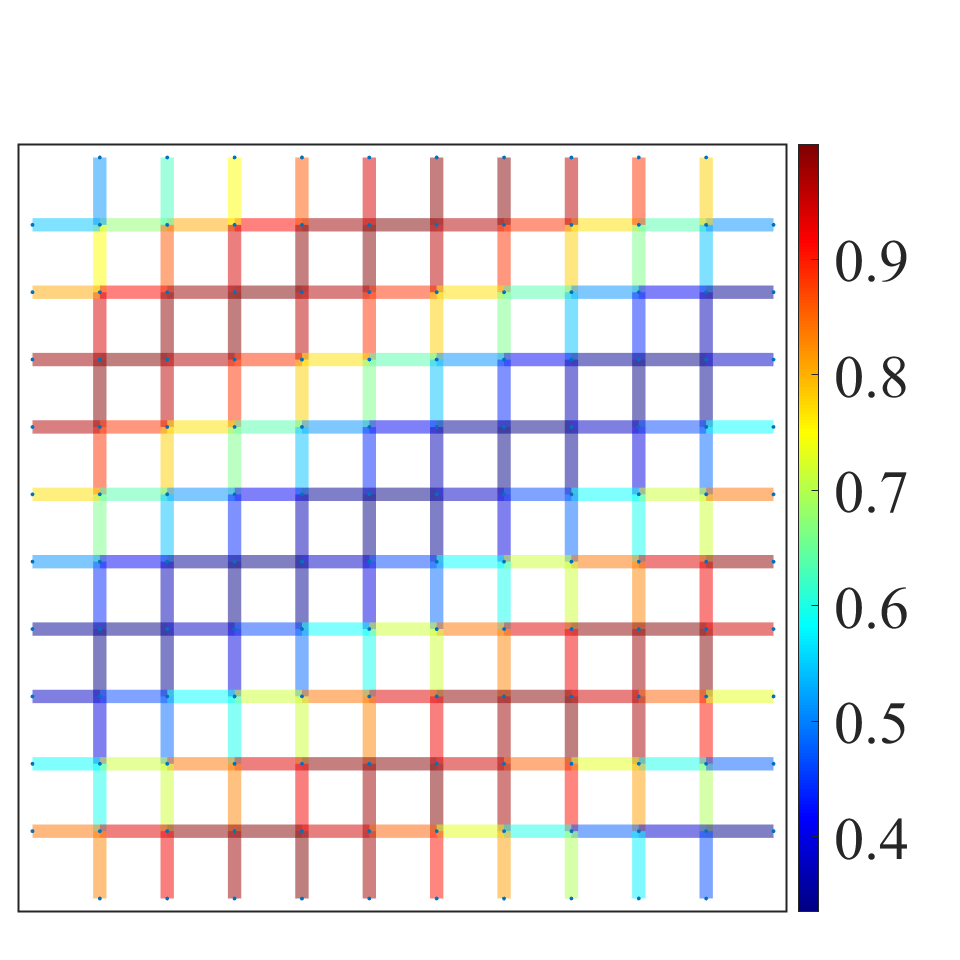}  &\includegraphics[width=0.2\linewidth]{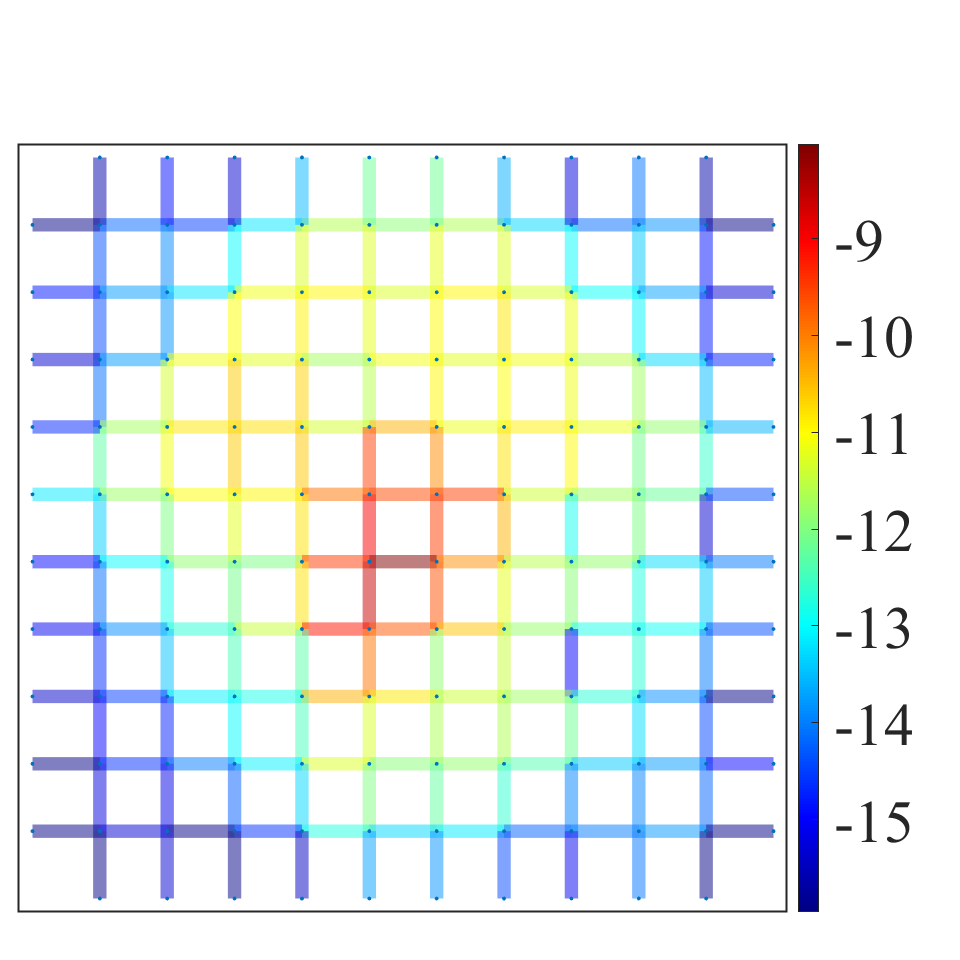}&\includegraphics[width=0.2\linewidth]{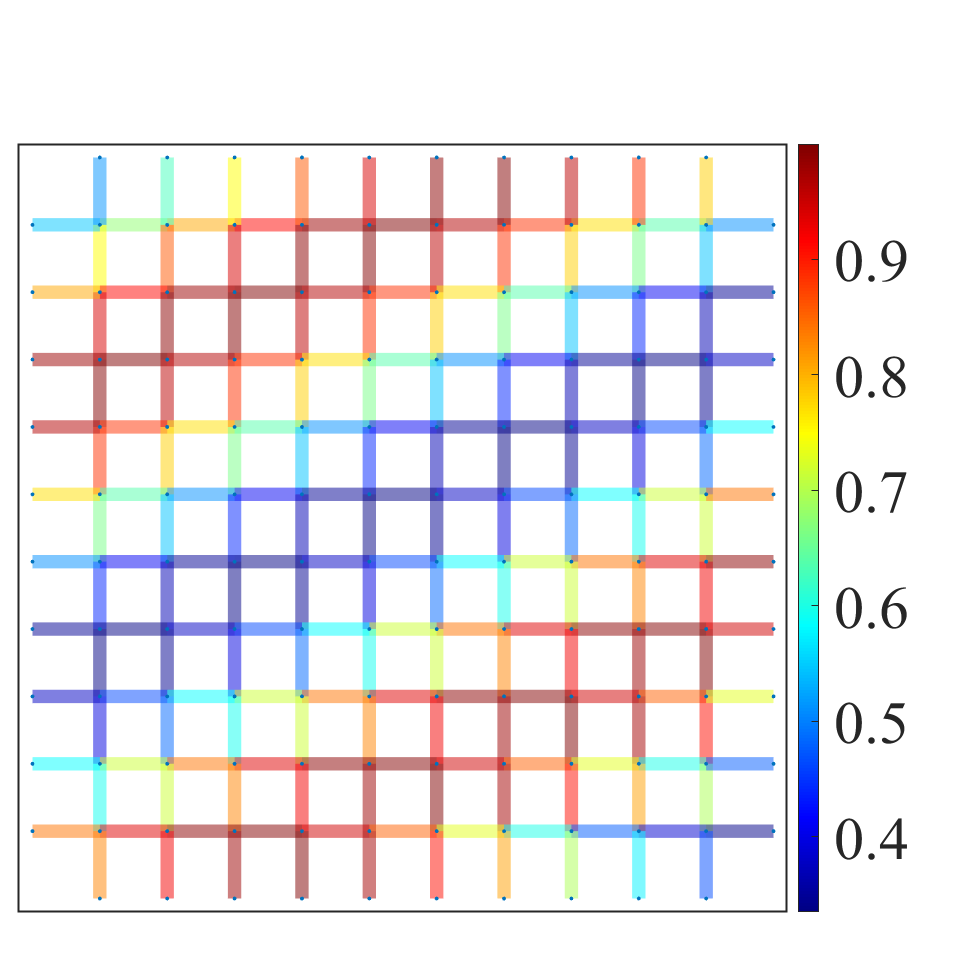} &\includegraphics[width=0.2\linewidth]{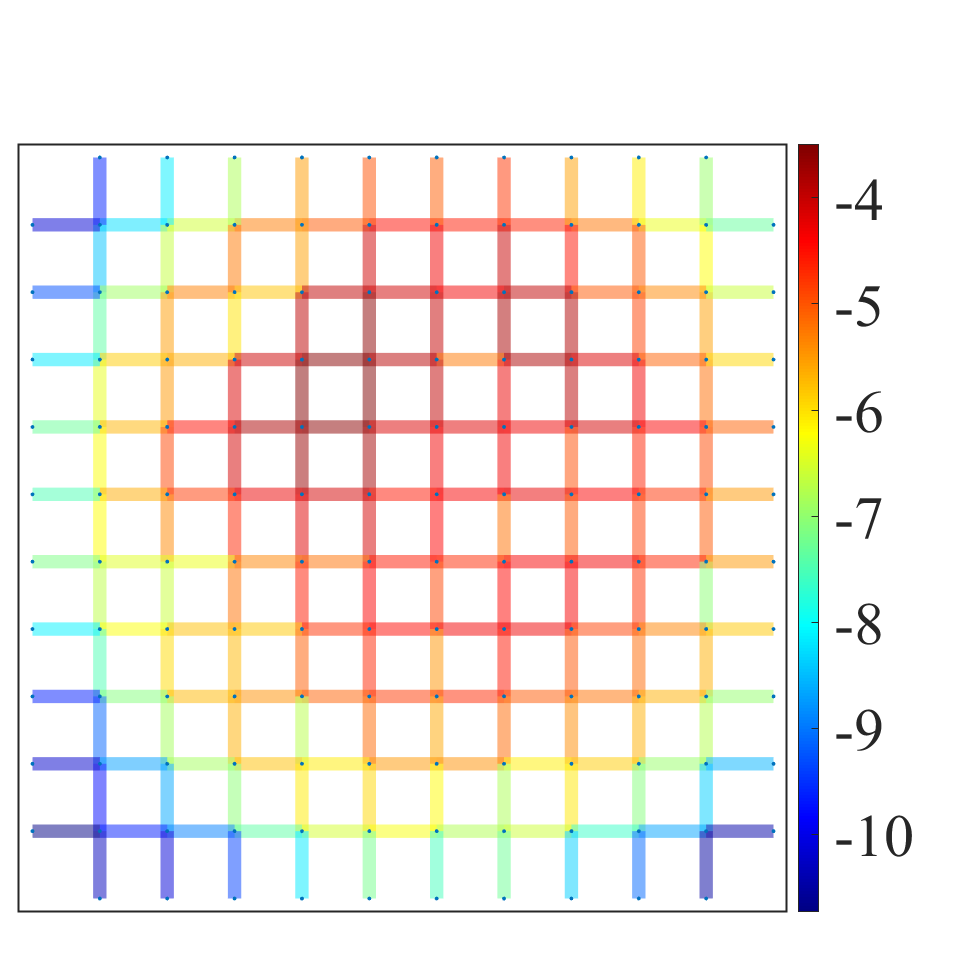}\\ 
    \includegraphics[width=0.21\linewidth]{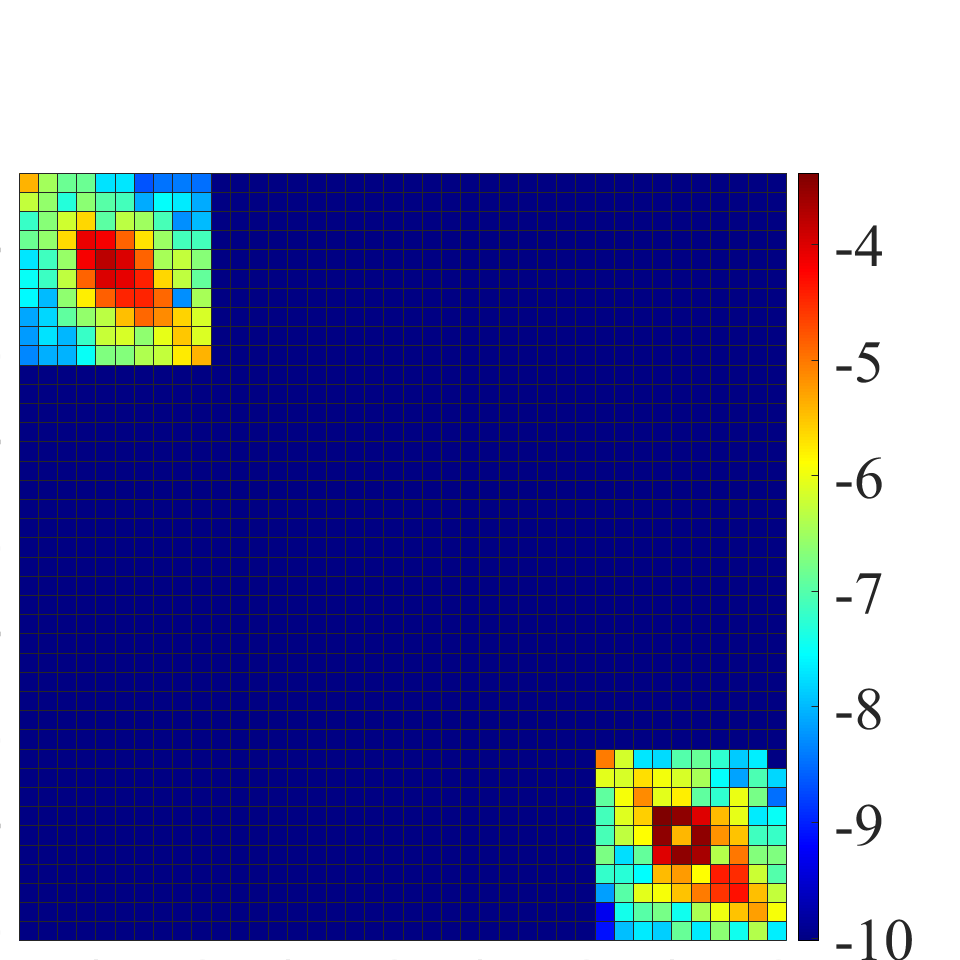}&\includegraphics[width=0.2\linewidth]{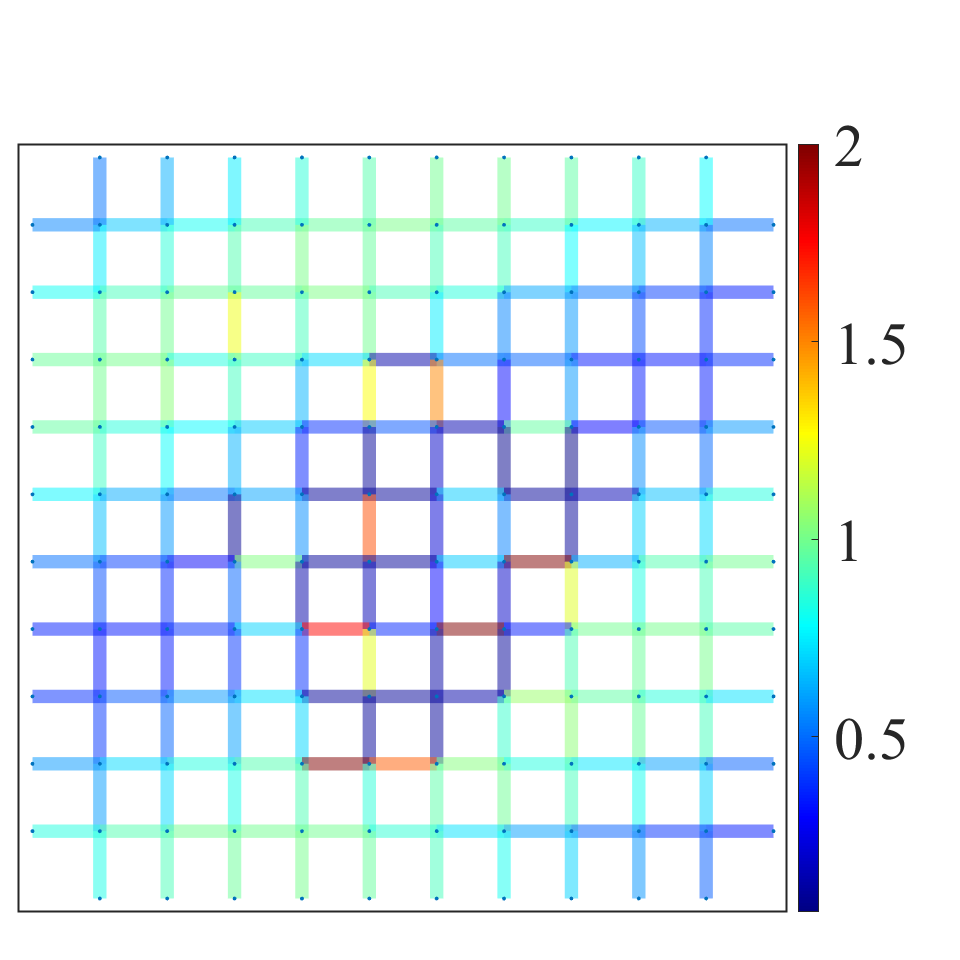}  &\includegraphics[width=0.2\linewidth]{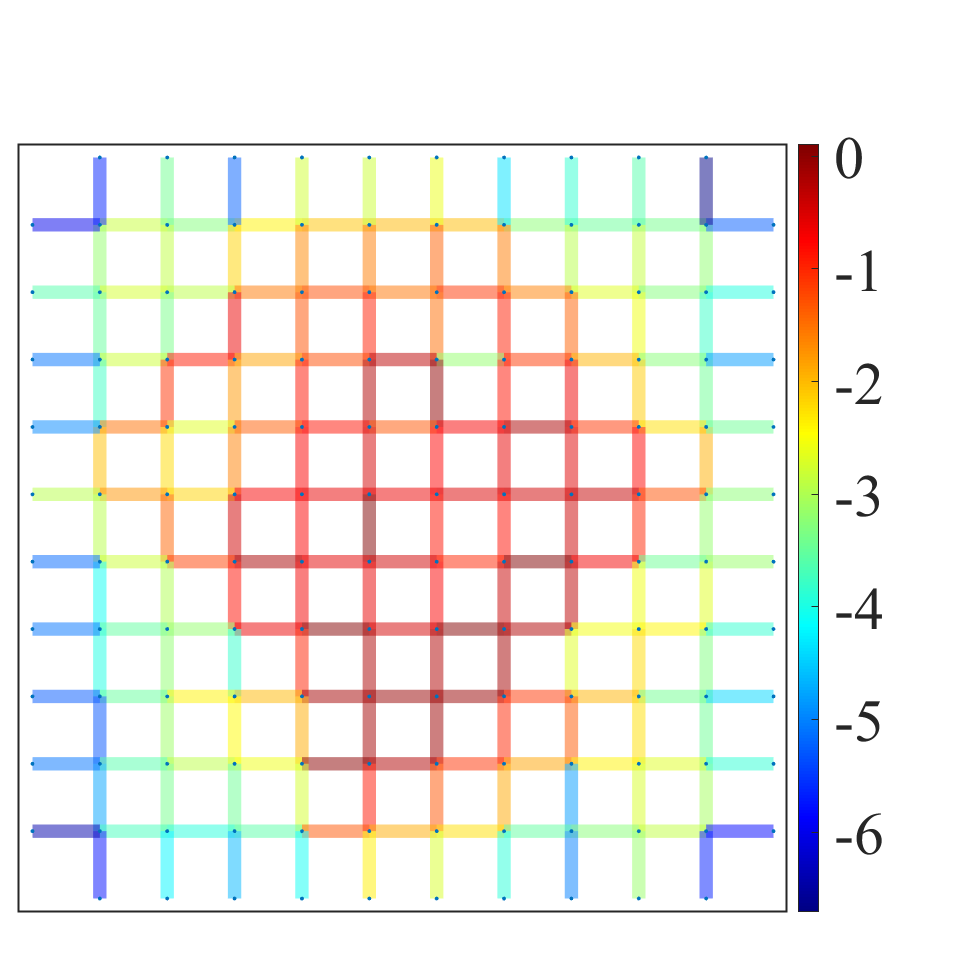}&\includegraphics[width=0.2\linewidth]{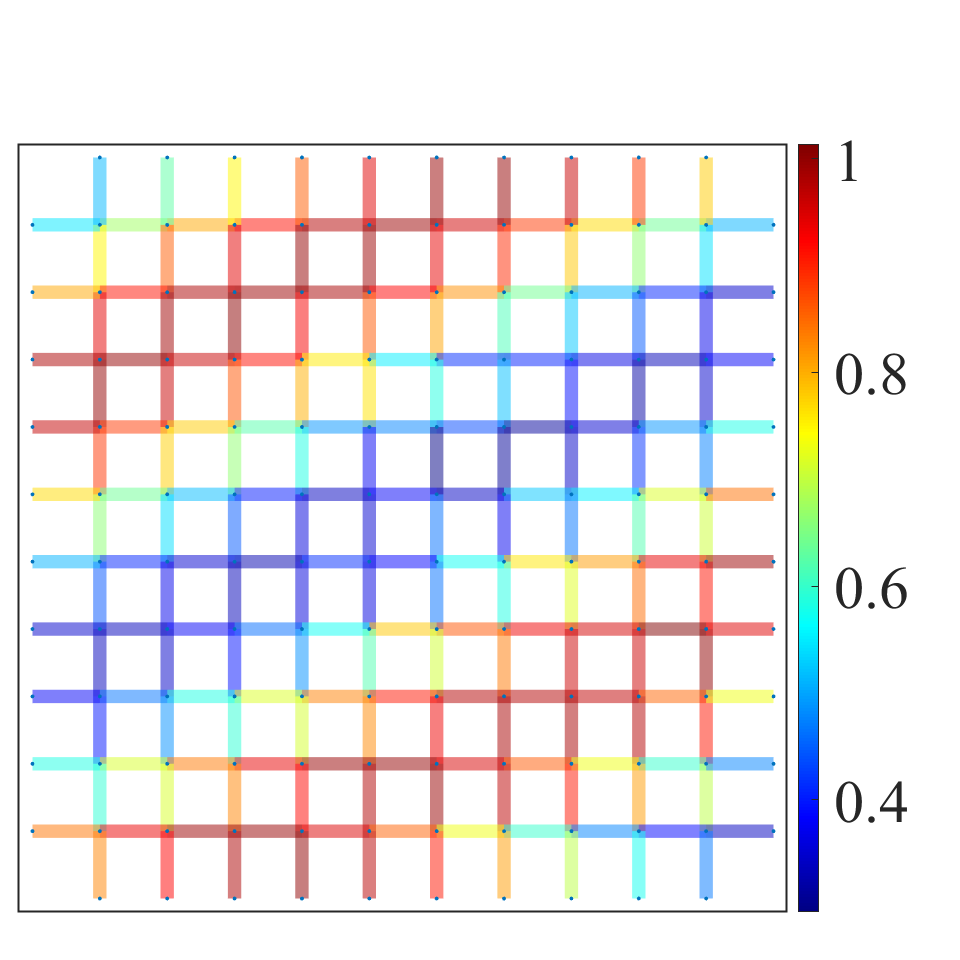} &\includegraphics[width=0.2\linewidth]{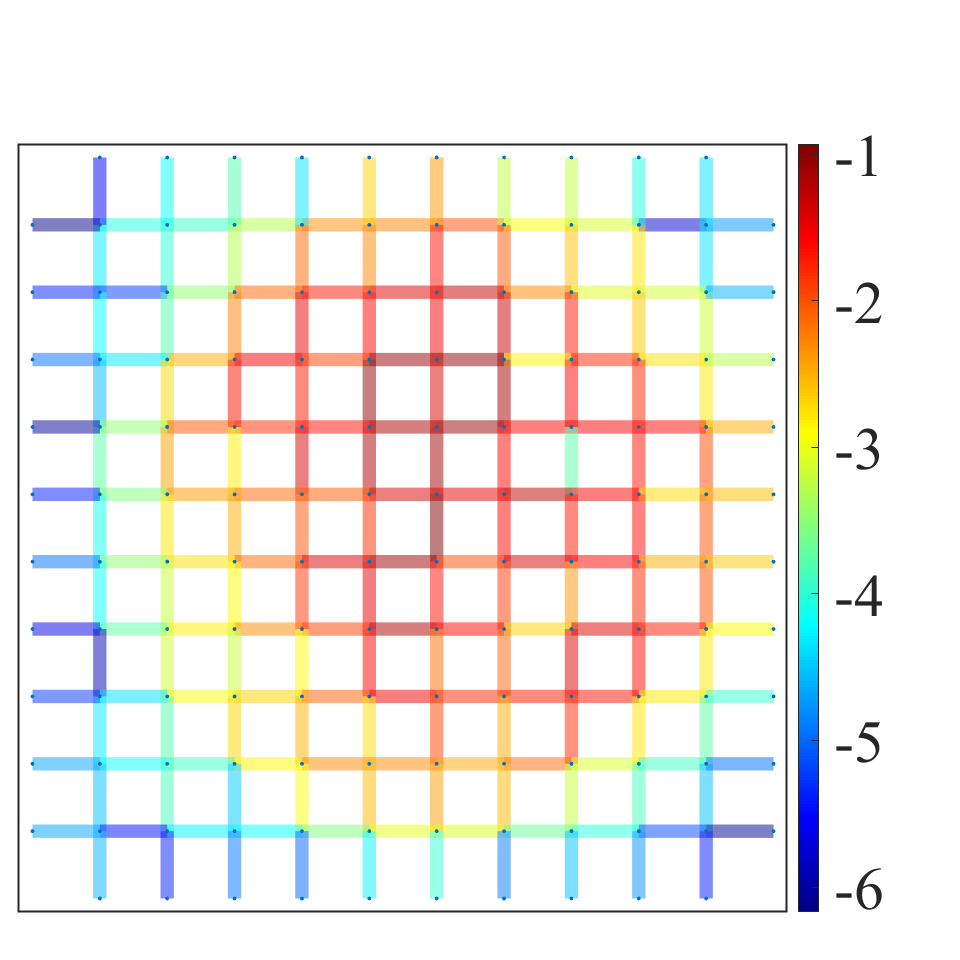}\\ 
    \includegraphics[width=0.21\linewidth]{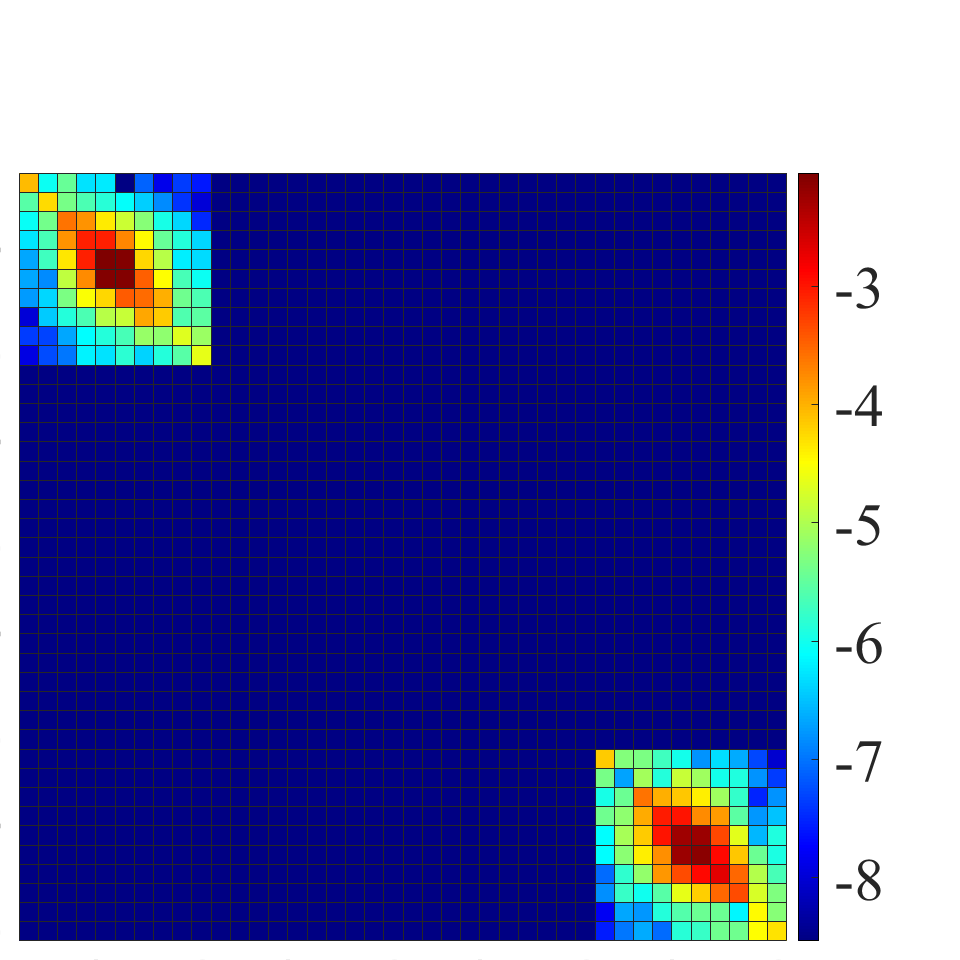} &\includegraphics[width=0.2\linewidth]{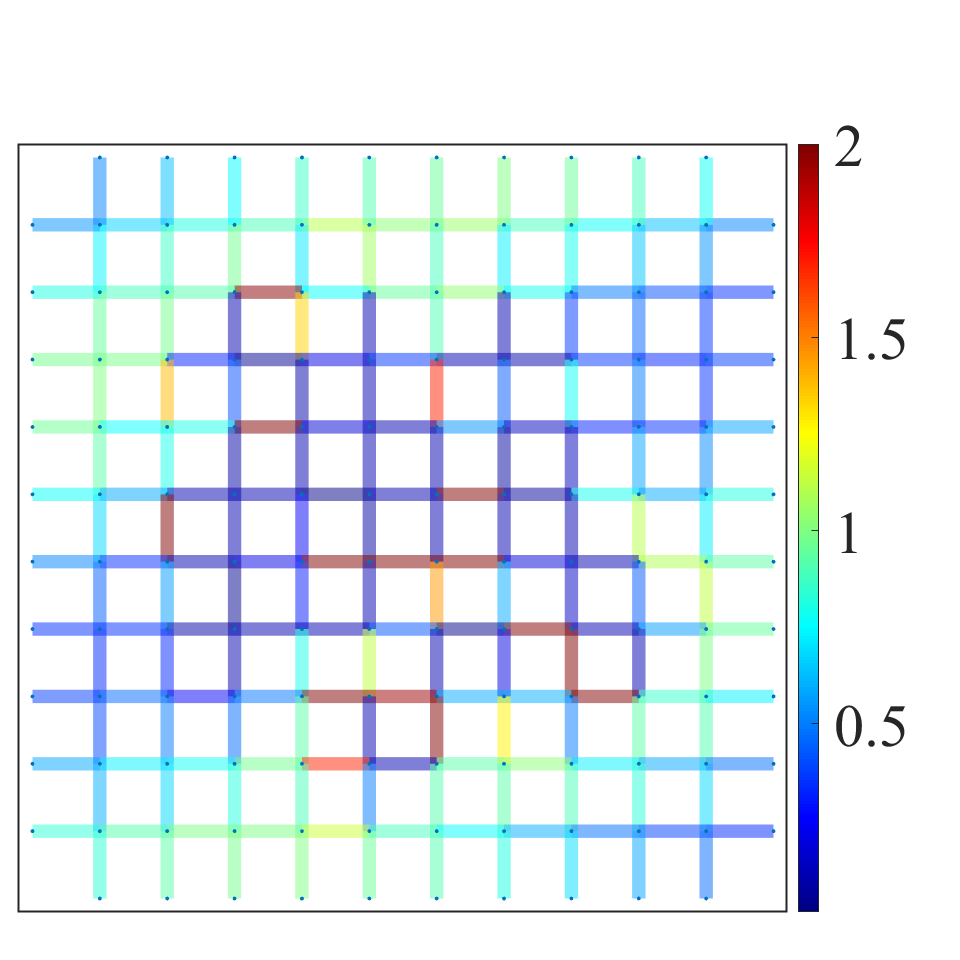}  &\includegraphics[width=0.2\linewidth]{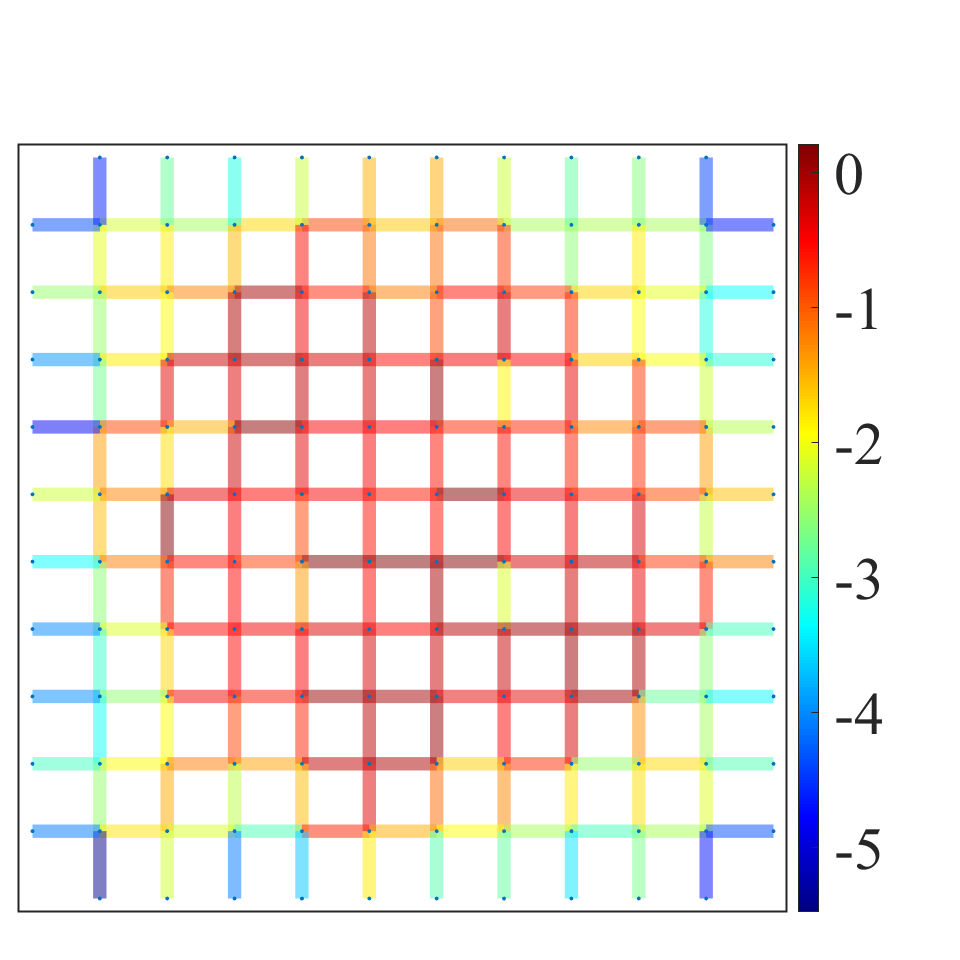}&\includegraphics[width=0.2\linewidth]{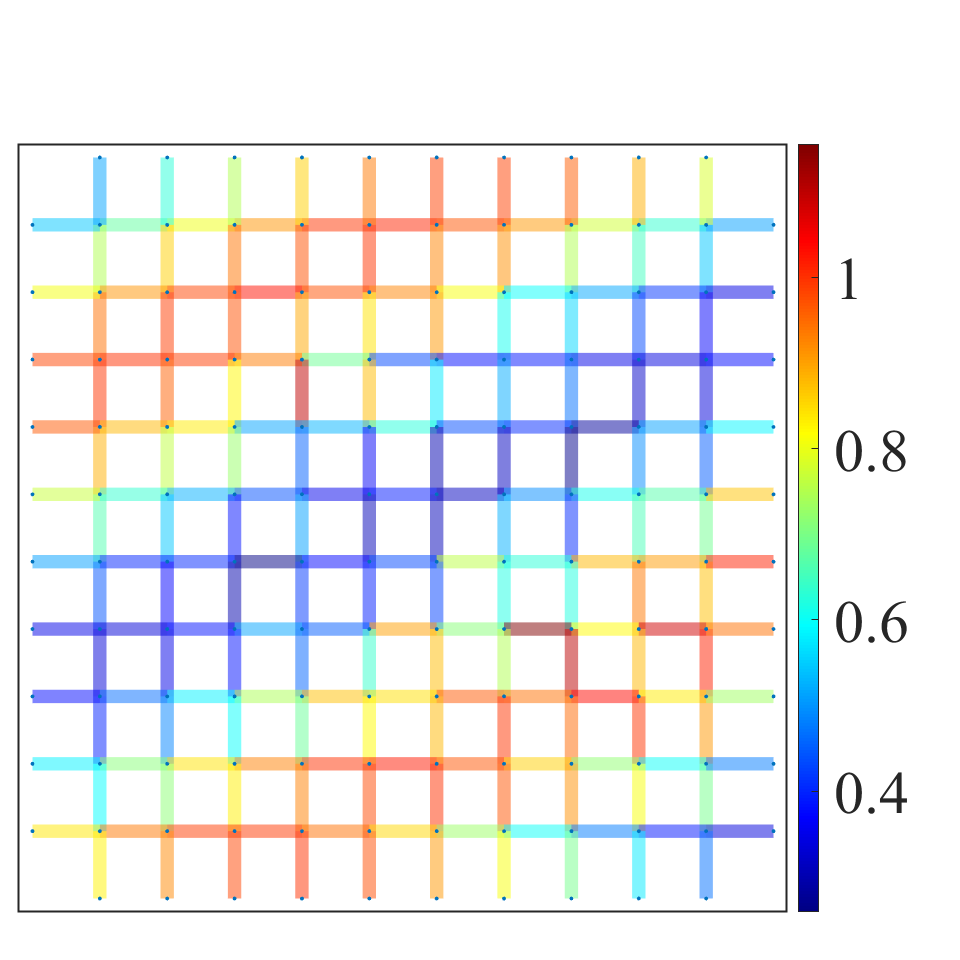} &\includegraphics[width=0.2\linewidth]{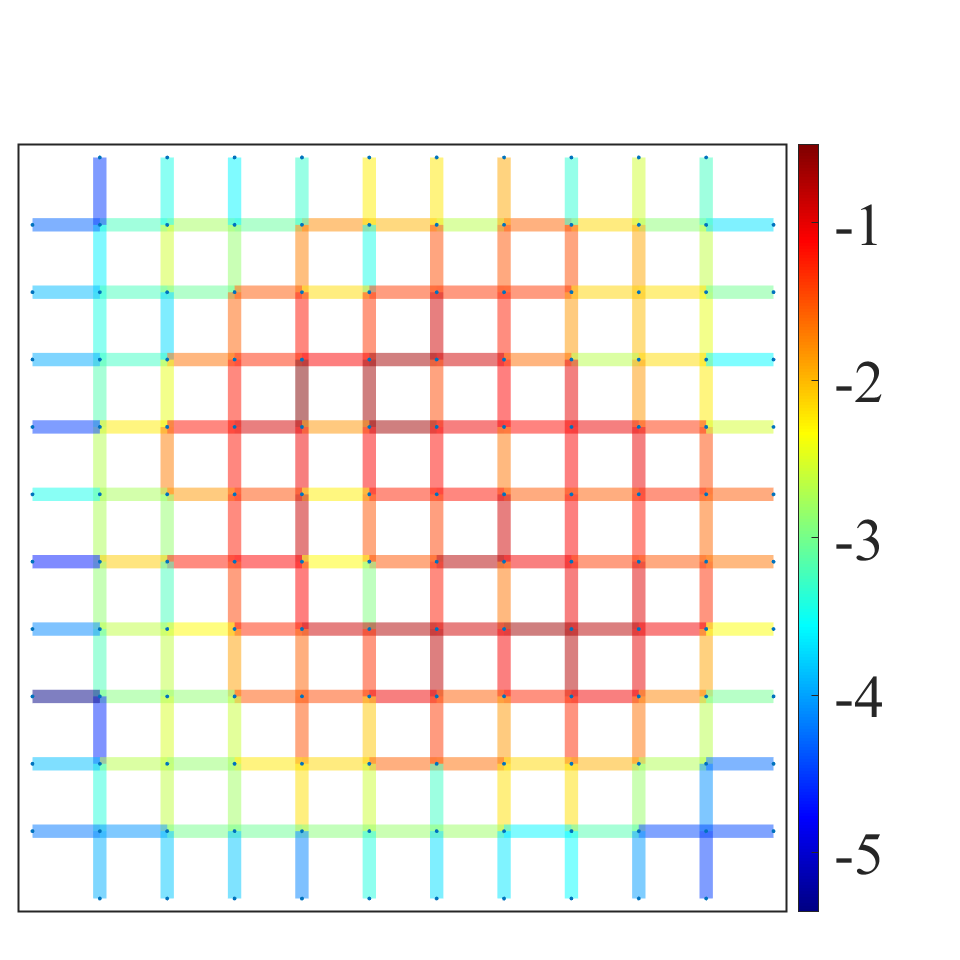}\\
   $\log_{10} |e_{\Lambda}|$ & $\widehat {\bm\gamma}$ & $\log_{10}|e_{\bm\gamma}|$ &  $\widehat {\bm\gamma}$ & $\log_{10}|e_{\bm\gamma}|$
    \end{tabular}
    \caption{
    Numerical results by the Curtis-Morrow algorithm and NN approach for partial data over the set $(1:3 n)\times(n+1:4 n)$ at three noise levels, $\epsilon=0\%$ (top), $\epsilon=0.001\%$ (middle) and $\epsilon=0.01\%$ (bottom). }
    \label{fig:Partial-square}
\end{figure}

The numerical results for the cases of partial data (i.e., $T_{\rm sq}$ and $T_{\rm re}$) are shown in Figs. \ref{fig:Partial-square} and \ref{fig:partial-Rectangle}. The overall observations are similar to the case of incomplete DtN map: for exact data, the Curtis-Morrow algorithm works for both datasets and is more accurate than the NN approach, but for noisy data, it completely fails to produce acceptable reconstruction (for $\epsilon=0.001\%$). The case $T_{\rm re}$
is even more challenging: it poses challenge to the NN approach, which manages to recover the overall shape, but it causes the Curtis-Morrow algorithm to fail completely. These results again show the superior performance of the NN approach for incomplete / partial DtN data.

\begin{figure}[hbt!]
    \centering
    \setlength{\tabcolsep}{0pt}
    \begin{tabular}{ccc|cc}
 \toprule
\multicolumn{3}{c}{Curtis-Morrow algorithm}&\multicolumn{2}{c}{NN approach} \\
 \midrule
\includegraphics[width=0.21\linewidth]{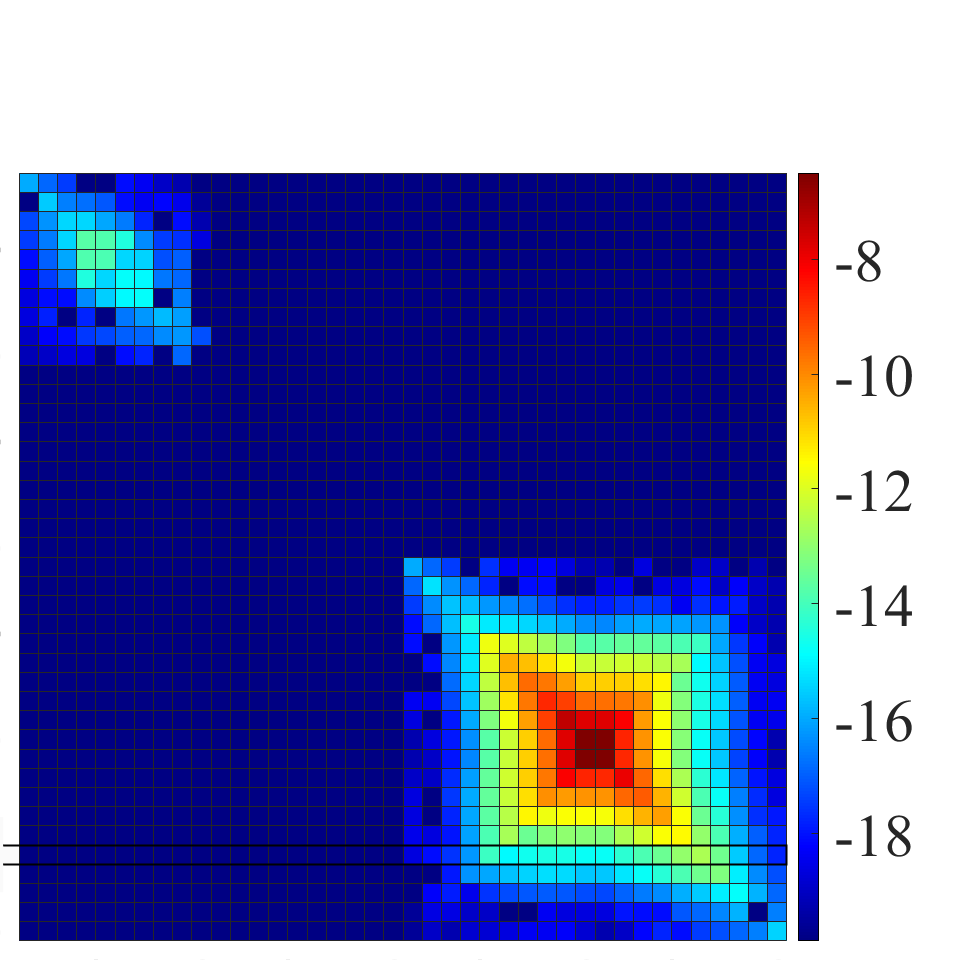}&\includegraphics[width=0.2\linewidth]{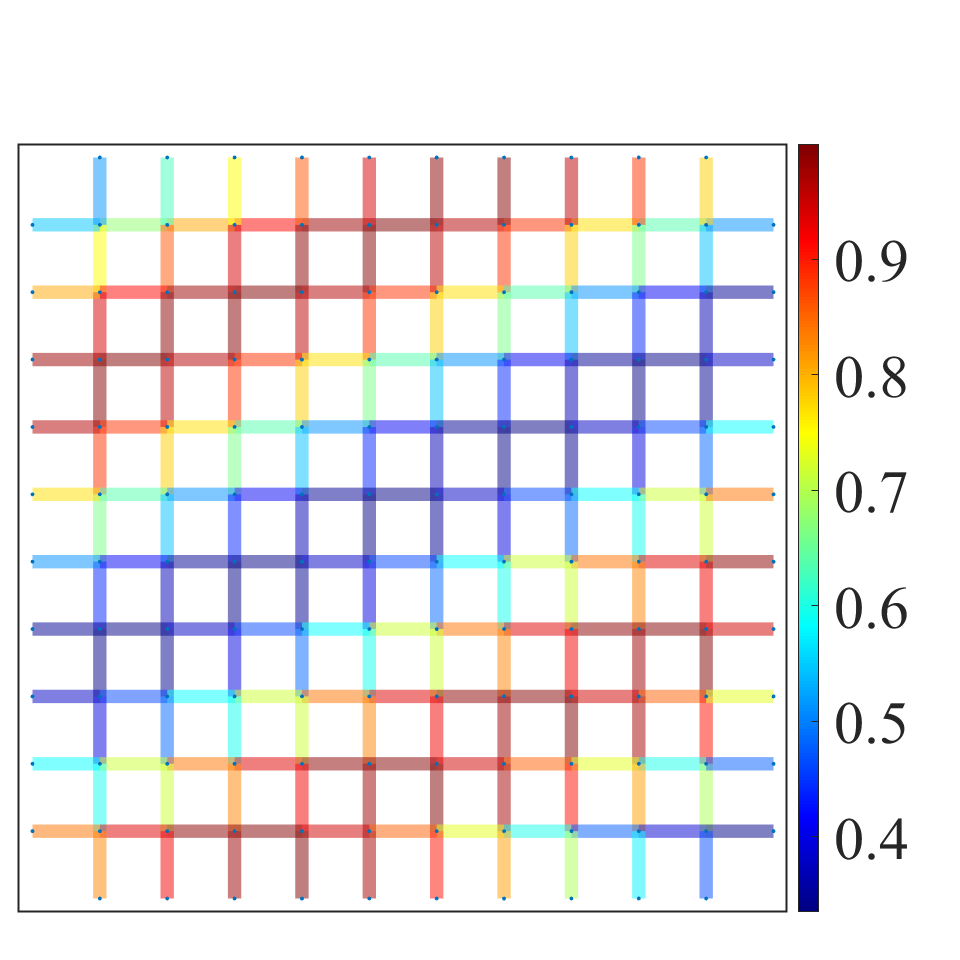}  &\includegraphics[width=0.2\linewidth]{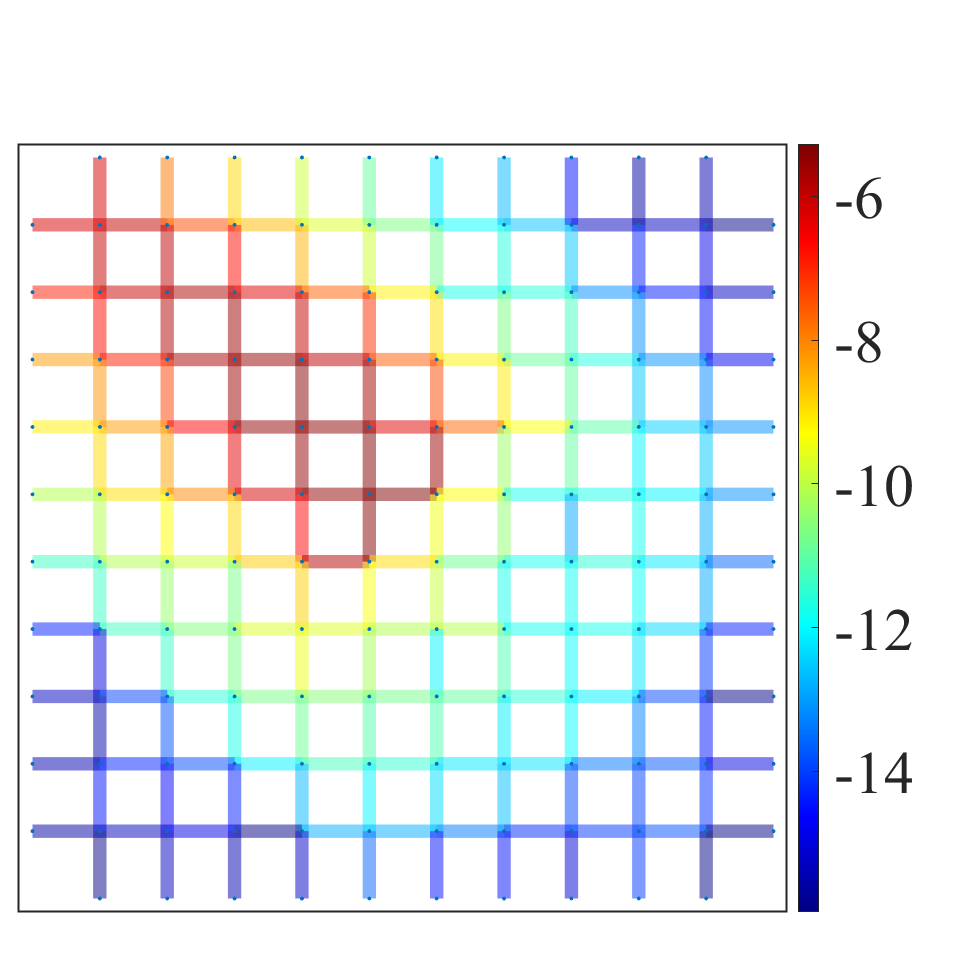}&\includegraphics[width=0.2\linewidth]{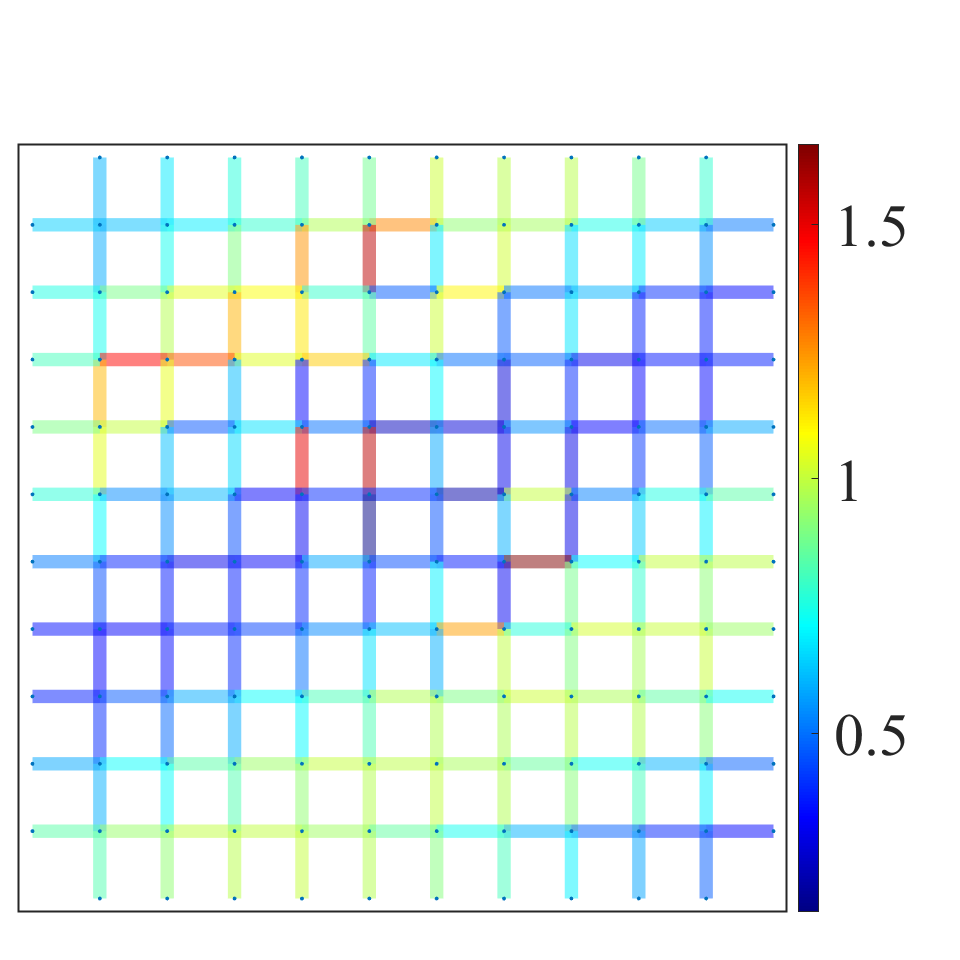} &\includegraphics[width=0.2\linewidth]{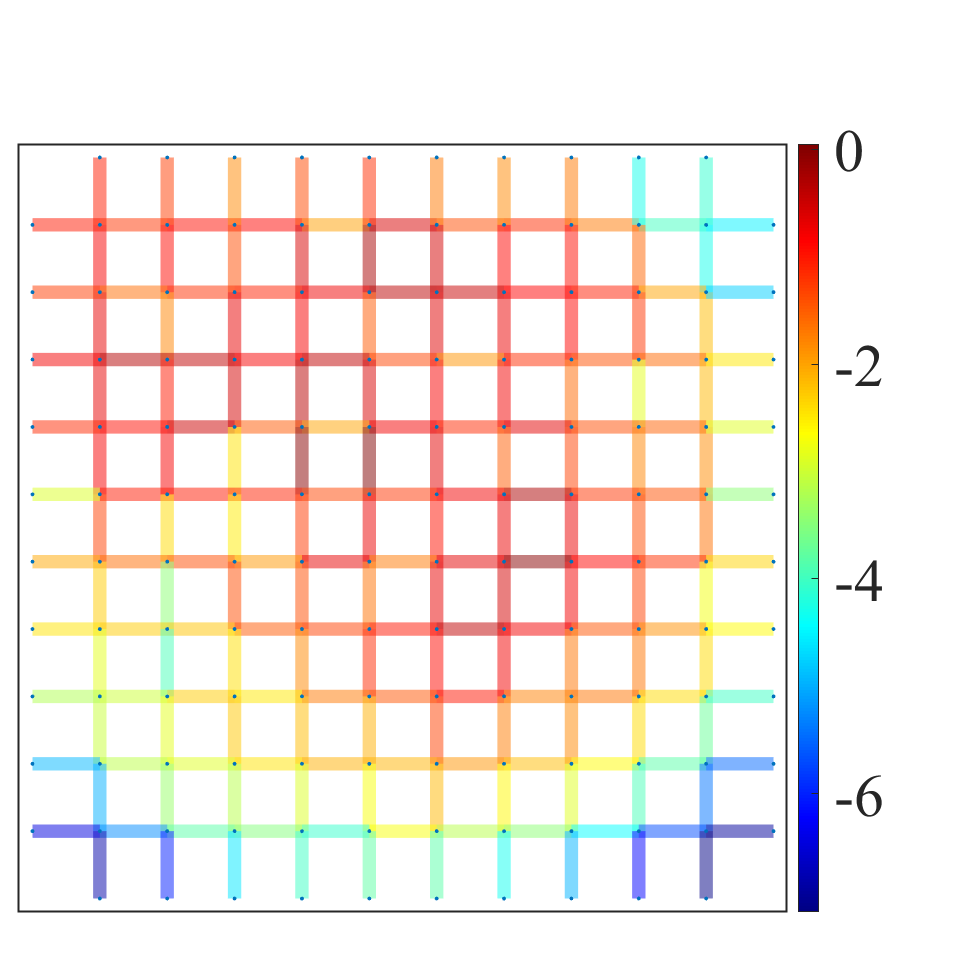}\\ 
\includegraphics[width=0.21\linewidth]{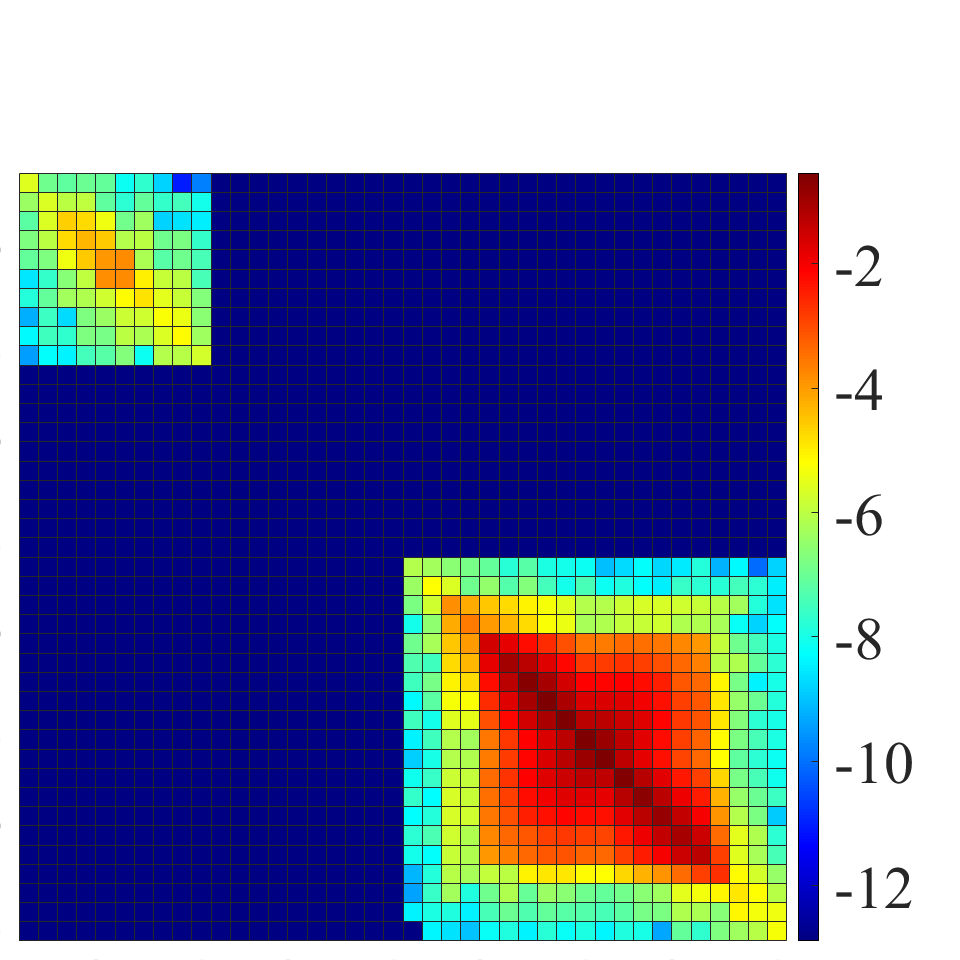} &\includegraphics[width=0.2\linewidth]{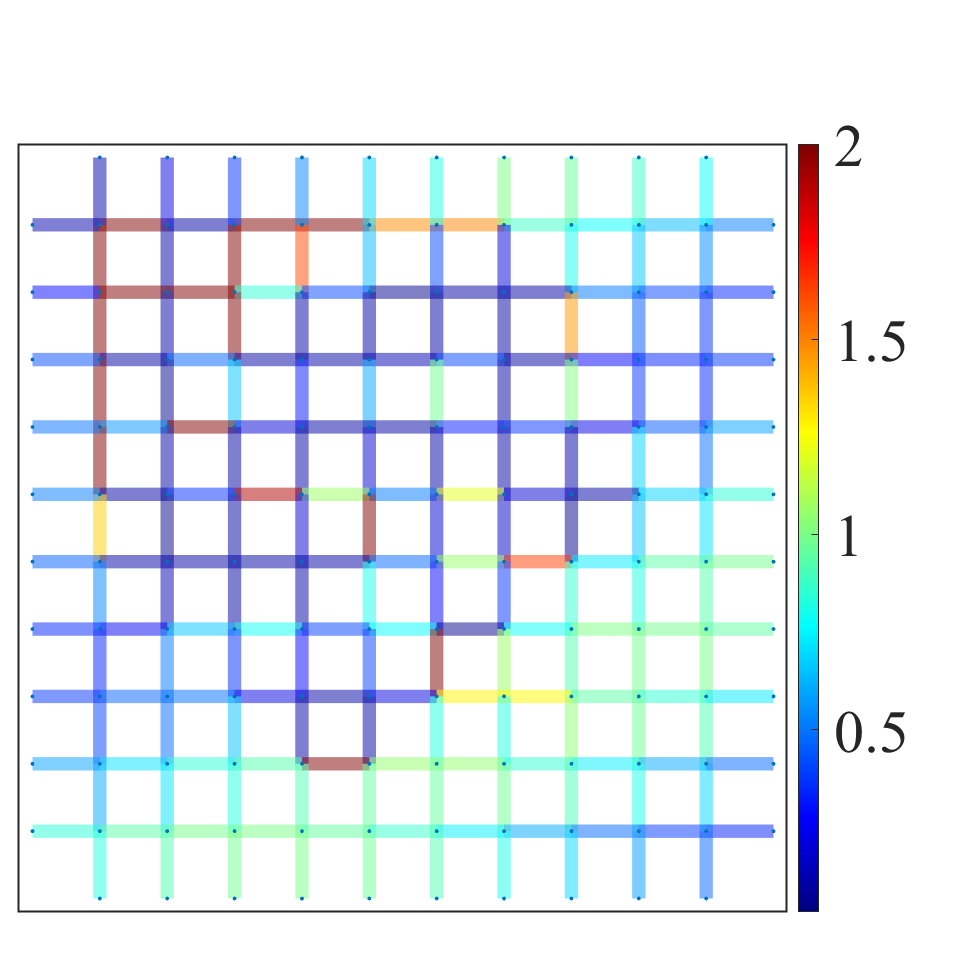} &\includegraphics[width=0.2\linewidth]{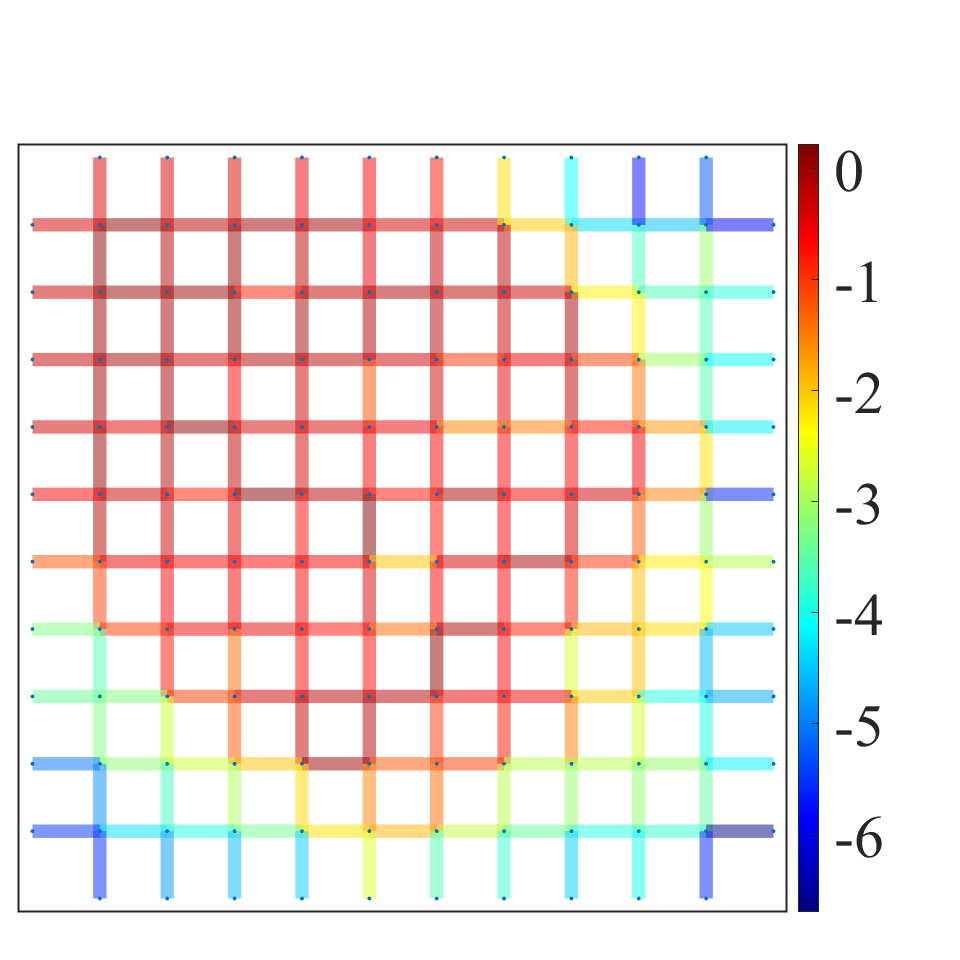}&\includegraphics[width=0.2\linewidth]{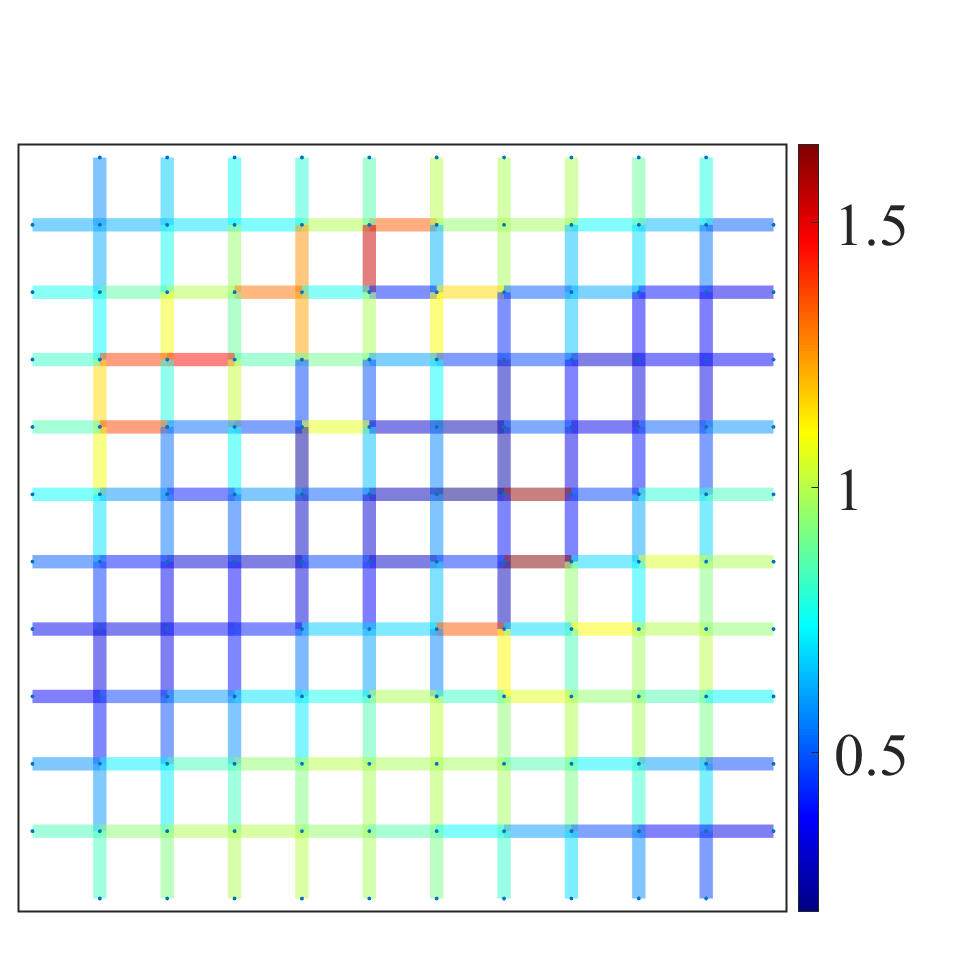} &\includegraphics[width=0.2\linewidth]{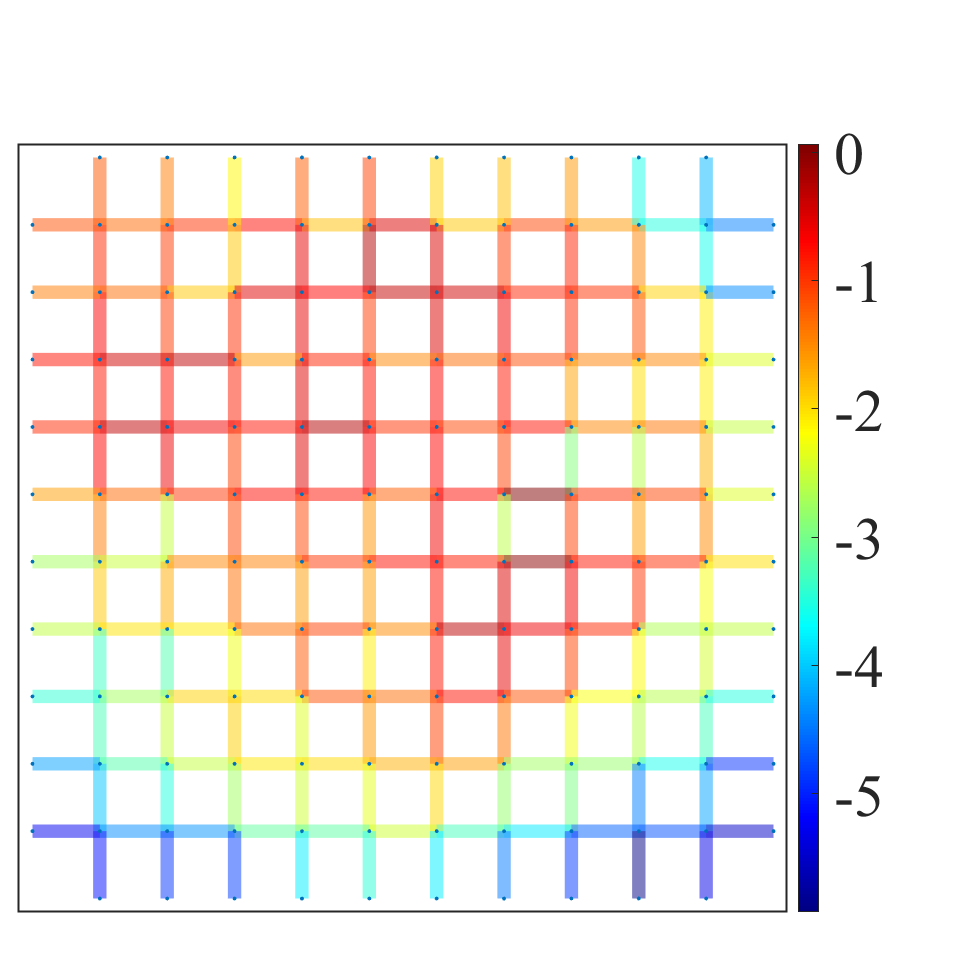}\\
  \includegraphics[width=0.21\linewidth]{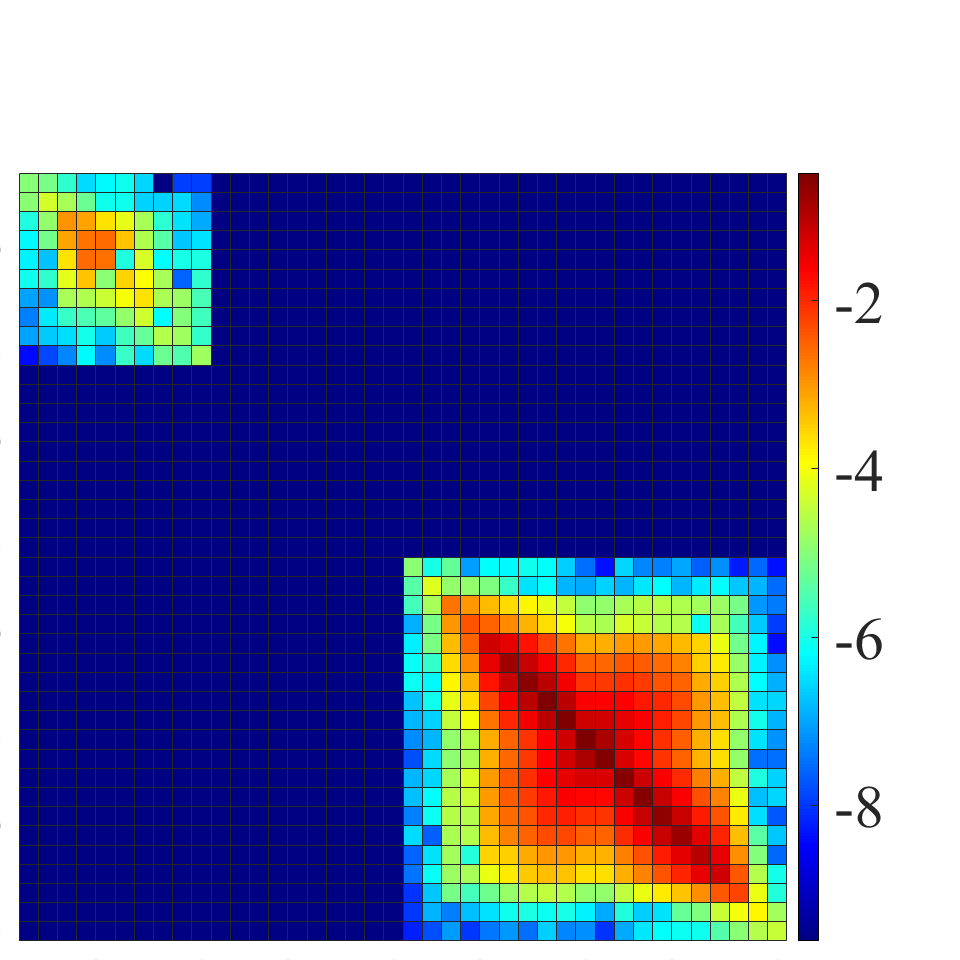} &\includegraphics[width=0.2\linewidth]{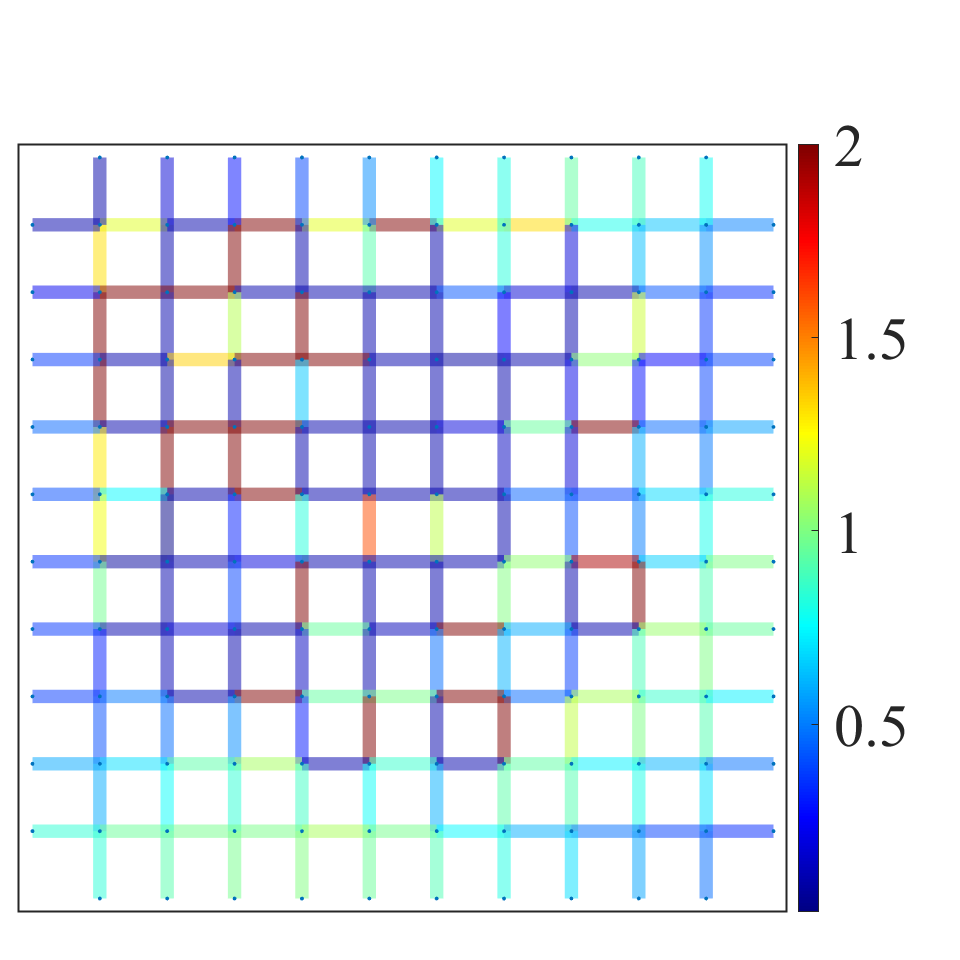}  &\includegraphics[width=0.2\linewidth]{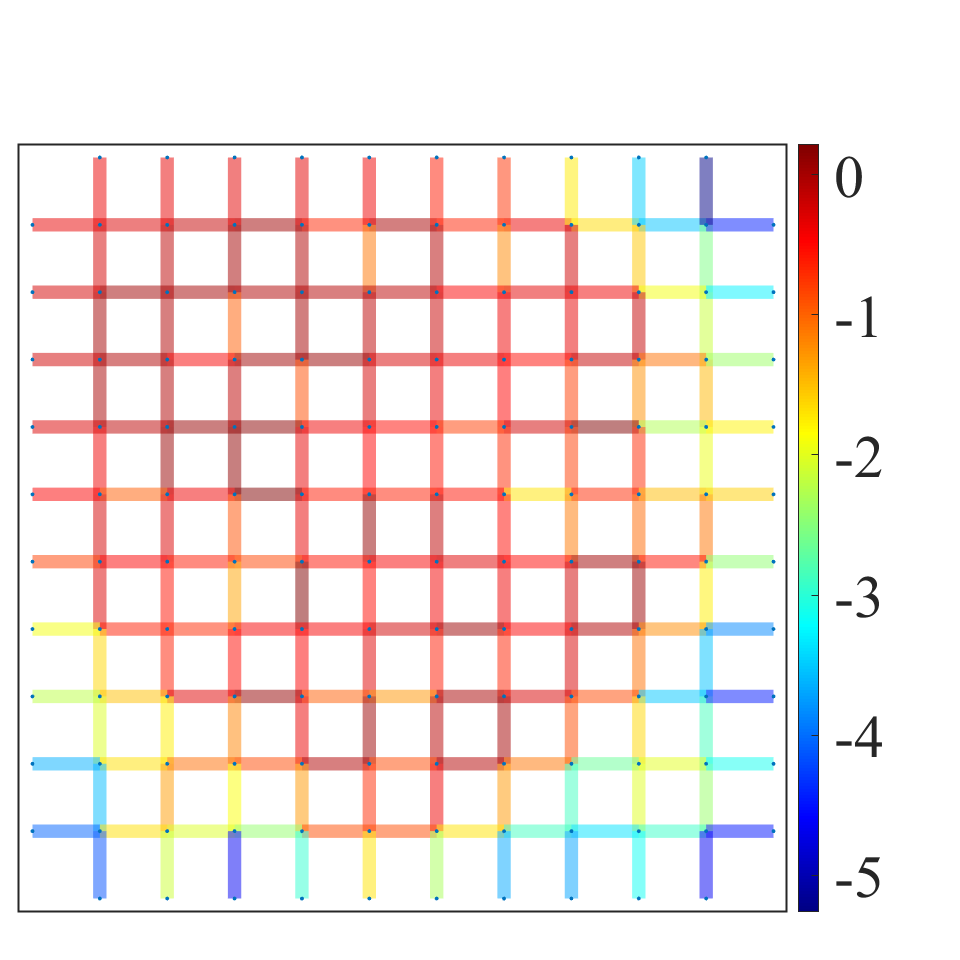}&\includegraphics[width=0.2\linewidth]{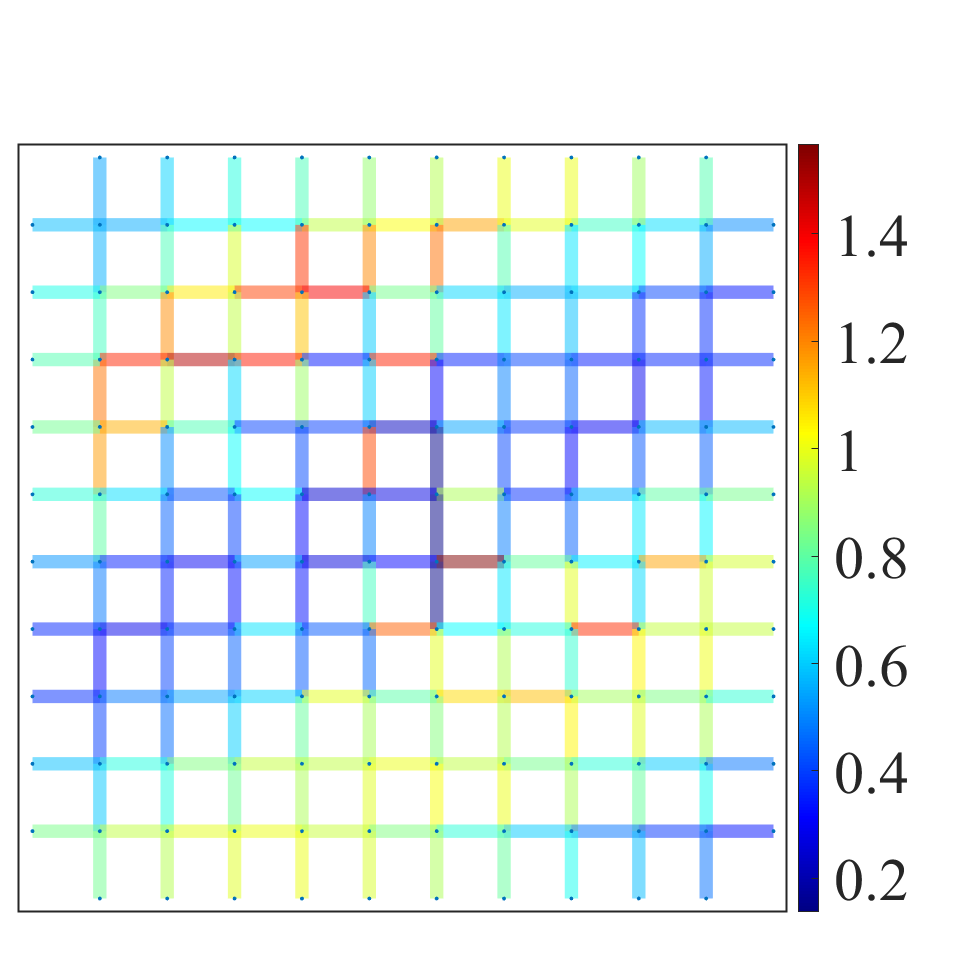} &\includegraphics[width=0.2\linewidth]{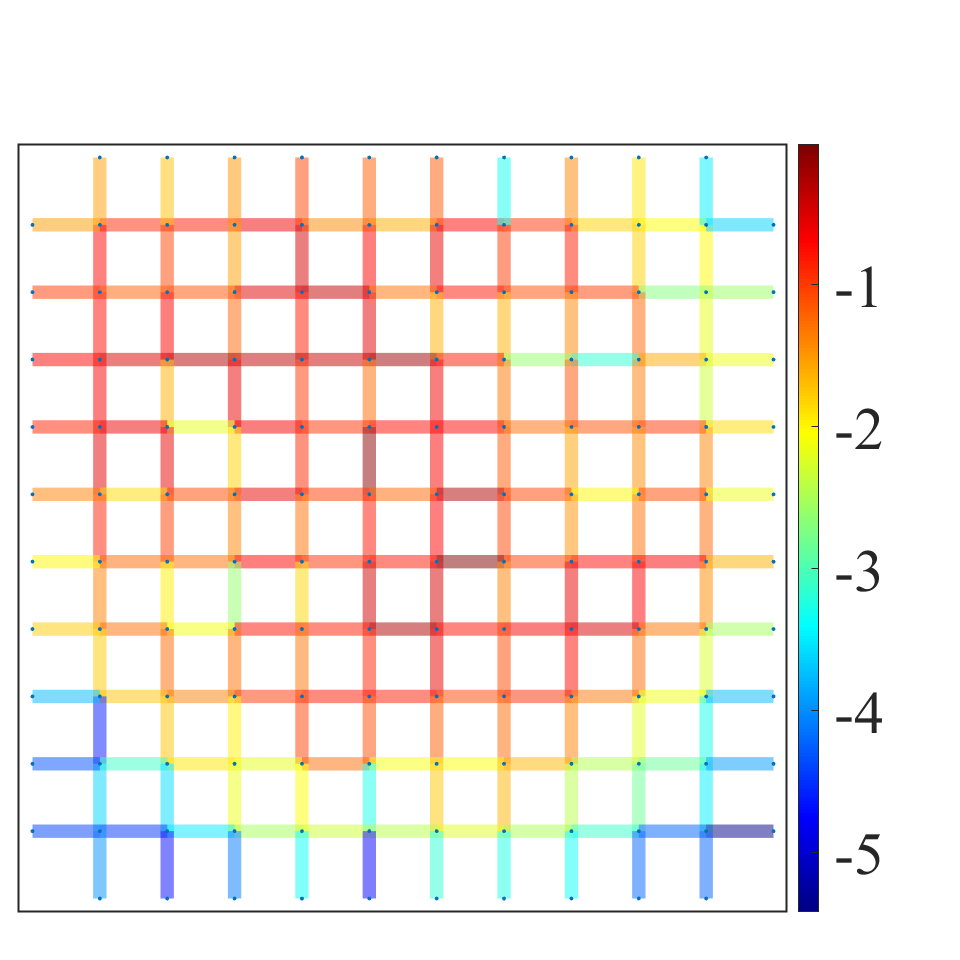}\\
    $\log_{10} |e_{\Lambda}|$ & $\widehat {\bm\gamma}$ & $\log_{10}|e_{\bm\gamma}|$ &  $\widehat {\bm\gamma}$ & $\log_{10}|e_{\bm\gamma}|$\\
\end{tabular}
\caption{Numerical results by the Curtis-Morrow algorithm and NN approach with partial data over the set $(1:2n)\times(n+1:4n)$ at three noise levels, $\epsilon=0\%$ (top), $\epsilon=0.001\%$ (middle) and $\epsilon=0.01\%$ (bottom).}\label{fig:partial-Rectangle}
\end{figure}

 \begin{remark}\label{rmk:interp21}
These results  illustrate the advantages of the  NN approach over the Curtis-Morrow algorithm for both full and partial data, in the presence of data noise. One possible explanation of its superior performance
lies in having selected the weights
$\widetilde W^{(1)}$ in the first layer independently
of $\widetilde  W^{(2)}$. This choice, encoded in the NN architecture,
seems to achieve a balance between traditional
reconstruction models, e.g., the Curtis-Morrow algorithm, which deal with noisy {\rm(}and hence out of range{\rm)} data by
projecting them on the range, and traditional
NN models which operate independently of the range of the conductivity $\bm\gamma$
to the DtN matrix $\Lambda_{\bm\gamma}$; the interpretable
NN operates like traditional methods if the data
is in the range, but enjoys the flexibility of
traditional NNs when
the DtN data is out of range. It seems that the additional flexibility
allows for better approximations in noisy and
partial data. Indeed, the Curtis-Morrow algorithm
keeps track of the dependence of the discrete
Green kernel from the variables
estimating the conductivities, and
forces the approximation to lie on the
appropriate nonlinear manifold, while the NN is not constrained to that.
\end{remark}
\section{Sensitivity analysis}\label{SENS}

Now we conduct a sensitivity analysis of the minimizer(s) to problem \eqref{loss}, and study the impact of noise in the data on the recovered conductivity $\widehat{{\bm\gamma}}\in \mathbb{R}^{2n(2n+1)}$.
The analysis is based on the observation that for the exact DtN data $\Lambda^\dag \equiv \Lambda_{{\bm\gamma}^\dag}$, upon ignoring the multiplicative constant $\frac1{2m}$ (which does not influence the minimizer), we can express the loss $C_\alpha(\widetilde W)$ equivalently as
\begin{equation} \label{eqn:opt}
C_{\alpha}(\widetilde  W) = \sum_{k=1}^{m}
\Big(\sum_{p\in \partial D} (   y_{p}(k))^2
+ \alpha \sum_{p\in  D} ( y_{p}(k))^2\Big).
\end{equation}

Numerically, the following constrained optimization problem in the limit case $\alpha=+\infty$ yields results that are consistent with the case $\alpha<\infty$:
\begin{align}\label{eqn:opt-constrain}    C_\infty(\widetilde W) %&=\sum_{k=1}^m \|\widetilde W^{(2)}(\widetilde W^{(1)}\mathbf  {\widetilde x}^{(1)}(k)) \|_2^2
&= \sum_{k=1}^{m}\sum_{p
\in \partial D} ( y_{p}(k))^2
 \\\quad &\mbox{with }\widetilde W\mbox{ s.t. } \sum_{p\in  D} (  y_{p}(k))^2=0,\quad k=1,\ldots,m.\nonumber
\end{align}
%where $\|\cdot\|_2$ denotes Euclidean norm of vectors.
Given the weight matrix \( \widetilde W^{(2)} \), the constraint
$\sum_{p \in D} (y_{p}(k))^2 =0$ for $k=1,\dots,m$
is equivalent to solving the first two equations in  \eqref{W2}. If we fix $\widetilde W^{(2)}$ and optimize with respect to $ \widetilde W^{(1)}$ using the input $
{\widetilde {\bf x}^{(1)}}(k)$, the input
$\widetilde W^{(1)}{\widetilde {\bf x}^{(1)}}(k)$
of the upper part of the second layer can be determined by \( \widetilde W^{(2)} \) and \(
{\widetilde{\bf  x}^{(1)}}(k) \)
and can achieve the zero loss. This is equivalent to solving for the potential \(\mathbf{u}=  \widetilde W^{(1)}{\widetilde {\bf x}^{(1)}}(k)  \) based on the conductivity \( {\bm\gamma} ={\bm \gamma}(\widetilde W^{(2)}) \) given  by
$\gamma_{\{q, p\}}=\widetilde w^{(2)}_{qp}$ for any $p,q \in D$ if $\{p,q\}\in B$, and the boundary potential \(
{\widetilde {\bf x}^{(1)}}(k) \) by solving the discrete Dirichlet problem, which has a unique solution. That is, we can determine \(\widetilde W^{(1)}{\widetilde {\bf x}^{(1)}}(k)  =: \mathbf{u}(\widetilde W^{(2)}, \mathbf{\widetilde x}^{(1)}(k)) \) based on \( \widetilde W^{(2)} \). Then the constrained loss
\( C_\infty(\widetilde W) \) can be reduced to
\[
C_\infty(\widetilde W^{(2)}) = \sum_{k=1}^m\| \widetilde W^{(2)}(\mathbf{u}(\widetilde W^{(2)}, {\widetilde {\bf x}^{(1)}}(k))) - \widehat{\mathbf{x}}^{(1)}(k)\|_2^2,
\]
where $\|\cdot\|_2$ denotes the Euclidean norm of vectors. Note that the reformulation does not change the global minimizer. With the Dirichlet data \(\{{\widetilde {\bf x}^{(1)}}(k)\}_{k=1}^m\), problem \eqref{eqn:opt-constrain} can be expressed equivalently as
\begin{equation} \label{eqn:opt-const-d}
\sum_{k=1}^{m} \|\left(\Lambda_{\bm\gamma} -\Lambda^\dag \right){\widetilde{\bf  x}^{(1)}}(k)\|^2_2,\quad\text{ with } {\bm\gamma}={\bm\gamma}(\widetilde W^{(2}).
\end{equation}

Now consider the Dirichlet data ${\bf \overline u}(k)=\mathbf e_{q^{(k)}}$ with $q^{(k)}$, given by  \eqref{order_boundary}.
Then the $k$-th column of the DtN matrix $\Lambda$
corresponds to the Neumann datum for ${\bf \overline u}(k)$. We can identify partial data
by a subset \(T' \subset \{1, 2, \ldots, 4n\}^2\). In addition, the entries may be perturbed by noise $\varepsilon$. For partial data over $T^\prime$, the corresponding loss $C_\infty$ is given in \eqref{eqn:argmin**_partial}.
The condition on $T^\prime$ is specified in Theorem  \ref{thm:sensitivity} below:  it must contain the index set $T$ (given in Fig. \ref{fig:network} (b)) to ensure uniqueness of the discrete inverse conductivity problem and also the full rank of the Jacobian to obtain the sensitivity in Theorem  \ref{thm:sensitivity} below.  Since there is a one-to-one correspondence between $\bm\gamma$ and $\widetilde W^{(2)}$, the sensitivity with respect to $\bm\gamma$ is equivalent to that with respect to $\widetilde W^{(2)}$. Then problem \eqref{eqn:opt-const-d} reads
\begin{equation}
    {\bm\gamma} _\varepsilon^{T^\prime}=\arg \min_{\bm\gamma} \sum_{(i,j)\in{T^\prime}} \left((\Lambda_{\bm\gamma})_{ij}-\Lambda_{ij}^\dag-\varepsilon_{ij}\right)^2.\label{eqn:argmin**_partial}
\end{equation}

Now we give some useful notation for the analysis of problem \eqref{eqn:argmin**_partial}.
Throughout, we identify the DtN matrix $\Lambda_{\bm \gamma}$ with a vector in $\mathbb{R}^{16n^2}$ for notational brevity. Note that $\nabla_{\bm\gamma}(\Lambda_{\bm\gamma})_{ij} \in \mathbb{R}^{2n(n+1)}$, and by stacking the vectors $\nabla_{\bm\gamma}(\Lambda_{\bm\gamma})_{ij}$ as columns for the indices $(i,j) \in \{1, \ldots, 4n\}^2$, we obtain the full Jacobian matrix $\nabla_{\bm\gamma}(\Lambda_{\bm\gamma}) \in \mathbb{R}^{16n^2 \times 2n(n+1)}$. The Jacobian $\nabla_{\bm\gamma}\Lambda_{\bm\gamma}$ can be efficiently computed.
The proof is given in Appendix \ref{Proof:Jacobian}.
\begin{proposition}\label{prop:Jacobian}
Let $u^{(i)}$ denote the solution to problem \eqref{EQ1} with Dirichlet boundary data ${\bf \overline u}(i)$. Then the Jacobian matrix $\nabla_{\bm \gamma} (\Lambda_{\bm \gamma})$ is given by
\begin{equation}\label{eqn:UUrepresent}
\nabla_{\bm\gamma}(\Lambda_{\bm\gamma})_{ij}=\Delta u^{(i)}\odot\Delta u^{(j)},
\end{equation}
where $\odot$ denotes Hadamard product, i.e., $ (\Delta u^{(i)}\odot\Delta u^{(j)})_\ell=\Delta _\ell u^{(i)}\Delta _\ell u^{(j)}$, and $\Delta_\ell u^{(i)}$ denotes taking the difference along the edge $\ell$ connecting nodes $p$ and $q$, i.e., $\Delta_\ell u^{(i)}= u^{(i)}_p-u^{(i)}_q$.  %\elena{maybe it is better to write what this difference is}.
\end{proposition}

The analysis employs the following result. Let $T\subseteq\{1, \ldots, 4n\}^2$ be of size $2n(n+1)$, i.e., the entries marked with $\ast$ in Fig. \ref{fig:network}(b). Note that the DtN matrices \(\Lambda_{\bm\gamma}\), associated with all admissible conductivities \({\bm\gamma} \in \mathcal{A}\) belong to a subset of \(\mathbb{R}^{16n^2}\). When we project these matrices into \(\mathbb{R}^{|T|}\), we obtain the set of exact partial data, denoted by \(\mathcal{P} = \{ w \in \mathbb{R}^{2n(n+1)} : w = (\Lambda_{\bm\gamma})_T, {\bm\gamma} \in \mathcal{A}\}\).

\begin{lemma}
There exists a neighborhood $\widetilde{\mathcal{P}}$ of $\mathcal{P}$ and a continuously differentiable map $\widetilde {R}: \widetilde{\mathcal{P}} \rightarrow \mathbb{R}^{16n^2}$ such that $\tilde{R}((\Lambda_\gamma)_T)=\Lambda_{\bm\gamma}$.
\end{lemma}
\begin{proof}
By \cite[Theorem 5.1]{curtis1991DNmap}, there exists a mapping \(R: \mathcal{P} \rightarrow \mathbb{R}^{16n^2}\) that allows recovering the full DtN matrix $\Lambda_{\bm\gamma}$ from its partial data, i.e., $R((\Lambda_{\bm\gamma})_T) = \Lambda_{\bm\gamma}$. The recovery process is iterative and proceeds as follows.
Suppose that we have recovered the entries indexed by $Q_k$, the algorithm then extends it into a larger set $Q_{k+1}:=Q_k\cup (i,j)$. The expanded entry $\Lambda_{ij} $ is recovered  by
\begin{equation*}
\Lambda_{ij} = A_{Q_k} B_{Q_k}^{-1} C_{Q_k},
\end{equation*}
where $A_{Q_k}$, $B_{Q_k}$ and $C_{Q_k}$ are some submatrices of $\Lambda_{\bm\gamma}|_{Q_k}$. Moreover, for exact data, the Curtis-Morrow theory \cite[Theorem 5.1]{curtis1991DNmap} ensures the invertibility of each matrix $B_{Q_k}$. Provided that $\Lambda_{ij}$, $(i,j)\in Q_k$, depends smoothly on the partial data $(\Lambda_\gamma)_T$ in a neighborhood $\tilde{P}_i$, then also  $\Lambda_{ij} $, $(i,j)\in Q_{k+1}$, depends smoothly on the partial data $(\Lambda_\gamma)_T$ in a neighborhood $\tilde{P}_{k+1}$. Since this process terminates in a finite number of steps, we get the desired map.
\end{proof}

Moreover, we define \(\widetilde{R}_{ij}: \widetilde{\mathcal{P}} \rightarrow \mathbb{R}\) to represent the \((i,j)\)th component of the map \(\widetilde{R}\). Let $S:=\nabla_{\bm\gamma}\Lambda_{\bm\gamma}\in\mathbb{R}^{16n^2\times2n(n+1)}$ be the Jacobian matrix, and denote by $S_Q$ the submatrix of $S$ formed by the rows indexed by the set $Q$.
The next result gives the expression of the Jacobian, i.e., the sensitivity.
\begin{theorem}\label{thm:sensitivity}
If $T\subset T'$, then the matrix $S_{T^\prime}$ is of full column rank, and $\nabla_\varepsilon {\bm\gamma}_\varepsilon^{T'} |_{\varepsilon = 0}= S_{T'}^\dag$.
\end{theorem}
\begin{proof}
By \cite[Theorem 4.2]{curtis1990determining}, the matrix \(S\)
%\elena{ \(S\) is never defined and so is \(S_{T'}\)}
has full column rank, with a rank $2n(n+1)$.
Now for any index $(i', j') \not\in T$, we have
%\elena{the following notation is not clear. Why are we using $s$ we should use $(i,j)$}
\[
\nabla_{\bm\gamma}(\Lambda_{\bm\gamma})_{i' j'} = \nabla_{\bm\gamma}\big(\tilde{R}_{i' j'}((\Lambda_{\bm\gamma})_T)\big) = \sum_{(i,j) \in T} (\nabla_{(i,j)} \tilde{R}_{i' j'})((\Lambda_{\bm\gamma})_T) \nabla_{\bm\gamma}(\Lambda_{\bm\gamma})_{ij}.
\]
Thus, the gradient of $\Lambda_{\bm\gamma}$ at the index \((i',j')\) is a linear combination of the derivatives at the indices in $T$. Consequently,
$\mathrm{rank}(S_T) = \mathrm{rank}(S) = 2n(n+1)$.
This shows the first assertion.  Next, by differentiating \eqref{eqn:argmin**_partial} with respect to ${\bm\gamma}$ and taking value at the minimum point ${\bm\gamma}_\varepsilon$, we obtain
\begin{equation}\label{eqn:kkt}
\sum_{(i,j)\in{T'}}\left((\Lambda_{{\bm\gamma}_\varepsilon})_{ij}-\Lambda_{ij}^\dag-\varepsilon_{ij}\right)\nabla_{\bm\gamma}(\Lambda_{\bm\gamma})_{ij}\big|_{{\bm\gamma}={\bm\gamma}_\varepsilon }=0,\quad \forall \varepsilon\in \mathbb R^{16n^2},
\end{equation}
where we use the subscript $ij$ to denote the entries of the vector $\varepsilon$ (and also $v$ below), similar to $(\Lambda_{\bm\gamma})_{ij}$. Now consider a directional noise $\varepsilon(s) = sv$ for $s \in \mathbb{R}$ and any fixed $v \in \mathbb{R}^{16n^2}$.  Let ${\bm\gamma}^\prime(s) = [\nabla_\varepsilon {\bm\gamma}_\varepsilon] v \big|_{\varepsilon = \varepsilon(s)}$, with $\nabla_\varepsilon {\bm\gamma}_\varepsilon\in \mathbb{R}^{2n(n+1)\times 16n^2}$.
Then differentiating \eqref{eqn:kkt} with respect to $s$ gives
\begin{align} &\sum_{(i,j)\in{T^\prime}}\bigg(\left((\Lambda_{{\bm\gamma}_\varepsilon})_{ij}-\Lambda_{ij}^\dag-\varepsilon_{ij}\right)[\nabla_{\bm\gamma}^2(\Lambda_{\bm\gamma})_{ij}]{{\bm\gamma}^\prime}\big|_{{\bm\gamma}={\bm\gamma}_\varepsilon }+(\nabla_{\bm\gamma}(\Lambda_{\bm\gamma})_{ij}\cdot {\bm\gamma}^\prime )\big|_{{\bm\gamma}={\bm\gamma}_\varepsilon }\nabla_{\bm\gamma}(\Lambda_{\bm\gamma})_{ij}\big|_{{\bm\gamma}={\bm\gamma}_\varepsilon }\nonumber\\
    &- v_{ij}\nabla_{\bm\gamma}(\Lambda_{\bm\gamma})_{ij}\big|_{{\bm\gamma}={\bm\gamma}_\varepsilon }\bigg)=0,\quad \text{ with } \varepsilon = \varepsilon(s),{\bm\gamma}^\prime ={\bm\gamma}^\prime(s),
\label{eqn:first_deri}
\end{align}
where $[\nabla_{\bm\gamma}^2(\Lambda_{\bm\gamma})_{ij}]\in\mathbb{R}^{2n(n+1)\times 2n(n+1)}$ is the Hessian matrix of the function $(\Lambda_{\bm\gamma})_{ij}$ with respect to ${\bm\gamma}$, and the notation $\cdot$ denotes taking the inner product between two vectors.
%\elena{thank you for clarifying. Still I think we need to change the term $v_{ij}\nabla_{\bm\gamma}(\Lambda_{\bm\gamma})_{ij}\big|_{{\bm\gamma}}$ in $v\cdot \nabla_{\bm\gamma}(\Lambda_{\bm\gamma})_{ij}\big|_{{\bm\gamma}}$ }
Since \(\Lambda_{{\bm\gamma}_\varepsilon}\big|_{\varepsilon=0}=\Lambda^\dag\), setting $s=0$ in \eqref{eqn:first_deri} yields that the first term vanishes. Hence,
\begin{align*}
&\sum_{(i,j)\in{T^\prime}}\big((\nabla_{\bm\gamma}(\Lambda_{\bm\gamma})_{ij}\cdot {\bm\gamma}^\prime)\big|_{{\bm\gamma}={\bm\gamma}_\varepsilon } \nabla_{\bm\gamma}(\Lambda_{\bm\gamma})_{ij}\big|_{{\bm\gamma}={\bm\gamma}_\varepsilon } - v_{ij}\nabla_{\bm\gamma}(\Lambda_{\bm\gamma})_{ij}\big|_{{\bm\gamma}={\bm\gamma}_\varepsilon }\big)=0,
\end{align*}
with $\varepsilon= \varepsilon(0)$ and ${\bm\gamma}^\prime={\bm\gamma}^\prime(0)$.
Equivalently, this reads
$S_{T^\prime}^\top S_{T^\prime}{\bm\gamma}^\prime-S_{T^\prime}^\top v=0$.
Since $S_{T'}$ is of full column rank, $ S_{T'}^\top S_{T'}$ is invertible. Hence, $ {\bm\gamma}^\prime = S_{T'}^\dag v$ and $\nabla_\varepsilon {\bm\gamma}_\varepsilon^{T'} = S_{T'}^\dag$.
\end{proof}

In a similar manner, we can derive higher-order sensitivity, by differentiating equation \eqref{eqn:first_deri} $k$ times ($k \geq 2$) with respect to $s$. Note that the coefficient of the highest-order derivative remains unchanged. Hence, the $k$th derivative $
{\bm\gamma}^{(k)} := \partial_s^k {\bm\gamma}_{\varepsilon(s)}|_{s=0}$
 satisfies
\[
{\bm\gamma}^{(k)} = S_{T'}^{\dag} H(v, {\bm\gamma}^{(1)}, \ldots, {\bm\gamma}^{(k-1)}),
\]
for some $H$ depending on $v$ and the low-order derivatives $\{{\bm\gamma}^{(j)}\}_{j=0}^{k-1}$.

Now we illustrate the sensitivity for (incomplete) DtN data. We denote by $C_{3n}$ (respectively $C_{2n}$) a subset containing entries in the last \(3n\) (respectively $ {2n}$) columns of the matrix $\Lambda_\varepsilon \in \mathbb{R}^{4n\times 4n}$. Similarly, the set $D_{2n}$ contains the first $ {2n}$ columns. Also we define $T_{\text{sq}}=(1:3 n)\times(n+1:4 n)$, and $T_{\text{re}}=(1:2n)\times(n+1:4 n)$.

\begin{remark}
Since the DtN matrix $\Lambda_{\bm\gamma}$ is symmetric, we can relax the hypothesis $T\subset T'$ of Theorem \ref{thm:sensitivity} {\rm(}also Theorem \ref{Thm1}{\rm)} to $T \subset T' \cup \overline{T'}$, where $\overline{T'} = \{(i,j): (j,i) \in T'\}$. Further, by the symmetry of the network on a square lattice under a $90^\circ$ rotation, the partial data $C_{3n}$, $C_{2n}$, $D_{2n}$, $T_{\mathrm{sq}}$ and $T_{\mathrm{re}}$ satisfy the condition of the theorem.
\end{remark}

Fig. \ref{fig:err-sensitivity} displays the pointwise error \( e_{\bm\gamma}^{T'} \) of the reconstruction by the NN approach (at a noise level $\delta=1\times 10^{-5}$) and the sensitivity \( {\bm\gamma}_{T'}'(0) \) computed using Theorem \ref{thm:sensitivity}. Numerically, it is observed that the shapes of the pointwise error and the sensitivity ${\bm\gamma}_{T'}'(0)$ match each other very well for all cases, except $T_{\rm re}$. Thus, the sensitivity analysis does provide useful information when the noise $\varepsilon$ is small. The quantitative results are presented in Table \ref{tab:Quantitative}, which indicates that
the two errors are matching nicely in terms of both magnitude and shape for all cases, except the case $T'=T_{\rm re}$.

% For the full data set, \( \|e_{\bm\gamma}\|_2 = 9.96 \times 10^{-4} \) and \( \|{\bm\gamma}'(0)\|_2 = 1 \times 10^{-3} \). For the partial data, \( \|e^{C_{3n}}_{\bm\gamma}\|_2 = 1.0617 \times 10^{-3} \) and \( \|{\bm\gamma}_{C_{3n}}'(0)\|_2 = 1.0605 \times 10^{-3} \); \( \|e^{D_{2n}}_{\bm\gamma}\|_2 = 0.51 \) and \( \|{\bm\gamma}_{D_{2n}}'(0)\|_2 = 3.05 \); \( \|e^{C_{2n}}_{\bm\gamma}\|_2 = 0.55 \) and \( \|{\bm\gamma}_{C_{2n}}'(0)\|_2 = 1.24 \).  The angles between \( e_{\bm\gamma} \) and \( {\bm\gamma}'(0) \) is \( 0.11 \) radians, and between \( e^{C_{3n}}_{\bm\gamma} \) and \( {\bm\gamma}_{C_{3n}}'(0) \) is \( 0.1646 \) radians.

\begin{figure}[hbt!]
\centering
\setlength{\tabcolsep}{0pt}
\begin{tabular}{cccc}
    \includegraphics[width=0.25\linewidth]{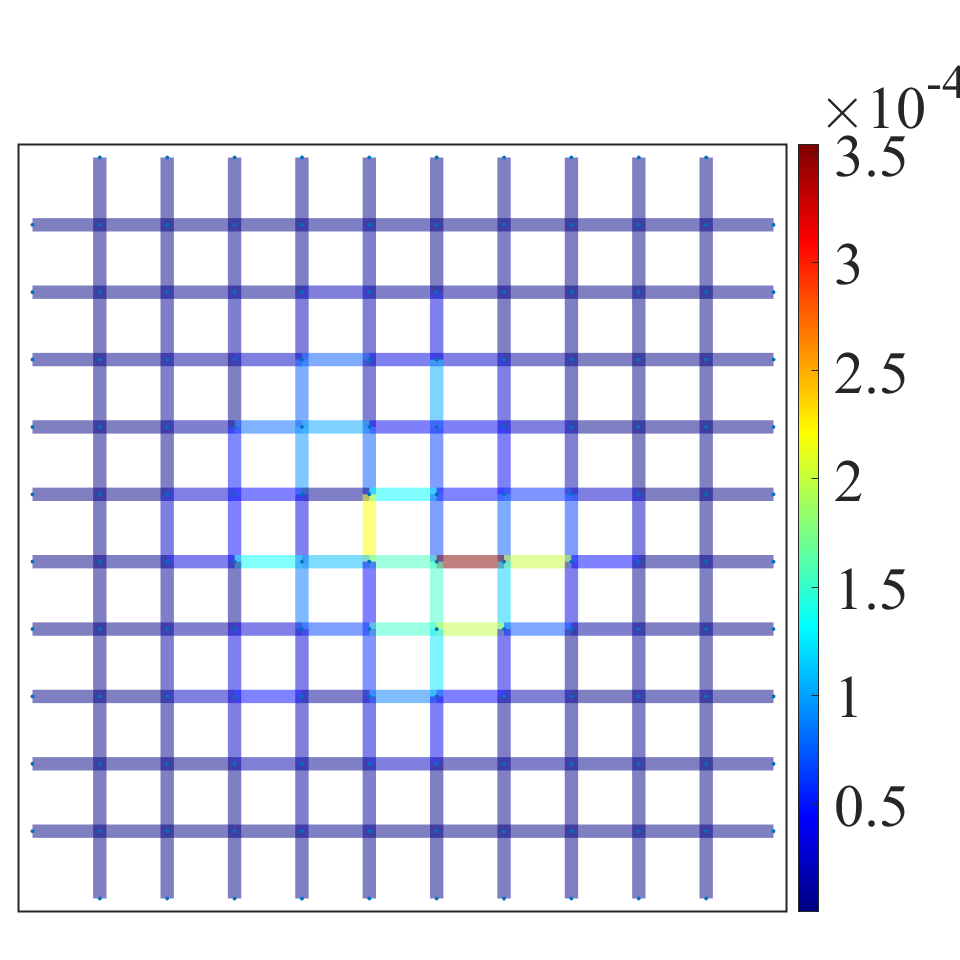} & \includegraphics[width=0.25\linewidth]{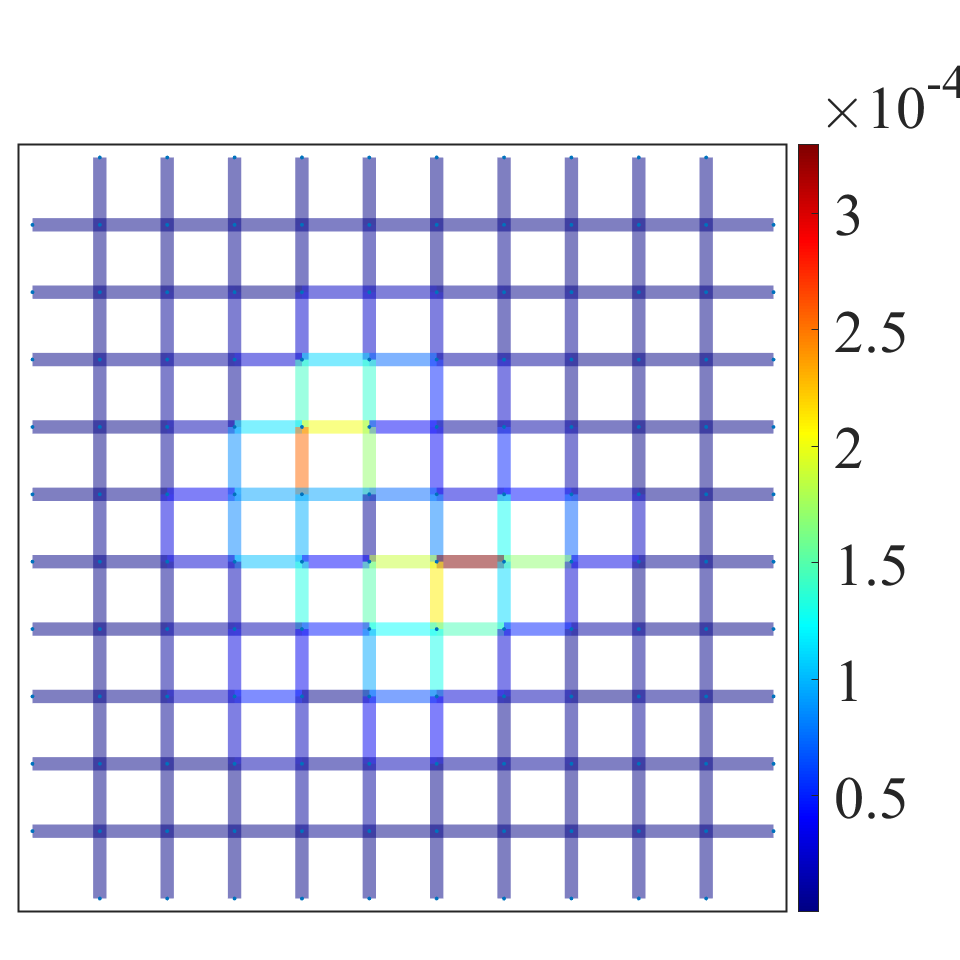} &\includegraphics[width=0.25\linewidth]{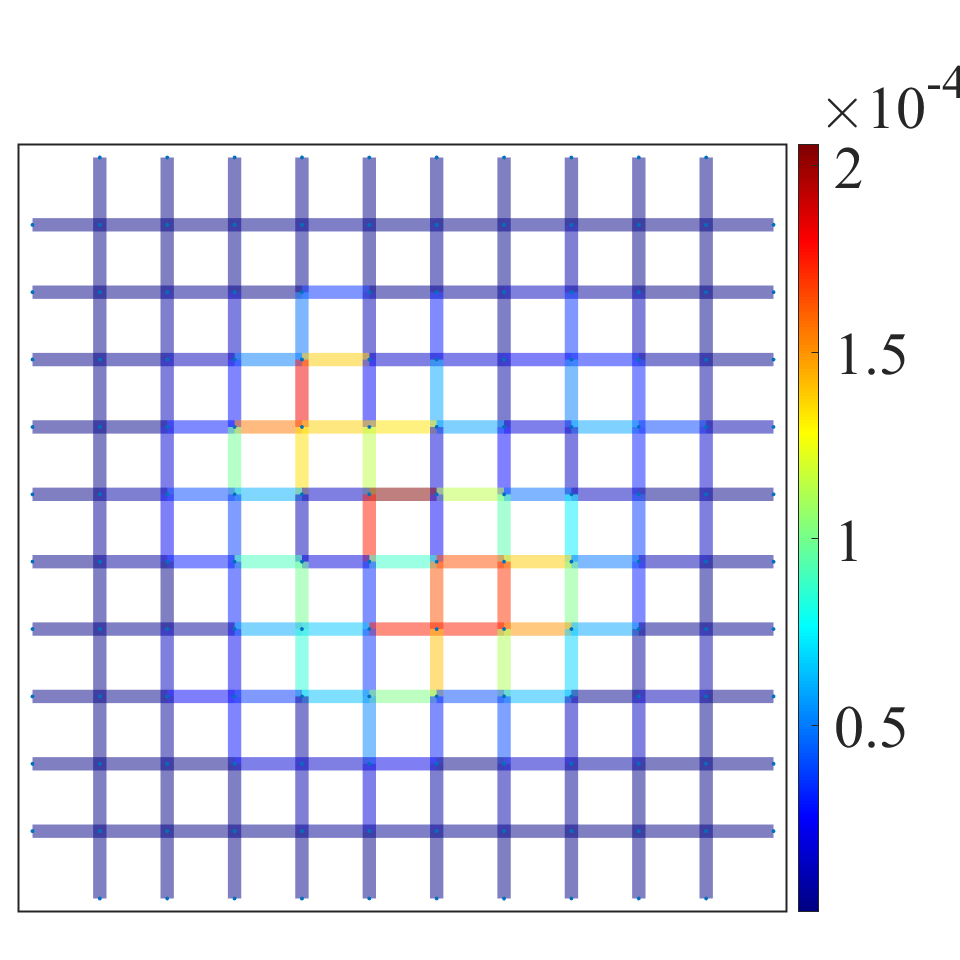} & \includegraphics[width=0.25\linewidth]{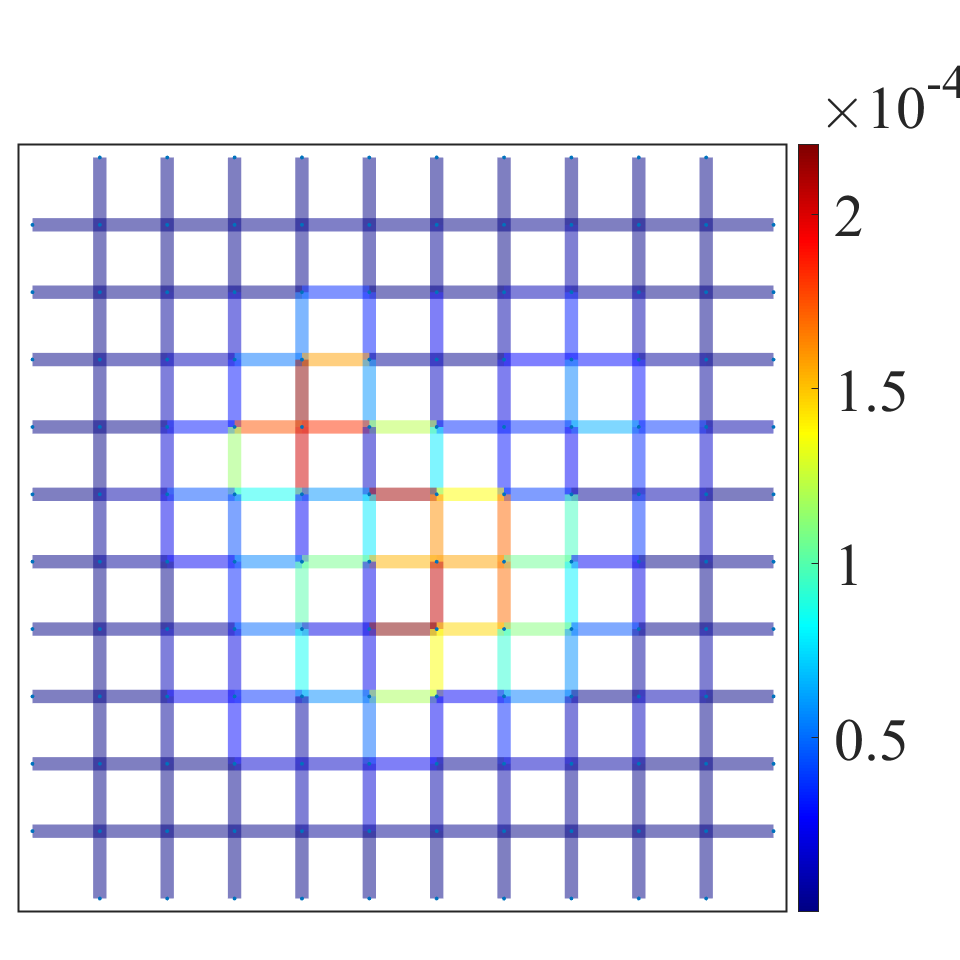}\\
   $ |e_{\bm\gamma}|$ & $ |{\bm\gamma}^\prime(0)|$ & $ |e^{C_{3n}}_{\bm\gamma}|$ & $ |{\bm\gamma}_{C_{3n}}^\prime(0)|$\\
 \includegraphics[width=0.25\linewidth]{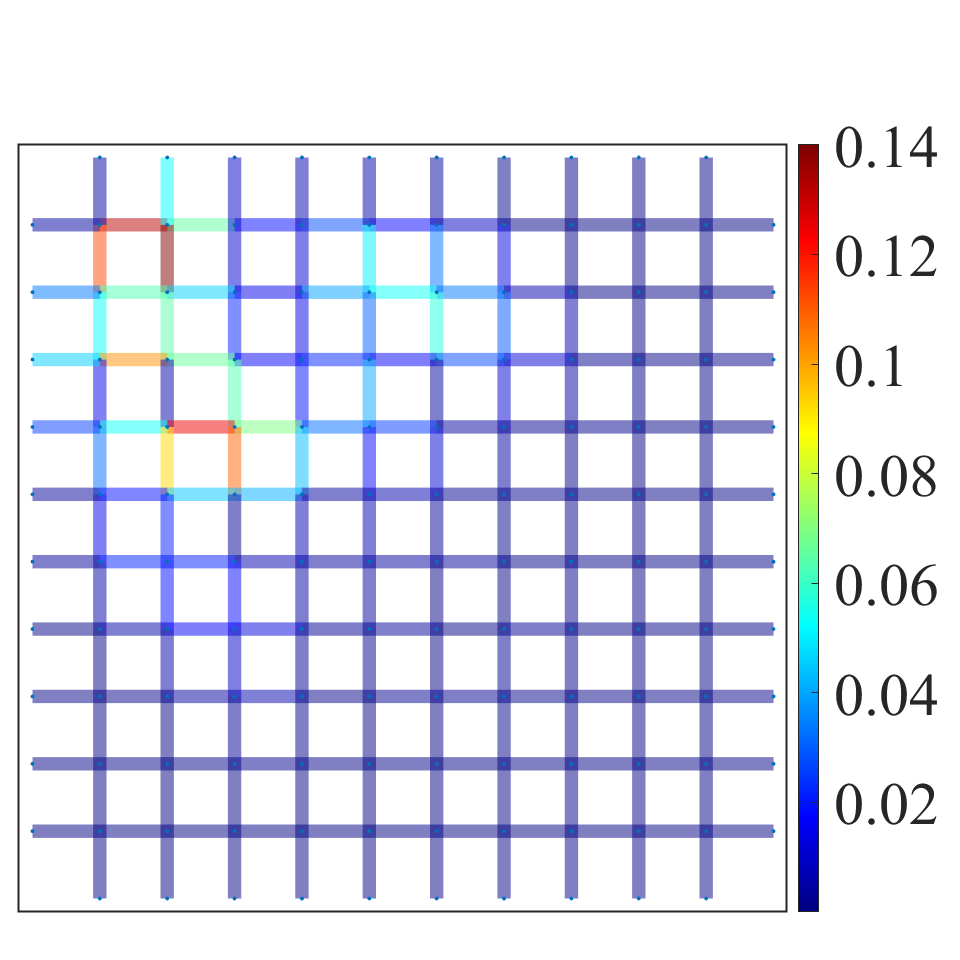} & \includegraphics[width=0.25\linewidth]{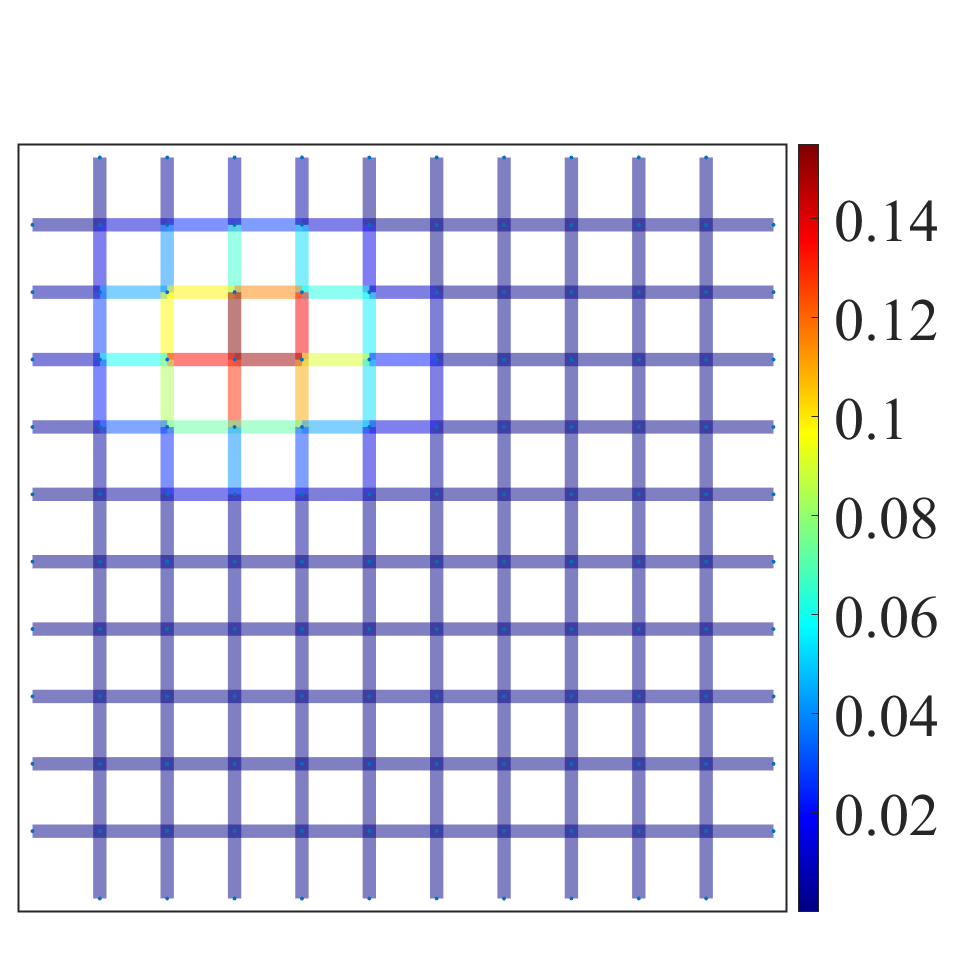} &\includegraphics[width=0.25\linewidth]{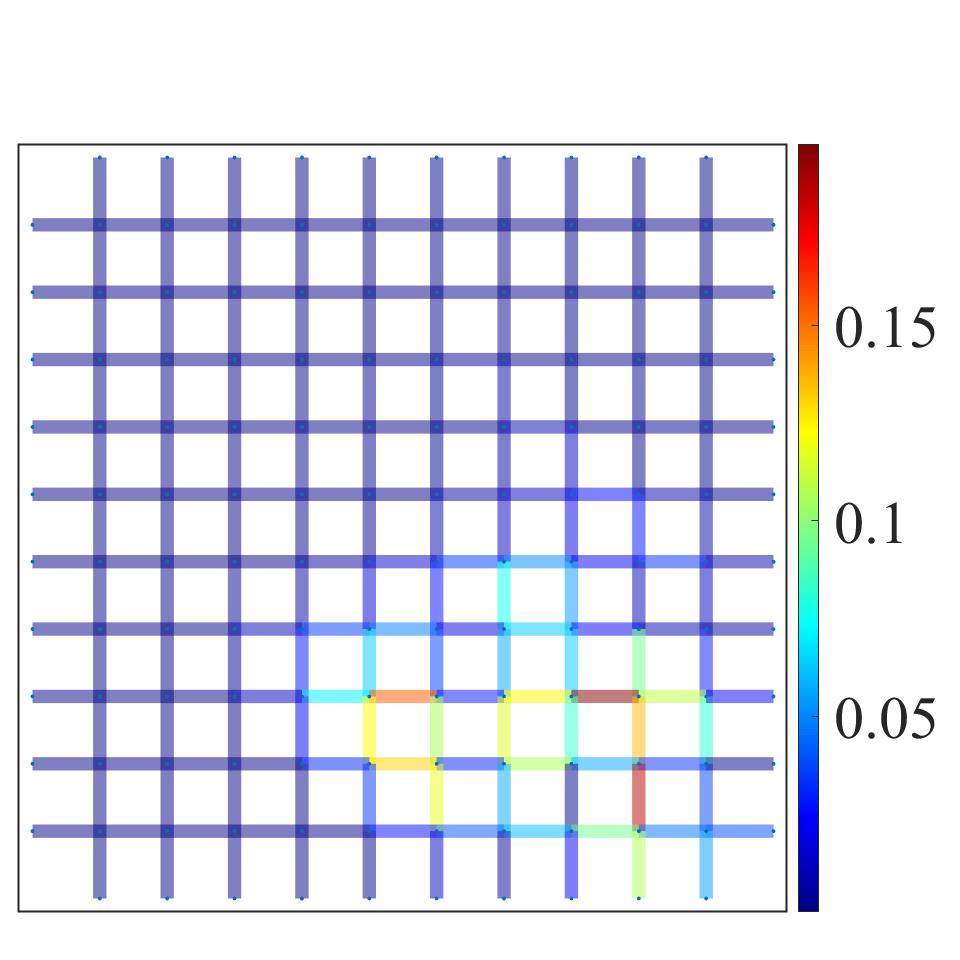} & \includegraphics[width=0.25\linewidth]{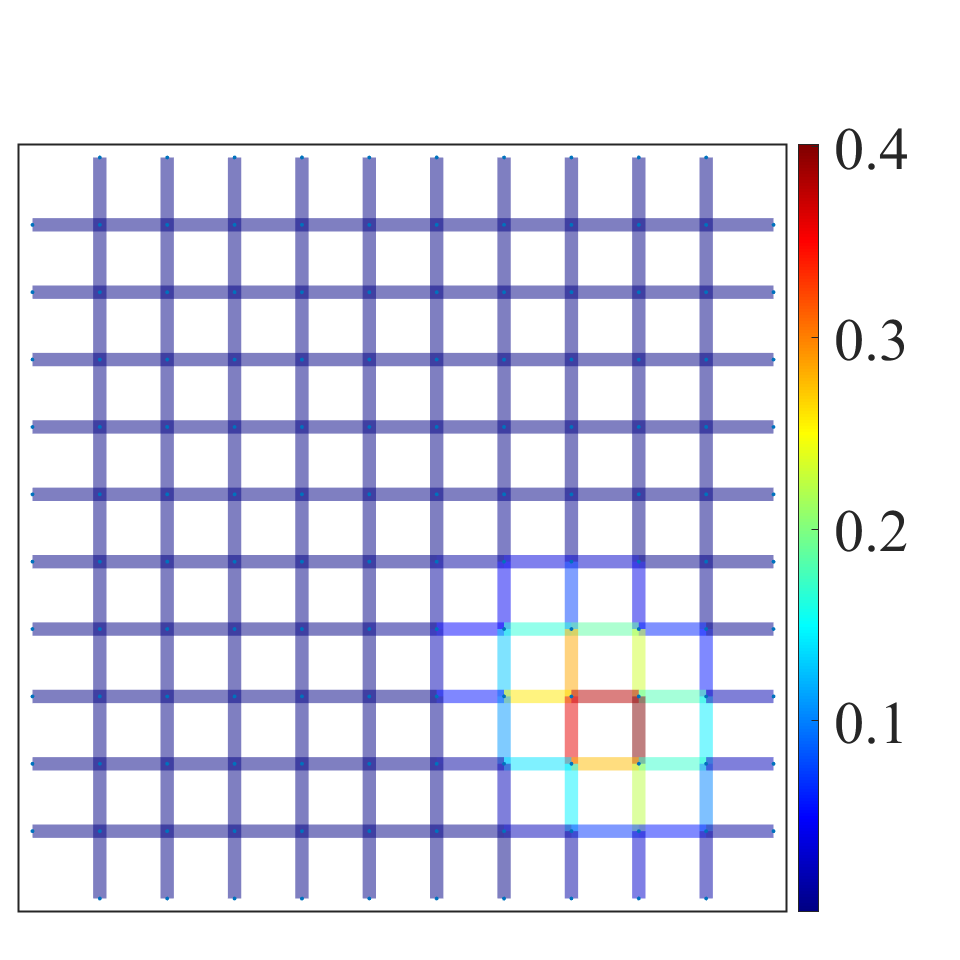}\\
   $ |e^{D_{2n}}_{\bm\gamma}|$ & $ |{\bm\gamma}_{D_{2n}}^\prime(0)|$ & $ |e^{C_{2n}}_{\bm\gamma}|$ & $ |{\bm\gamma}_{C_{2n}}^\prime(0)|$\\
 \includegraphics[width=0.25\linewidth]{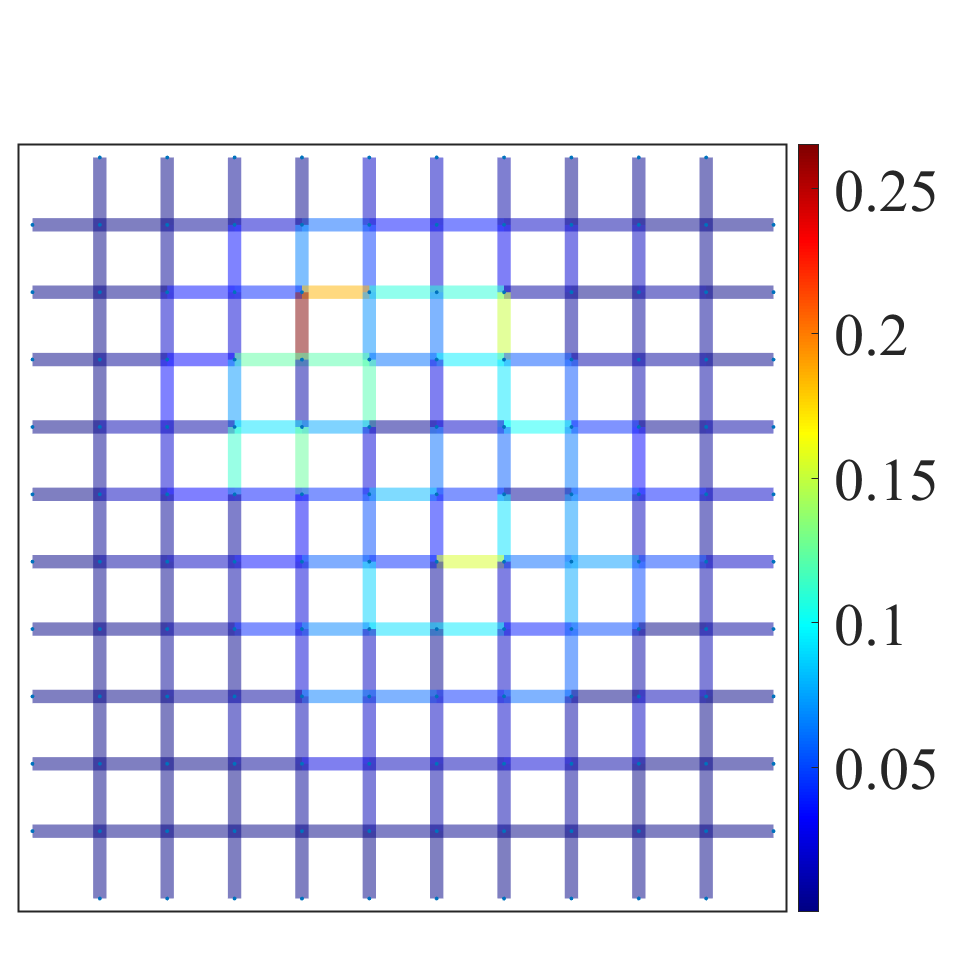} & \includegraphics[width=0.25\linewidth]{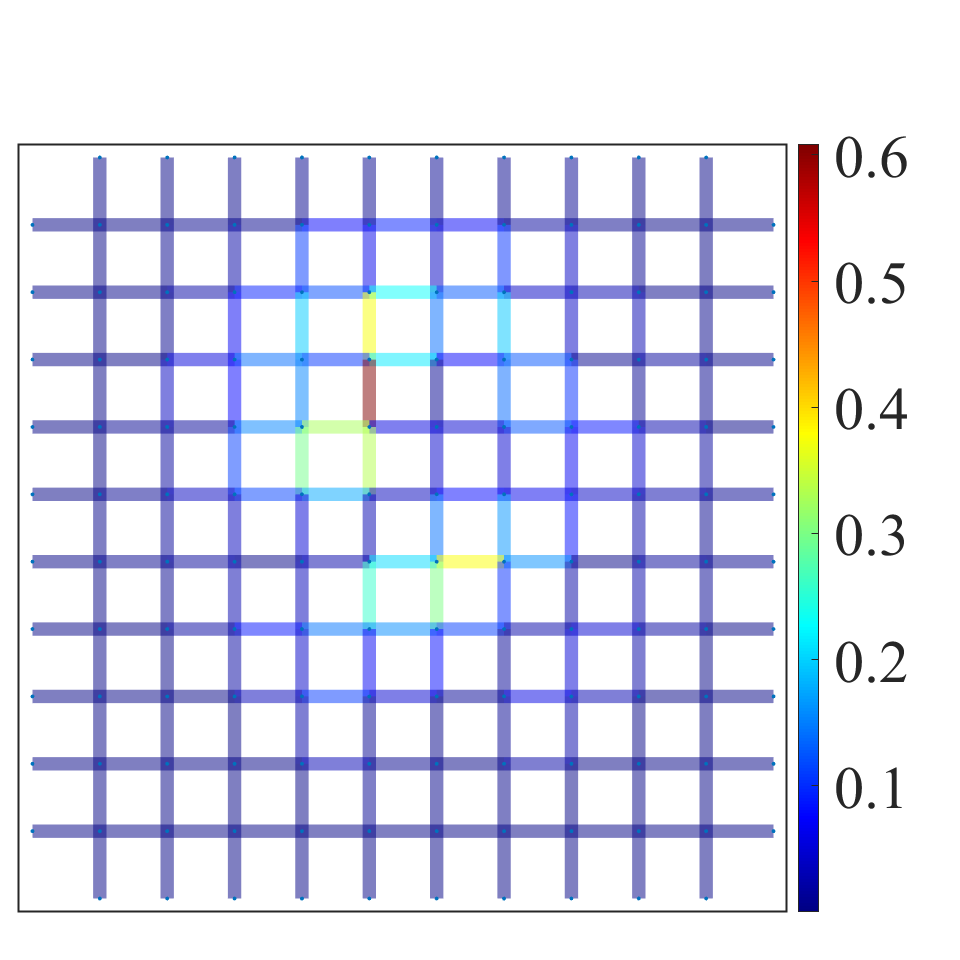} &\includegraphics[width=0.25\linewidth]{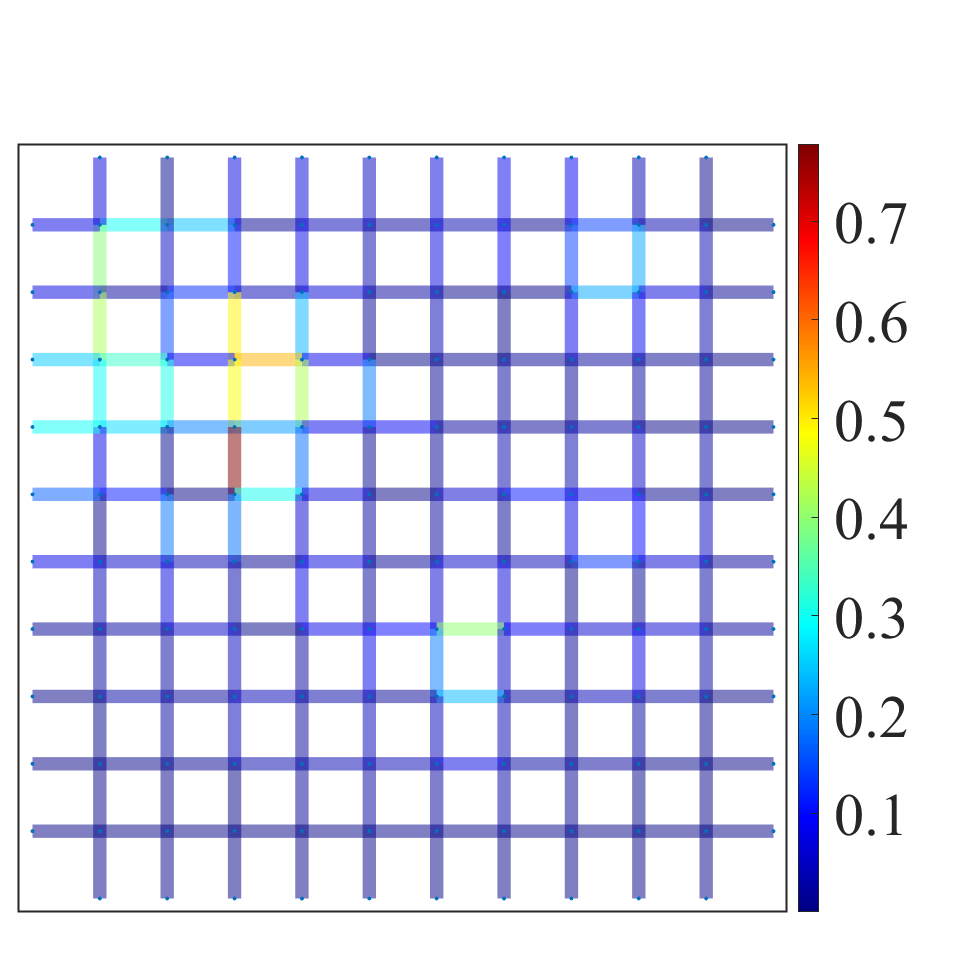} & \includegraphics[width=0.25\linewidth]{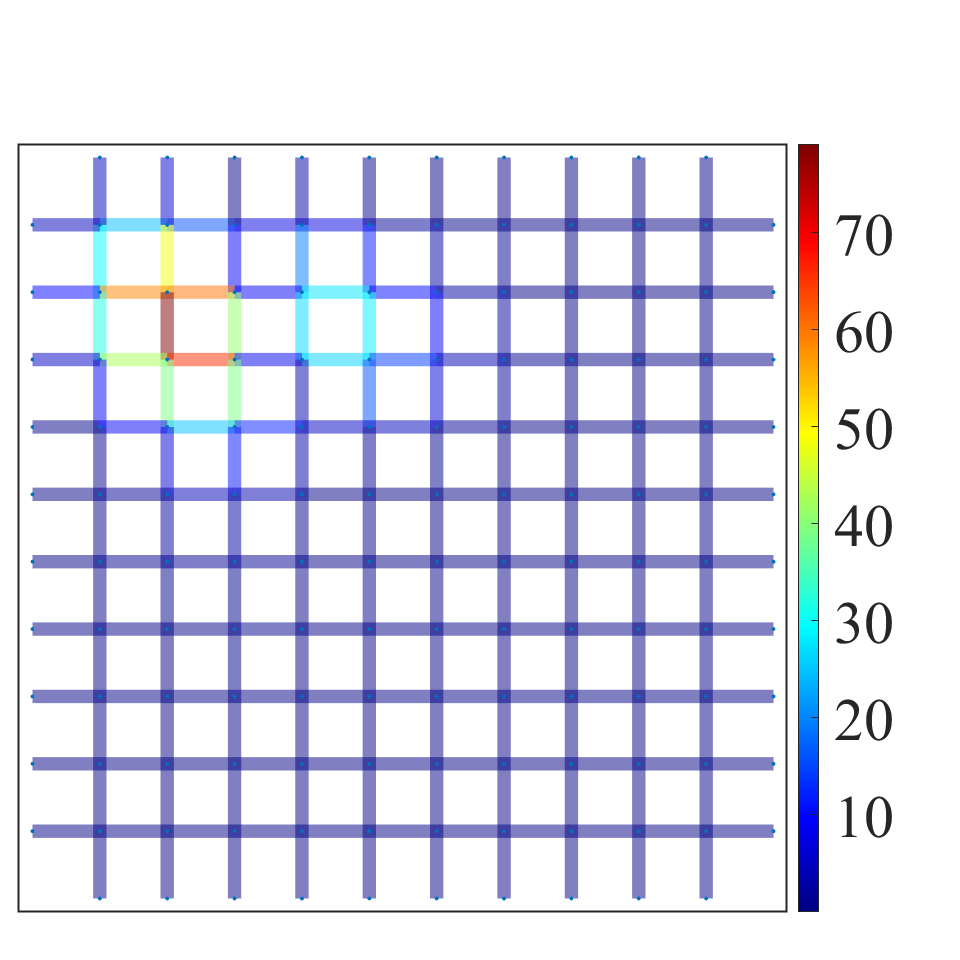}\\
   $ |e^{T_{\text{sq}}}_{\bm\gamma}|$ & $ |{\bm\gamma}_{T_{\text{sq}}}^\prime(0)|$ & $ |e^{T_{\text{re}}}_{\bm\gamma}|$& $ |{\bm\gamma}_{T_{\text{re}}}^\prime(0)|$\\
\end{tabular}
\caption{The comparison of error and first-order approximation with full data and partial data $C_{3n}$,  $D_{2n}$, $C_{2n}$, $T_{\text{sq}}$ and $T_{\text{re}}$. ($\alpha = \infty$, $\epsilon = 10^{-5}\%$)} \label{fig:err-sensitivity}
\end{figure}

\begin{table}[hbt!]
    \centering
    \begin{threeparttable}
    \caption{Quantitative results for the error and prediction by the sensitivity analysis for full, incomplete and partial data.}
    \label{tab:Quantitative}
    \begin{tabular}{c|c|c|c|c|c|c}
    \toprule
         Index $T^\prime$&  $\{1,\ldots,4n\}^2$&  $C_{3n}$&  $D_{2n}$& $C_{2n}$&$T_{\text{sq}}$&$T_{\text{re}}$\\
    \midrule  $\|e_{\bm\gamma}^{T^{\prime}}\|_2$&   $7.37 \times 10^{-4}$& $7.83 \times 10^{-3}$ &$0.41$  &$0.58$& $0.69$ &$1.93$\\
     \hline  $\|{\bm\gamma}_{T^\prime}'(0)\|_2$& $8.15\times10^{-4} $ &$8.47 \times 10^{-3}$  & $0.45$ &$1.24$&$1.38$&$187.45$\\
     \bottomrule
    \end{tabular}
    \end{threeparttable}
\end{table}

\section{Conclusions}

In this work we have proposed a novel interpretable neural network based approach for solving the discrete inverse conductivity problem on a square lattice. One distinct
feature of this approach is that the activation
functions are linear and the conductivity
is recovered by the weights of the second layer of the trained NN. The trained weights in
the first layer play the role of discrete Green's kernel, although  the nonlinear dependency on the
conductivities is not enforced during the actual
training. Further, we present a sensitivity analysis of the minimizer, which allows predicting the recovery error for small noise.

Numerical experiments indicate that the approach can give reasonable results for both full and incomplete / partial data, and outperforms the Curtis-Morrow algorithm for noisy data. The superior performance may lie in the role played by the weights
in the first layer. These weights can
approximate discrete Green's kernel
when the Dirichlet-to-Neumann data lies within the
range of the map $\bm\gamma\mapsto \Lambda_{\bm\gamma}$, and exploit the additional flexibility of not being constrained to adhere to discrete Green's kernel when the data is out of range.

This study leaves several interesting questions for further research. First, the proposed neural network employs the identity as the activation function, and it is natural to explore whether nonlinear ones may be profitably used. Second, the considered inverse problem enjoys Lipschitz stability due to its finite-dimensionality, but a precise quantification of the stability constant is still missing: the constant is expected to grow exponentially with the network size $n$, as indicated by the numerical experiments.
Third, it is of interest to explore the connection of circular resistor networks and optimal grids for EIT \cite{borcea2008electrical,BorceaDruskin:2010pyramid,BorceaMamonov:2010}, so as to obtain an effective neural network based solver for EIT: this naturally motivates extending the neural network approach to circular resistor networks and higher dimensions, including non-uniform optimal grids for improved reconstructions. Fourth and last, it is of interest to develop analogous approaches for similar inverse problems on graphs, e.g., optical tomography on graphs \cite{ChungGilbert:2017}.

\appendix
\section{Proof of Proposition \ref{prop:Jacobian}}\label{Proof:Jacobian}

The proof requires some notation. Given the conductivity $\bm\gamma$ on a square lattice (of size $n$), we define the linear map $L_{\bm\gamma}:\mathbb{R}^{n^2+4n}\rightarrow\mathbb{R}^{n^2+4n} $ by
\begin{equation*}
(L_{\bm\gamma} u)_p = \left\{ \begin{aligned}
 \sum_{q\sim p}\gamma_{pq}(u_p-u_q),   & \quad p\in D,\\
 u_p,   &\quad  p\in\partial D.
\end{aligned}\right.
\end{equation*}
We denote by $B_{\bm\gamma}:\mathbb{R}^{n^2+4n}\rightarrow\mathbb{R}^{ 4n}$ the Neumann operator
$(B_{\bm\gamma} u)_p =  \gamma_{pq}(u_p-u_q)$ for $p\in\partial D \text{ and } q= \mathcal  N(p)$. We can decompose the matrices $L_{\bm\gamma}$ and $B_{\bm\gamma}$ respectively into $$L_{\bm\gamma}=\begin{pmatrix}
 \Sigma & B\\
 0 & I_{4n}
\end{pmatrix}\quad \mbox{and}\quad B_{\bm\gamma} = \begin{pmatrix}
 B^\top & C
\end{pmatrix},$$
where $\Sigma\in \mathbb{R}^{n^2\times n^2}$ is symmetric and $C\in \mathbb{R}^{4n\times 4n}$ is diagonal. Then there hold
\begin{equation} \label{eqn:inverse}
L_{\bm\gamma}^{-1} = \begin{pmatrix}
 \Sigma^{-1}& -\Sigma^{-1}B\\
 0 & I_{4n}
\end{pmatrix} \quad \mbox{and}\quad  B_{\bm\gamma}L_{\bm\gamma}^{-1}=\begin{pmatrix}
 B^\top \Sigma^{-1} & -B^\top \Sigma^{-1}B+C
\end{pmatrix}.
\end{equation}
Also we define the extension operator $\Pi:\mathbb{R}^{ 4n}\rightarrow\mathbb{R}^{ n^2+4n}$ and restriction operator $P_r:\mathbb{R}^{n^2+4n}\to \mathbb{R}^{n^2+4n}$ by
$ \Pi = \begin{pmatrix}
        0\\ I_{4n}
    \end{pmatrix}$ and $P_r = \begin{pmatrix}
        I_{n^2} & 0 \\
        0 & 0
    \end{pmatrix}.$
The DtN map $\Lambda_{\bm\gamma}$ is then given by
\begin{equation*} \Lambda_{\bm \gamma} =- B_{\bm\gamma} L_{\bm\gamma}^{-1} \Pi.
\end{equation*}
Now fix any edge  $\ell$. The product rule and the identity ${\partial}_{\bm\gamma_\ell} L_{\bm\gamma}^{-1} = L_{\bm\gamma}^{-1}({\partial}_{\gamma_\ell} L_{\bm\gamma}  ) L_{\bm\gamma}^{-1}$ yield
\begin{align*}
   - {\partial}_{\gamma_\ell}\Lambda_{\bm\gamma}
    &=({\partial}_{\gamma_\ell} B_{\bm\gamma}) L_{\bm\gamma}^{-1} \Pi+ B_{\bm\gamma}  L_{\bm\gamma}^{-1}({\partial}_{\gamma_\ell} L_{\bm\gamma}  ) L_{\bm\gamma}^{-1}\Pi.
\end{align*}
Next we discuss the two cases $\ell$, i.e., an interior edge and a boundary edge, separately.\\
Case 1: The edge $\ell$ connects two interior nodes $p,q\in D$. Then we have ${\partial}_{\gamma_\ell} B_{\bm\gamma}=0$, and
\begin{align*}
  -  {\partial}_{\gamma_\ell}\Lambda_{\bm\gamma}& = B_{\bm\gamma}  L_{\bm\gamma}^{-1}({\partial}_{\gamma_\ell} L_{\bm\gamma}  ) L_{\bm\gamma}^{-1}\Pi.
\end{align*}
Note that $(({\partial}_{\gamma_\ell} L_{\bm\gamma}) u)_p = u_p-u_q $, $(({\partial}_{\gamma_\ell} L_{\bm\gamma}) u)_q = u_q-u_p$ and $(({\partial}_{\gamma_\ell} L_{\bm\gamma}) u)_r=0$ for $r\ne p,q$.
%Let $P_r:\mathbb{R}^{ n^2+4n}\rightarrow\mathbb{R}^{ n^2+4n}$ be the restriction operator given by $ [I^\top \ \bar{\mathbf{u}^\top}]^\top \mapsto [I^\top\ \mathbf{0}_{1\times 4n}]^\top$.
Direct computation with the identities in \eqref{eqn:inverse} gives
\begin{align*}
    B_{\bm\gamma} L_{\bm\gamma}^{-1}P_r=\begin{pmatrix}
 B^\top \Sigma^{-1} &0
\end{pmatrix} =-\Pi^\top L_{\bm\gamma}^{-\top} P_r.
\end{align*}
From the relation $P_r {\partial}_{\gamma_\ell} L_{\bm\gamma}   ={\partial}_{\gamma_\ell} L_{\bm\gamma} $, it follows that
\begin{align*}
  &-  {\partial}_{\gamma_\ell}\Lambda_{\bm\gamma}=  B_{\bm\gamma}  L_{\bm\gamma}^{-1}P_r({\partial}_{\gamma_\ell} L_{\bm\gamma}  ) L_{\bm\gamma}^{-1}\Pi\\
  =& -\Pi^\top L_{\bm\gamma}^{-\top}P_r({\partial}_{\gamma_\ell} L_{\bm\gamma}  ) L_{\bm\gamma}^{-1}\Pi
  = -\Pi^\top L_{\bm\gamma}^{-\top} ({\partial}_{\gamma_\ell} L_{\bm\gamma}  ) L_{\bm\gamma}^{-1}\Pi.
\end{align*}
Now by the definition $u^{(i)}= L_{\bm\gamma}^{-1}\Pi e_i$, we get $$
e_i^\top\partial_{\gamma_\ell}\Lambda_{\bm\gamma} e_j=(u^{(i)}(p)-u^{(i)}(q))(u^{(j)}(p)-u^{(j)}(q))=\Delta_\ell u^{(i)}\Delta_\ell u^{(j)}.$$
Case 2: The edge $\ell$ connects the boundary node $p$ and the interior node $q=\mathcal N(p)$. Then
$$ -  {\partial}_{\gamma_\ell}\Lambda_{\bm\gamma}=({\partial}_{\gamma_\ell} B_{\bm\gamma}) L_{\bm\gamma}^{-1} \Pi -\Pi^\top L_{\bm\gamma}^{-\top} ({\partial}_{\gamma_\ell} L_{\bm\gamma}  ) L_{\bm\gamma}^{-1}\Pi,$$
and $(({\partial}_{\gamma_\ell} L_{\bm\gamma})u  )_q =u_q-u_p $ and $(({\partial}_{\gamma_\ell} L_{\bm\gamma})u  )_r=0$ for $r\ne q$; $(({\partial}_{\gamma_\ell} B_{\bm\gamma})u  )_p =u_p-u_q  $ and $(({\partial}_{\gamma_\ell} B_{\bm\gamma})u  )_r = 0$ for $r\ne p$. Hence
$$e_i^\top\partial_{\gamma_\ell}\Lambda_{\bm\gamma} e_j= -\delta_{ip}(u^{(j)}(q)-u^{(j)}(p))+u^{(i)}(q)(u^{(j)}(q)-u^{(j)}(p)).$$
Since $u^{(i)}(p) = \delta_{ip}$, we have $$
e_i^\top\partial_{\gamma_\ell}\Lambda_{\bm\gamma} e_j=(u^{(i)}(p)-u^{(i)}(q))(u^{(j)}(p)-u^{(j)}(q))=\Delta_\ell u^{(i)}\Delta_\ell u^{(j)}.$$
Combining these two cases completes the proof of the proposition.

\bibliographystyle{abbrv}
\bibliography{references}

\begin{thebibliography}{10}

\bibitem{akhtar2023survey}
N.~Akhtar.
\newblock A survey of explainable ai in deep visual modeling: Methods and
  metrics.
\newblock Preprint, arXiv:2301.13445, 2023.

\bibitem{al2001feed}
A.~F. Al-Mudhaf.
\newblock {\em A feed forward neural network approach for matrix computations}.
\newblock PhD thesis, Brunel University, London, 2001.
\newblock Available at \url{https://bura.brunel.ac.uk/handle/2438/5010}.

\bibitem{Aless1988}
G.~Alessandrini.
\newblock Stable determination of conductivity by boundary measurements.
\newblock {\em Appl. Anal.}, 27(1-3):153--172, 1988.

\bibitem{A1989}
G.~Alessandrini.
\newblock Determining conductivity by boundary measurements, the stability
  issue.
\newblock In {\em Applied and Industrial Mathematics ({V}enice, 1989)},
  volume~56 of {\em Math. Appl.}, pages 317--324. Kluwer Acad. Publ.,
  Dordrecht, 1991.

\bibitem{AV}
G.~Alessandrini and S.~Vessella.
\newblock Lipschitz stability for the inverse conductivity problem.
\newblock {\em Adv. Appl. Math.}, 35(2):207--241, 2005.

\bibitem{alvarez2018towards}
D.~Alvarez~Melis and T.~Jaakkola.
\newblock Towards robust interpretability with self-explaining neural networks.
\newblock In {\em Advances in Neural Information Processing Systems},
  volume~31, pages 7786--7795, 2018.

\bibitem{BF}
E.~Beretta and E.~Francini.
\newblock Lipschitz stability for the electrical impedance tomography problem:
  the complex case.
\newblock {\em Comm. Partial Differential Equations}, 36(10):1723--1749, 2011.

\bibitem{BerMicPerSan18}
E.~Beretta, S.~Micheletti, S.~Perotto, and M.~Santacesaria.
\newblock Reconstruction of a piecewise constant conductivity on a polygonal
  partition via shape optimization in {EIT}.
\newblock {\em J. Comput. Phys.}, 353:264--280, 2018.

\bibitem{bishop1994neural}
C.~M. Bishop.
\newblock Neural networks and their applications.
\newblock {\em Rev. Sci. Instrum.}, 65(6):1803--1832, 1994.

\bibitem{BlastenIsozaki:2023}
E.~Bl{\aa}sten, H.~Isozaki, M.~Lassas, and J.~Lu.
\newblock Inverse problems for discrete heat equations and random walks for a
  class of graphs.
\newblock {\em SIAM J. Discrete Math.}, 37(2):831--863, 2023.

\bibitem{blazek2021explainable}
P.~J. Blazek and M.~M. Lin.
\newblock Explainable neural networks that simulate reasoning.
\newblock {\em Nature Comput. Sci.}, 1(9):607--618, 2021.

\bibitem{Borcea}
L.~Borcea.
\newblock Electrical impedance tomography.
\newblock {\em Inverse Problems}, 18(6):R99--R136, 2002.

\bibitem{BorceaMamonov:2010}
L.~Borcea, V.~Druskin, and A.~V. Mamonov.
\newblock Circular resistor networks for electrical impedance tomography with
  partial boundary measurements.
\newblock {\em Inverse Problems}, 26(4):045010, 30, 2010.

\bibitem{BorceaDruskin:2010pyramid}
L.~Borcea, V.~Druskin, A.~V. Mamonov, and F.~Guevara~Vasquez.
\newblock Pyramidal resistor networks for electrical impedance tomography with
  partial boundary measurements.
\newblock {\em Inverse Problems}, 26(10):105009, 36, 2010.

\bibitem{borcea2008electrical}
L.~Borcea, V.~Druskin, and F.~G. Vasquez.
\newblock Electrical impedance tomography with resistor networks.
\newblock {\em Inverse Problems}, 24(3):035013, 2008.

\bibitem{Borcea:2013}
L.~Borcea, F.~Guevara~Vasquez, and A.~V. Mamonov.
\newblock Study of noise effects in electrical impedance tomography with
  resistor networks.
\newblock {\em Inverse Probl. Imaging}, 7(2):417--443, 2013.

\bibitem{BoyerGarzella:2016}
J.~Boyer, J.~J. Garzella, and F.~Guevara~Vasquez.
\newblock On the solvability of the discrete conductivity and {S}chr\"odinger
  inverse problems.
\newblock {\em SIAM J. Appl. Math.}, 76(3):1053--1075, 2016.

\bibitem{CenJin:2023}
S.~Cen, B.~Jin, K.~Shin, and Z.~Zhou.
\newblock Electrical impedance tomography with deep {C}alder\'{o}n method.
\newblock {\em J. Comput. Phys.}, 493:112427, 2023.

\bibitem{chakraborty2024explainable}
M.~Chakraborty.
\newblock Explainable neural networks: achieving interpretability in neural
  models.
\newblock {\em Arch. Comput. Methods Eng.}, 31:3535--3550, 2024.

\bibitem{ChungGilbert:2017}
F.~J. Chung, A.~C. Gilbert, J.~G. Hoskins, and J.~C. Schotland.
\newblock Optical tomography on graphs.
\newblock {\em Inverse Problems}, 33(5):055016, 21, 2017.

\bibitem{CurtisMorrow:1994}
E.~Curtis, E.~Mooers, and J.~Morrow.
\newblock Finding the conductors in circular networks from boundary
  measurements.
\newblock {\em RAIRO Mod\'el. Math. Anal. Num\'er.}, 28(7):781--814, 1994.

\bibitem{curtis1990determining}
E.~B. Curtis and J.~A. Morrow.
\newblock Determining the resistors in a network.
\newblock {\em SIAM J. Appl. Math.}, 50(3):918--930, 1990.

\bibitem{curtis1991DNmap}
E.~B. Curtis and J.~A. Morrow.
\newblock The {Dirichlet} to {Neumann} map for a resistor network.
\newblock {\em SIAM J. Appl. Math.}, 51(4):1011--1029, 1991.

\bibitem{Fan2020}
Y.~Fan and L.~Ying.
\newblock Solving electrical impedance tomography with deep learning.
\newblock {\em J. Comput. Phys.}, 404:109--119, 2020.

\bibitem{GaoGuan:2023}
L.~Gao and L.~Guan.
\newblock Interpretability of machine learning: recent advances and future
  prospects.
\newblock {\em IEEE MultiMedia}, 30(4):105--118, 2023.

\bibitem{GardeKnudsen:2017}
H.~Garde and K.~Knudsen.
\newblock Distinguishability revisited: depth dependent bounds on
  reconstruction quality in electrical impedance tomography.
\newblock {\em SIAM J. Appl. Math.}, 77(2):697--720, 2017.

\bibitem{GuoJiang:2021}
R.~Guo and J.~Jiang.
\newblock Construct deep neural networks based on direct sampling methods for
  solving electrical impedance tomography.
\newblock {\em SIAM J. Sci. Comput.}, 43(3):B678--B711, 2021.

\bibitem{HamiltonHauptmann:2018}
S.~J. Hamilton and A.~Hauptmann.
\newblock Deep d-bar: real-time electrical impedance tomography imaging with
  deep neural networks.
\newblock {\em IEEE Trans. Med. Imag.}, 37(10):2367--2377, 2018.

\bibitem{han2015learning}
S.~Han, J.~Pool, J.~Tran, and W.~Dally.
\newblock Learning both weights and connections for efficient neural network.
\newblock In {\em Advances in Neural Information Processing Systems},
  volume~28, pages 1135--1143, 2015.

\bibitem{Ingerman:2010}
D.~V. Ingerman.
\newblock Discrete and continuous {D}irichlet-to-{N}eumann maps in the layered
  case.
\newblock {\em SIAM J. Math. Anal.}, 31(6):1214--1234, 2000.

\bibitem{JinLiQuanZhou:2024}
B.~Jin, X.~Li, Q.~Quan, and Z.~Zhou.
\newblock Conductivity imaging from internal measurements with mixed
  least-squares deep neural networks.
\newblock {\em SIAM J. Imaging Sci.}, 17(1):147--187, 2024.

\bibitem{Kingma2014AdamAM}
D.~P. Kingma and J.~Ba.
\newblock Adam: A method for stochastic optimization.
\newblock In {\em 3rd Inter- national Conference for Learning Representations},
  San Diego, 2015.

\bibitem{LamPylyavskyy:2012}
T.~Lam and P.~Pylyavskyy.
\newblock Inverse problem in cylindrical electrical networks.
\newblock {\em SIAM J. Appl. Math.}, 72(3):767--788, 2012.

\bibitem{li2022efficient}
L.~Li and J.~Hu.
\newblock An efficient second-order neural network model for computing the
  {Moore--Penrose} inverse of matrices.
\newblock {\em IET Signal Proc.}, 16(9):1106--1117, 2022.

\bibitem{li2019novel}
X.~Li, Y.~Zhou, J.~Wang, Q.~Wang, Y.~Lu, X.~Duan, Y.~Sun, J.~Zhang, and Z.~Liu.
\newblock A novel deep neural network method for electrical impedance
  tomography.
\newblock {\em Trans. Inst. Measure. Control}, 41(14):4035--4049, 2019.

\bibitem{M}
N.~Mandache.
\newblock Exponential instability in an inverse problem for the
  {S}chr\"{o}dinger equation.
\newblock {\em Inverse Problems}, 17(5):1435--1444, 2001.

\bibitem{ronneberger2015u}
O.~Ronneberger, P.~Fischer, and T.~Brox.
\newblock U-net: Convolutional networks for biomedical image segmentation.
\newblock In {\em Proceedings of the Eighteenth Conference on Medical Image
  Computing and Computer-Assisted Intervention}, pages 234--241. MICCAI, 2015.

\bibitem{rudin2019stop}
C.~Rudin.
\newblock Stop explaining black box machine learning models for high stakes
  decisions and use interpretable models instead.
\newblock {\em Nature Mach. Intell.}, 1(5):206--215, 2019.

\bibitem{rudin2022interpretable}
C.~Rudin, C.~Chen, Z.~Chen, H.~Huang, L.~Semenova, and C.~Zhong.
\newblock Interpretable machine learning: Fundamental principles and 10 grand
  challenges.
\newblock {\em Stat. Surv.}, 16:1--85, 2022.

\bibitem{sejnowski2020unreasonable}
T.~J. Sejnowski.
\newblock The unreasonable effectiveness of deep learning in artificial
  intelligence.
\newblock {\em Proc. Nat. Acad. Sci.}, 117(48):30033--30038, 2020.

\bibitem{seo2019learning}
J.~K. Seo, K.~C. Kim, A.~Jargal, K.~Lee, and B.~Harrach.
\newblock A learning-based method for solving ill-posed nonlinear inverse
  problems: A simulation study of lung {EIT}.
\newblock {\em SIAM J. Imaging Sci.}, 12(3):1275--1295, 2019.

\bibitem{TanyuMaass:2023}
D.~N. Tanyu, J.~Ning, A.~Hauptmann, B.~Jin, and P.~Maass.
\newblock Electrical impedance tomography: A fair comparative study on deep
  learning and analytic-based approaches.
\newblock In T.~Bubba, editor, {\em Data-Driven Models in Inverse Problems},
  pages 437--470. de Gruyter, Berlin, 2024.

\bibitem{wang1990structured}
L.~Wang and J.~M. Mendel.
\newblock Structured trainable networks for matrix algebra.
\newblock In {\em International Joint Conference on Neural Networks}, pages
  125--132, San Diego, CA, 1990.

\bibitem{yuan2006model}
M.~Yuan and Y.~Lin.
\newblock Model selection and estimation in regression with grouped variables.
\newblock {\em J. R. Stat. Soc. Ser. B Stat. Method.}, 68(1):49--67, 2006.

\bibitem{Zhang:2021}
Y.~Zhang, P.~Ti\~{n}o, A.~Leonardis, and K.~Tang.
\newblock A survey on neural network interpretability.
\newblock {\em IEEE Trans. Emerg. Topics Comput. Intell.}, 5(5):726--742, 2021.

\end{thebibliography}

\end{document}